\documentclass[11pt]{article}

\usepackage[margin=1in]{geometry}

\usepackage{natbib}

\usepackage[utf8]{inputenc}
\usepackage{amsmath, amssymb, amsthm}

\usepackage{hyperref}
\usepackage{xcolor} 

\hypersetup{
    colorlinks=true,
    linkcolor=red,
    citecolor=blue,
    filecolor=blue,
    urlcolor=magenta,
    linktocpage=true
}

\usepackage{subcaption}
\usepackage{graphicx}
\graphicspath{{figures/}}
\usepackage{float}

\usepackage[linesnumbered,ruled,vlined]{algorithm2e}

\SetCommentSty{mycommfont}
\SetKwInput{KwInput}{Input}                
\SetKwInput{KwOutput}{Output}              

\usepackage{mathrsfs}
\usepackage[scr=esstix]{mathalfa}

\newcommand{\spa}[1]{^{(#1)}}

\newcommand{\penabr}{FCLS}

\newcommand{\opnorm}[3]{||#1||_{(#2 \to #3)}} 
\newcommand{\opone}[1]{||#1||_{\text{op}-1} } 
\newcommand{\opn}[1]{||#1||_{\text{op}} } 

\newcommand{\ind}[1]{\mathbf{1}\left(#1 \right)}

\newcommand{\MatSuppSym}{\mathcal{S}}
\newcommand{\matsupp}[1]{\MatSuppSym\left(#1 \right)}
\newcommand{\blocksupp}[1]{\mathcal{BS}\left(#1 \right)}

\newcommand{\prob}[1]{\mathcal{P}\left(#1\right)} 
\newcommand{\expect}[1]{\mathbb{E}\left[#1\right]} 

\newcommand{\fclsSym}{s}
\newcommand{\fclsFunc}[1]{\fclsSym_{#1}}
\newcommand{\fcls}[2]{\fclsSym_{#1}\left(#2 \right)}

\newcommand{\lapSpecFuncSym}{t}
\newcommand{\lapSpecFunc}[1]{\lapSpecFuncSym_{#1}}
\newcommand{\lapSpec}[2]{\lapSpecFuncSym_{#1}\left(#2 \right)}

\newcommand{\lossSym}{\ell}
\newcommand{\loss}[1]{\lossSym\left(#1 \right)}

\newcommand{\est}{\beta}
\newcommand{\estorc}{\widehat{\est}^{\text{oracle}}}
\newcommand{\estinit}{\widehat{\est}^{\text{initial}}}
\newcommand{\estzero}{\widehat{\est}^{(0)}} 
\newcommand{\estok}{\widehat{\est}^{\text{ok}}} 

\newcommand{\estOrcEntry}{\estorc_{\text{entrywise}}} 
\newcommand{\minNnzMagTarg}{||\esttarg_{\nnzSetTarg}||_{\text{min}}} 

\newcommand{\estOrcLasso}{\widehat{\est}^{\text{orc, lasso}}}

\newcommand{\estCurrent}{b} 

\newcommand{\estone}{\widehat{\est}^{(1)}}

\newcommand{\esttwo}{\widehat{\est}^{(2)}}

\newcommand{\estLasso}{\widehat{\est}_{\text{lasso}}}
\newcommand{\estLassoSubscr}[1]{\widehat{\est}_{\text{lasso}, #1}}

\newcommand{\esttarg}{\est^{*}}

\newcommand{\suppSet}{\mathcal{S}} 
\newcommand{\suppSize}{s} 
\newcommand{\suppSetTarg}{\suppSet^*} 
\newcommand{\suppSetTargPow}[1]{\suppSet^{*^{#1}}} 
\newcommand{\suppSizeTarg}{\suppSize^*} 

\newcommand{\binGraphTarg}{\mathcal{G}^*}

\newcommand{\nnzSize}{\mathscr{s}} 
\newcommand{\nnzSizeTarg}{\nnzSize^*} 
\newcommand{\nnzSet}{\mathbb{S}} 
\newcommand{\nnzSetTarg}{\nnzSet^*} 
\newcommand{\nnzSetTargPow}[1]{\nnzSet^{*^{#1}}} 

\newcommand{\maxDegTarg}{\mathscr{d}^*}
\newcommand{\maxDegTargPow}[1]{\mathscr{d}^{*^{#1}}}

\newcommand{\supp}{\mathcal{S}} 
\newcommand{\suppTarg}{\mathcal{S}^*} 

\newcommand{\ncc}{K_{\text{cc}}}
\newcommand{\nccTarg}{\ncc^*}
\newcommand{\nccTargPow}[1]{\ncc^{*^{#1}}}

\newcommand{\nccNz}{K_{\text{cc, nz}}}
\newcommand{\nccNzTarg}{\nccNz^*}

\newcommand{\ccSym}{\mathcal{C}}
\newcommand{\cc}[1]{\ccSym_{#1}}
\newcommand{\ccTarg}[1]{\ccSym^*_{#1}}

\newcommand{\maxccTarg}{d^*_{\text{max}}}
\newcommand{\maxccTargPow}[1]{d^{*^{#1}}_{\text{max}}}

\newcommand{\minccTarg}{d^*_{\text{min}}}

\newcommand{\NumNonIso}{d_{\text{non-iso}}}

\newcommand{\NumIso}{\mathscr{i}^*}

\newcommand{\numMissingEdges}{E_{\text{missing}} }
\newcommand{\graphDiam}{\text{diam}}

\newcommand{\maxMissingTarg}{E^*_{\text{max missing}} }

\newcommand{\maxDiamTarg}{\text{diam}^*_{\text{max}} }

\newcommand{\supptarg}{s^*}

\newcommand{\supptargpow}[1]{s^{*^{#1}}}

\newcommand{\tailclass}[2]{\mathcal{T}(#1, #2 )}
\newcommand{\tailconst}{v_*}

\newcommand{\lapcoeff}[2]{\mathcal{M}(#1, #2)}

\newcommand{\residDiffGapRatio}{\rho}

\newcommand{\initRatioConsts}{C_{\residDiffGapRatio} }

\newcommand{\tuneparam}{\tau}

\newcommand{\ThreshTuneParam}{\gamma}

\newcommand{\ParmMax}{\xi}
\newcommand{\lassoOrcSet}[1]{\mathcal{B}_{#1}^{\text{lasso orc}}}
\newcommand{\PrevlassoOrcSet}[1]{\mathcal{B}_{#1}^{\text{orc}}}   

\newcommand{\lassoOrcSetOp}[1]{\mathcal{B}_{#1}^{\text{lasso orc, op}}} 
\newcommand{\lassoOrcSetFrob}[1]{\mathcal{B}_{#1}^{\text{lasso orc, frob}}} 

 \newcommand{\lapSmEval}[2]{\lambda_{(#1)} \left(\mathcal{L}(#2)\right)}

\newcommand{\RestrEvalLinReg}{ \kappa_{\text{linear}}}

\newcommand{\setNiceGrad}[1]{\mathcal{B}^{\text{nice grad}}_{#1}} 
\newcommand{\setBigEvalComp}[1]{\mathcal{B}^{\text{big eval, component-wise}}_{#1}} 
\newcommand{\setEvalGap}[1]{\mathcal{B}^{\text{gap}}_{#1}} 

\newcommand{\probGoodInitLLA}{\delta^{\text{good init}}}
\newcommand{\probGoodInitLLAOp}{\delta^{\text{good init, op}}}
\newcommand{\probGoodInitLLAFrob}{\delta^{\text{good init, frob}}}

\newcommand{\probNiceGradLLA}{\delta^{\text{nice grad}}}

\newcommand{\probSmallResidLLA}{\delta^{\text{small lap resid}}}

\newcommand{\probNiceGradLLAOp}{\delta^{\text{nice grad, op}}}
\newcommand{\probSmallResidLLAOp}{\delta^{\text{small lap resid, op}}}

\newcommand{\probNiceGradLLAFrob}{\delta^{\text{nice grad, frob}}}
\newcommand{\probSmallResidLLAFrob}{\delta^{\text{small lap resid, frob}}}

\newcommand{\probNiceGradStatPt}{\probNiceGradLLA_{\text{orc}}}

\newcommand{\probSmallResidStatPt}{\probSmallResidLLA_{\text{orc}}}

\newcommand{\targgap}{\triangle^*}

\newcommand{\minSuppMagTarg }{||\esttarg_{\suppSetTarg}||_{\text{min}}}

\newcommand{\qTarg}{q^*} 
\newcommand{\qTargPow}[1]{q^{*^{#1}}}

\newcommand{\maxColNorm}{M}
\newcommand{\maxevalxexpl}{ \lambda_{\text{max}} \left(\frac{1}{n}  X^T X \right)}
\newcommand{\minevalxsexpl}{ \lambda_{\text{min}} \left(\frac{1}{n}  X_{\suppSetTarg}^T X_{\suppSetTarg} \right)}
\newcommand{\minevalxsshort}{\lambda_{\text{min}, \suppSetTarg} }
\newcommand{\maxevalxshort}{\lambda_{\text{max}}}

 \newcommand{\LogisitcLink}[1]{\psi\left(#1\right)} 
\newcommand{\LogisitcLinkGrad}[1]{\psi'\left(#1\right)} 
\newcommand{\LogisitcLinkGG}[1]{\psi''\left(#1\right)} 
\newcommand{\LogisitcLinkGGG}[1]{\psi'''\left(#1\right)} 

\newcommand{\LogisticGradCond}[1]{G_{\suppSetTarg}\left(#1\right)} 
 
\newcommand{\LogisitcMu}[1]{\mu\left(#1\right)}
\newcommand{\LogisitcMuS}[1]{\mu_{\suppSetTarg}\left(#1\right)}

\newcommand{\LogisitcH}[1]{h\left(#1\right)}
\newcommand{\LogisitcHS}[1]{h_{\suppSetTarg}\left(#1\right)}

\newcommand{\LogisitcT}[1]{t\left(#1\right)}
\newcommand{\LogisitcTS}[1]{t_{\suppSetTarg}\left(#1\right)}

\newcommand{\LogisticEvalMax}{Q_1} 
 
\newcommand{\LogisticXsHXscOne}{Q_2} 

\newcommand{\LogisticXshXsInv}{Q_3}   

\newcommand{\LogisticMaxVarNorm}{\maxColNorm} 

\newcommand{\LogisticMaxX}{m} 

\newcommand{\LogisticRestrEval}{ \kappa_{\text{logistic}}} 

\newcommand{\LogisticTempResidConst}{A}

\newcommand{\tempMatFrank}{\Psi}
 
\newcommand{\threshop}[2]{T_{#2}\left( #1\right)}

\newcommand{\ProofSection}[1]{\noindent \textbf{\underline{\text{#1}}}}

\newcommand{\sameAsPrev}{ \textendash\textendash\textendash \| \textendash\textendash\textendash}

\newcommand{\subGResidExpectBound}{\Xi}

\newcommand{\estHT}[1]{\widehat{\est}^{\text{HT}}_{#1}}

\newcommand{\probHTBigTrue}{ \delta_{\text{big true}}^{\text{HT}} }
\newcommand{\probHTSmallNoise}{ \delta_{\text{small noise}}^{\text{HT}} }

\newcommand{\opoEBoundAbs}{\mu_{op-1, abs}} 
\newcommand{\opoEBound}{\mu_{op-1}}

\newcommand{\LassoKillerBound}{\ThreshTuneParam_{\text{lasso-killer-lbd}}}

\newcommand{\FCLSKillerBound}{\tuneparam_{\text{FCLS-killer-lbd}}}

\newcommand{\HyperGraph}{\mathcal{H}}
\newcommand{\VetSet}{\mathcal{V}}
\newcommand{\EdgeSet}{\mathcal{E}}
\newcommand{\EdgeWeights}{\mathcal{W}}
\newcommand{\HyperAdjMat}{\mathcal{HA}}

\newcommand{\DimSumMA}{s}
\newcommand{\DimProdMA}{p}

\newcommand{\frobBound}{\mathbb{F}}
\newcommand{\opBound}{\mathbb{P}}

\newcommand{\arbNorm}[1]{\left| \left | \left | #1 \right | \right | \right |}

\newtheorem{proposition}{Proposition}[section]
\newtheorem{definition}{Definition}[section]
\newtheorem{corollary}{Corollary}[section]
\newtheorem{theorem}{Theorem}[section]
\newtheorem{remark}{Remark}[section]
\newtheorem{lemma}{Lemma}[section]

\newtheorem{assumption}{Assumption}[section]
\newtheorem{fact}{Fact}[section]

\title{The folded concave Laplacian spectral penalty learns block diagonal sparsity patterns with the strong oracle property}

\author{Iain Carmichael\footnote{idc9@uw.edu}}
\date{\today}

\begin{document}

\maketitle

\begin{abstract}
Structured sparsity is an important part of the modern statistical toolkit.
We say a set of model parameters has \textit{block diagonal sparsity up to permutations} if its elements can be viewed as the edges of a graph that has multiple connected components.
For example, a block diagonal correlation matrix with $K$ blocks of variables corresponds to a graph with $K$ connected components whose nodes are the variables and whose edges are the correlations.
This type of sparsity captures clusters of model parameters.
To learn block diagonal sparsity patterns we develop the \textit{folded concave Laplacian spectral penalty} and provide a majorization-minimization algorithm for the resulting non-convex problem.
We show this algorithm has the appealing computational and statistical guarantee of converging  to the oracle estimator after two steps with high probability, even in high-dimensional settings. 
The theory is then demonstrated in several classical problems including covariance estimation, linear regression, and logistic regression.
\end{abstract}


\tableofcontents

\section{Introduction} \label{s:intro}

Sparsity plays a crucial role in statistics and machine learning \citep{hastie2015statistical}.
Estimators that incorporate sparsity help the data analyst obtain interpretable results and the theoretician establish favorable estimation properties.
Some problems have additional structure beyond \textit{entrywise sparsity}.
For example, in the \textit{group sparse} setting the analyst knows that predefined groups of variables should be included or excluded together \citep{yuan2006model}.
There is a growing literature on \textit{structured sparsity} inducing penalties including: the group lasso \citep{yuan2006model}, fused lasso/total variation/edge lasso \citep{tibshirani2005sparsity, sharpnack2012sparsistency}, graph total variation \citep{li2020graph}, and generalized Lasso \citep{tibshirani2011solution}.
Many of these established sparsity inducing penalties -- and their concave extensions -- come with appealing statistical and computational guarantees \citep{negahban2012unified, fan2020statistical}. 
This paper adds a new structured sparsity inducing penalty to this list that aims to discover groupings of variables.

We say a vector, $\est$, has \textit{block diagonal sparsity up to permutations} if $\est$ can be viewed as the edges of a graph and this graph has multiple connected components.
This is made precise below.
For example, $\est \in \mathbb{R}^{{d \choose 2}}$ might be the upper-triangular elements of a block diagonal correlation matrix for $d$ variables.
Our goal is to learn the parameters of a statistical model when these parameters have an unknown block diagonal sparsity structure.
Since the blocks are not known ahead of time, learning the block structure amounts to learning clusters of model parameters.
This setting is reminiscent of \textit{community detection} \citep{porter2009communities} in networks, but we simultaneously learn both the (real valued) edges of the network and the underlying communities, which are perfectly separated.

\begin{figure}[H]
 \centering
\begin{subfigure}[t]{0.3\textwidth}
\centering
\includegraphics[width=\linewidth, height=\linewidth]{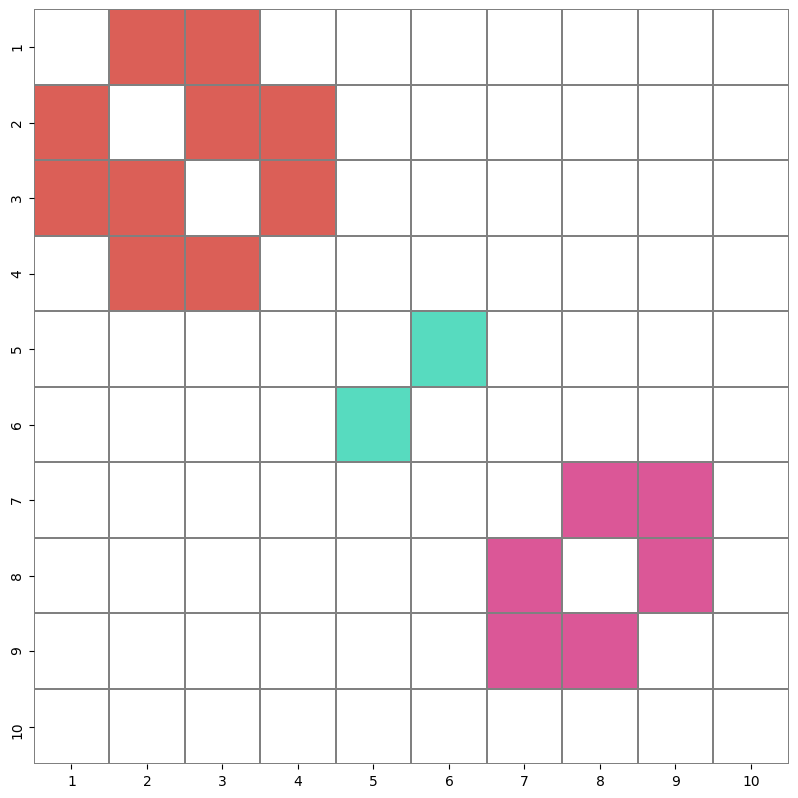}
\caption{
Block diagonal hollow symmetric matrix.
}
\label{fig:bd_sym_ex}
\end{subfigure}
\hfill
\begin{subfigure}[t]{0.3\textwidth}
\centering
\includegraphics[width=.6\linewidth, height=\linewidth]{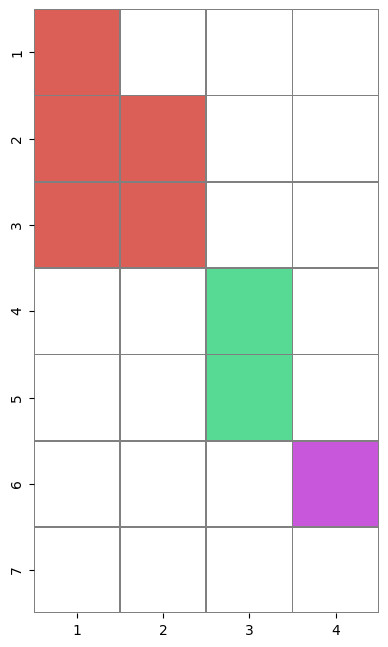}
\caption{
Block diagonal rectangular matrix. 
}
\label{fig:bd_rect_ex}
\end{subfigure}
\hfill
\begin{subfigure}[t]{0.3\textwidth}
\centering
\includegraphics[width=1.2\linewidth, height=\linewidth]{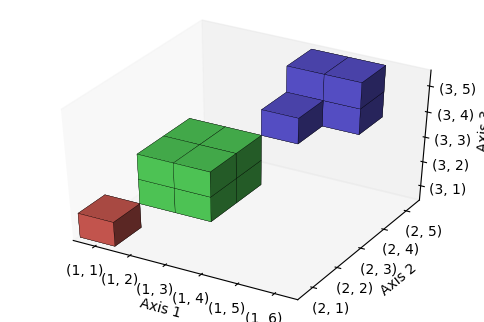}
\caption{
A 3 dimensional block diagonal multi-array.
}
\label{fig:bd_multi_array_ex}
\end{subfigure}
\caption{
Examples of block diagonal arrays.
Each corresponds to a graph that has 4 connected components, one of which is an isolated vertex.
Permuting the rows/columns of the matrices in Figures \ref{fig:bd_sym_ex}/\ref{fig:bd_rect_ex} (or the axes in Figure \ref{fig:bd_multi_array_ex}) results in graphs with the same number of connected components.
}
\label{fig:block_diag_examples}
\end{figure}

Block diagonal sparsity arises in several settings that are illustrated in Figure \ref{fig:block_diag_examples}.
\begin{enumerate}
\item Hollow symmetric matrices (0s on the diagonal) in $\mathbb{R}^{d \times d}$; the entries of $\est \in \mathbb{R}^{{d \choose 2}}$ parameterize the edges of the adjacency matrix of a graph with $d$ nodes.
For example, $\est$ might be the upper-triangular elements of a $d \times d$ dimensional correlation matrix.

\item Rectangular matrices in $\mathbb{R}^{m \times d}$; the entries $\est \in \mathbb{R}^{m \cdot d}$ parameterize the edges of a bipartite graph whose vertex sets are the rows and columns of the matrix. 
For example, $\est$ might be the regression coefficient matrix for a multiple response regression problem with $m$ responses and $d$ covariates.

\item Multi-arrays in $\mathbb{R}^{d\spa{1} \times \dots \times d\spa{V}}$; the entries of $\est \in \mathbb{R}^{\prod_{v=1}^V d\spa{v}}$ parameterize the hyperedges of a hypergraph  \citep{berge1984hypergraphs, zhou2006learning}.
See Appendix \ref{a:bd_rect_multi_array}. 

\end{enumerate}

There is a growing literature on block diagonal estimation including:
covariance matrix and \textit{graphical model} estimation \citep{marlin2009sparse, pavlenko2012covariance, tan2015cluster, hyodo2015testing, sun2015inferring, egilmez2017graph,  devijver2018block, kumar2019unified, broto2019block}, 
community detection \citep{nie2016constrained},
\textit{co-clustering} \citep{han2017bilateral, nie2017learning},
\textit{subspace clustering} \citep{feng2014robust, lu2018subspace},
\textit{principal components analysis} \citep{asteris2015sparse},
bipartite \textit{cross-correlation clustering} \citep{dewaskar2020finding},
neural network regularization \citep{tam2020fiedler},
and 
\textit{multi-view clustering} \citep{carmichael2020learning}.

A variety of approaches are used to estimate block diagonally structured parameters.
Some methods exploit the structure of particular statistical models \citep{asteris2015sparse, tan2015cluster, devijver2018block}.
Bayesian approaches to block diagonal estimation are based on priors that promote block diagonal structure \citep{mansinghka2006structured, marlin2009sparse}.
Constrained optimization approaches have been developed that constrain the eigenvalues of the \textit{graph Laplacian}  \citep{nie2016constrained, nie2017learning, egilmez2017graph, kumar2019unified} or the \textit{symmetric, normalized graph Laplacian} \citep{carmichael2020learning}.
We follow a similar optimization strategy by developing estimators that penalize Laplacian eigenvalues; see Problem \eqref{prob:con_comp_pen} below.

\subsection{Contributions and outline}

Inspired by the success of folded concave penalties for sparse vector estimation \citep{fan2001variable, zou2008one, zhang2012general, loh2015regularized, liu2017folded} -- especially \citep{fan2014strong} -- we develop the \textit{folded concave Laplacian spectral} (\penabr) penalty for block sparsity estimation (Section \ref{s:fcls_pen_opt}).
While this penalty is non-convex, we show it can be majorized by a (positively) weighted $L_1$ penalty.
Section \ref{ss:lla} develops a \textit{majorization-minimization} algorithm for this penalty, which we call a \textit{local linear approximation} algorithm (LLA) in the spirit of \citep{zou2008one}.

By making use of concave penalties such as the \textit{smoothly clipped absolute deviation} (SCAD) penalty, Section \ref{s:theory} shows  the output after \textbf{two steps} of the LLA algorithm obeys the \textit{strong oracle property}, as long as we have a ``good enough" initializer.
In addition to being an appealing statistical guarantee, this result has the computational implication that we may only need to take a few (e.g. 2) steps of the LLA algorithm.

Sections \ref{s:bd_shrink} and \ref{s:regression_examples} illustrate the theory in several settings, including block diagonal shrinkage, linear regression and logistic regression.
These sections provided strong non-asymptotic guarantees that often apply in ultra-high dimensional settings when $\log(d) = O(n^{\alpha})$ for some $\alpha \in (0, 1)$.
Section  \ref{s:sim} shows empirically that our proposed estimator outperforms competing entrywise penalties \citep{fan2014strong}.

Additional details such as proofs, technical results and additional computational details are provided in the appendix.
While the focus of this paper is on symmetric block diagonal matrices (Figure \ref{fig:bd_sym_ex}), our methodology extends naturally to rectangular matrices (Figure \ref{fig:bd_rect_ex}) and multi-arrays (Figure \ref{fig:bd_multi_array_ex}).
These extensions are sketched in Appendix \ref{a:bd_rect_multi_array}.
Section 
A python package implementing the methods in this paper is provided at \url{https://github.com/idc9/fclsp}
and code to reproduce the simulations in Section \ref{s:sim} is available at \url{https://github.com/idc9/repro\_lap_reg}.

\subsection{Related literature}
To our knowledge, we provide the first general purpose framework for block diagonal estimation that comes with strong statistical and computational guarantees.
Some existing bespoke methods for particular problems \citep{tan2015cluster, asteris2015sparse, devijver2018block} come with strong guarantees, but these methods do not naturally extend to other statistical models.
The existing literature on Laplacian spectral constraints typically provides weak computational guarantees (e.g. eventual convergence to a stationary point), but no statistical guarantees  \citep{nie2016constrained, egilmez2017graph, nie2017learning, kumar2019unified, carmichael2020learning}.
The recent work of \cite{tam2020fiedler} proposes a Laplacian spectral penalty for neural network regularization that is a special case of the \penabr\ penalty for a single Laplacian eigenvalue.

The \penabr\ penalty should not be confused with other penalties that make use of graph information \citep{smola2003kernels, ando2007learning, sharpnack2012sparsistency, jiang2013graph, li2020graph}. 
In these approaches a fixed graph known ahead of time and is used to incorporate prior graph structural information.
In our setting the graph is learned and the Laplacian is used to uncover unknown graph structure.

\subsection{Network constructions and notation}

Suppose $A\in \mathbb{R}^{d \times d}$ is the adjacency matrix of a graph with real valued edges and no self-loops (equivalently $A$ is a hollow symmetric matrix).
Let $\matsupp{A} \in \{0, 1\}^{d \times d}$ be the support of $A$ i.e. $\matsupp{A}_{ij} = \ind{A_{ij} \neq 0}$.
We say the connected components of $A$ are the connected components of $\matsupp{A}$.
In this case, the rows/columns $A$ can be permuted to form a block diagonal matrix.
We therefore refer to the \textit{blocks} of $A$ as these connected components.
 Note the number of connected components of $A$ is equal to the number of (non-zero) blocks plus the number of zero rows.
We will need the following notion of \textit{block support} i.e. the superset of the support that includes all within-block edges.
\begin{definition}
For a hollow symmetric matrix $A\in \mathbb{R}^{d \times d}$ the block support matrix $\blocksupp{A} \in \mathbb{R}^{d \times d}$ is the hollow symmetric matrix such that
$$
\blocksupp{A}_{ij} = \ind{\text{nodes } i \text{ and } j \text{ are in the same connected component of } \matsupp{A}}.
$$
\end{definition}

We next review a few basic facts about the graph Laplacian that can be found in \citep{von2007tutorial}.
For a symmetric matrix  $A \in \mathbb{R}^{d \times d}$ let 
$$L(A) :=  \text{diag}(A\mathbf{1}_d  )- A,$$
be the Laplacian matrix.
For a weighted graph with positive edges whose adjacency matrix is given by $A \in \mathbb{R}^{d \times d}_+$, the number of zero eigenvalues of $L(A)$ is equal to the number of connected components of $A$.
Furthermore, the indicator vectors (i.e. the vectors in $\{0, 1\}^d$ such that there is a 1 in the $j$th entry of the $k$ eigenvector if the $j$th node is a member of the $k$th connected component) of these connected components form an orthonormal basis for the kernel of $L(A)$.

For a vector $x \in \mathbb{R}^{ {d \choose 2}}$, let $\mathcal{A}(x) \in \mathbb{R}^{d \times d}$ be the hollow symmetric matrix whose upper triangular elements are given by $x$.
We will also write  $\mathcal{L}(x) := L(\mathcal{A}(x)) \in \mathbb{R}^{d \times d}$. 
To lighten notation, we will often write $x_{(ij)} := \mathcal{A}(x)_{ij}$ and $D := {d \choose 2}$.
To index the edges we will write either $\{ x_{\ell} | \ell \in [D] \}$ or $\{ x_{(ij)} | (ij) \in [D] \}$ i.e. $(ij)$ refers to the unique entry of $x$ corresponding to the edge between nodes $i$ and $j$.

\subsection{Notation} \label{ss:notation}

Let $[n] := \{1, \dots, n\}$ and let $\ind{\cdot}$ be the indicator function i.e. it is 1 if $\cdot$ is true and 0 otherwise.
Let $\mathbf{1}_d \in \mathbb{R}^d$ be the vector of ones.
For a vector $x$, $||x||_{q}$ refers to the usual $L_q$ norm.
We write $||x||_{\text{max}}$ and $||x||_{\text{min}}$ for the maximal (minimal) entry in absolute value.

For a symmetric matrix $A \in \mathbb{R}^{d \times d}$ let $\lambda(A) \in \mathbb{R}^{d}$ be the vector of eigenvalues of $A$.
We will write $\lambda_j$ for the $j$th largest eigenvalue and $\lambda_{(j)}$ for the $j$th smallest eigenvalue.
For the largest and smallest eigenvalues we will also use the notation $\lambda_{\text{max}} $ and $\lambda_{\text{min}}$.
For a matrix $A \in \mathbb{R}^{R \times C}$, the $(\kappa \to \rho)$ operator norm is 
$$\opnorm{A}{\kappa}{\rho} := \sup_{0 \neq x \in \mathbb{R}^C} \frac{||Ax||_{\rho}}{||x||_{\kappa}}$$
We will write $\opn{A} := \opnorm{A}{2}{2}$ for the standard operator norm (i.e. the largest singular value) and $\opone{A} := \opnorm{A}{1}{1} =  \max_{r \in [R]} \sum_{c=1}^C |X_{rc}|$ for the operator one norm.

For a matrix $X \in \mathbb{R}^{n \times d}$, $X(i, :) \in \mathbb{R}^{d}$ denotes the $i$th row and $X_j \in \mathbb{R}^n$ the $j$th column.
For a subset $A \subseteq [d]$ we let $X_A \in \mathbb{R}^{n \times |A|}$ be the sub-matrix of the columns in $A$.
When $x \in \mathbb{R}^d$ is a vector and $f: \mathbb{R} \to \mathbb{R}$ is a function we will overload notation and write $f(x) \in \mathbb{R}^d$ for $f$ applied element-wise to $x$ e.g. $|x| \in \mathbb{R}^d$ means the entrywise absolute value of $x$.

We say a random variable $X$ has a $\sigma$ sub-Gaussian distribution if $\prob{|X| \ge t} \le 2 \exp( - \frac{t}{2\sigma^2 })$ for all $t > 0$ e.g. see \citep{vershynin2018high}.
We say an event $E$ occurs \textit{with overwhelming probability} if there exist positive constants $N, c, C$ such that $\prob{E} \ge 1 - C e^{-c n}$ for all $n \ge N$ where $n$ is typically the sample size.
We will use $c, C, c_i$ to refer to absolute constants that may change from line to line.

\section{Block diagonal sparsity with Laplacian spectral regularization}  \label{s:fcls_pen_opt}

Consider estimating a model parameter, $\est \in \mathbb{R}^D$, by minimizing some loss function $\loss{\est}$.
Throughout the body of the paper we assume the model parameter can be viewed as the edges of an adjacency matrix, $\mathcal{A}(\est) \in \mathbb{R}^{d \times d}$ where $D = {d \choose 2}$. 
Suppose we know a priori that the target model parameter has a block diagonal sparsity structure.
We would then like to solve the following block diagonally constrained problem,
\begin{equation} \label{prob:con_comp_constr}
\begin{aligned}
& \underset{\est \in \mathbb{R}^{D}}{\textup{minimize}}  & & \ell(\est) \\ 
& \text{subject to } & & \mathcal{A}(|\est|) \text{ has at least K connected components},
\end{aligned}
\end{equation}
where $K$ is a tuning parameter.
This constraint set is non-convex and, even worse, seems combinatorial.

Problem \eqref{prob:con_comp_constr} can be made amenable to continuous optimization approaches by recalling that the number of zero eigenvalues of the Laplacian count the number of connected components of a graph.
Thus the (combinatorial) sparsity constraint of \eqref{prob:con_comp_constr} is equivalent to the (continuous) spectral constraint, $\lambda_{(K)}(\mathcal{L}(|\est|)) = 0$  \citep{nie2016constrained, nie2017learning, kumar2019unified, carmichael2020learning}.
Instead of solving such a constrained problem, we develop a related penalized approach that leads to an efficient algorithm with strong statistical and computational guarantees. 

\subsection{Folded concave Laplacian spectral penalty}

For a function $g_{\tuneparam}: \mathbb{R} \to \mathbb{R}$ and matrix $X \in \mathbb{R}^{d \times d}$ we write $ g_{\tuneparam} \circ \lambda(X) =  \sum_{i=1}^d g_{\tuneparam}(\lambda_i(X))$.
\begin{definition}
For a vector $\est \in  \mathbb{R}^{{d \choose 2}}$ and concave increasing function  $g_{\tuneparam}: \mathbb{R} \to \mathbb{R}$, the folded concave Laplacian spectral (\penabr) penalty is given by
\begin{equation} \label{eq:fcls}
 \fcls{\tuneparam}{\est} := \frac{1}{2} g_{\tuneparam} \circ \lambda (\mathcal{L}(|\est|))  =: \frac{1}{2} \lapSpec{\tuneparam}{|\est|}.
\end{equation}
\end{definition}

\begin{remark}
We can check that $ \lapSpec{\tuneparam}{\est}$ is concave when $g_{\tuneparam}$ is concave e.g. by Theorem 7.17 of \citealt{beck2017first}.
Therefore  $\fclsFunc{\tuneparam}(\cdot)$ is concave in $|\est|$  hence the name ``folded concave".
Note that if $g_{\tuneparam}$ is the identity then $\fcls{\tuneparam}{\est} = ||\est||_1$ i.e. the \penabr\ penalty reduces to the LASSO.
\end{remark}

The \penabr\ penalty encourages sparsity in the Laplacian eigenvalues thus encouraging $\mathcal{A}(\est)$ to have multiple connected components.
We therefore consider the following penalized problem
\begin{equation} \label{prob:con_comp_pen}
\underset{\est \in \mathbb{R}^{D}}{\textup{minimize}}  \;\; \ell(\est) + \frac{1}{2} g_{\tuneparam} \circ \lambda (\mathcal{L}(|\est|)) .
\end{equation}
Informally, the penalized Problem \eqref{prob:con_comp_pen} is to  the constrained Problem \eqref{prob:con_comp_constr} as Lasso/SCAD penalized problems are to $L_0$ constrained problems.

\subsection{Local linear approximation algorithm} \label{ss:lla}

We next derive a majorization-minimization (MM) algorithm for Problem \eqref{prob:con_comp_pen} \citep{lange2000optimization, sun2016majorization}. 
The key observation is that $\fclsFunc{\tuneparam}$  can be majorized by (positively) weighted Lasso function.

\begin{definition}
We say a surrogate function $Q(x | y)$ majorizes $f(x)$ if
$$
f(x) = Q(x | x) \text{ for all } x \text{ and }  f(x) \le Q(x | y) \text{ for all } x, y.
$$
\end{definition}

Since $\lapSpecFunc{\tuneparam}(\cdot)$ is concave we can construct a surrogate function via a standard linearization approach i.e. $\overline{Q}(\est | \estCurrent) := \nabla \lapSpec{\tuneparam}{\estCurrent}^T (\est- \estCurrent) +   \lapSpec{\tuneparam}{\estCurrent} =  \nabla \lapSpec{\tuneparam}{\estCurrent}^T  \est + \text{constants}(\estCurrent)$ is a surrogate for $\lapSpecFunc{\tuneparam}(\cdot)$.
The composition of this surrogate function with the absolute value function,  $Q(\est | \estCurrent) := \overline{Q}( |\est|  \Big| |\estCurrent| )$, then gives a surrogate for $\fclsFunc{\tuneparam}(\cdot)$.

Next we give an explicit form for the gradient, $ \nabla \lapSpec{\tuneparam}{\estCurrent}$.
For $w \in \mathbb{R}^K$ and $V \in \mathbb{R}^{d \times K}$ let $\lapcoeff{V}{w} \in \mathbb{R}_+^{D}$ be the vector such that
\begin{equation} \label{eq:lap_coeff_formula}
\lapcoeff{V}{w}_{(ij)} =   ||V(i, :) - V(j, :)||_{2, w}^2,
\end{equation}
where $||y||_{2, w}^2 := y^T \text{diag}(w) y$.
We call this vector the \textit{Laplacian coefficient}. 
We can check that
\begin{equation} \label{eq:lap_spect_grad}
 \nabla \lapSpec{\tuneparam}{\estCurrent} =   \lapcoeff{V^{\estCurrent}}{g'_{\tuneparam}(\lambda^{\estCurrent})},
\end{equation}
where  $V^{\estCurrent} \in \mathbb{R}^{D \times D}$ is any orthonormal matrix of eigenvectors of $\mathcal{L}(\estCurrent)$ with corresponding eigenvalues $\lambda^{\estCurrent} \in \mathbb{R}^{D}_+$ and $g'_{\tuneparam}$ denotes any super-gradient. 
Thus the surrogate function for $\fclsFunc{\tuneparam}$ at $\estCurrent$ is given by
\begin{equation} \label{eq:fcls_maj}
Q(\est | \estCurrent) =    \sum_{\ell=1}^D \lapcoeff{V^{|\estCurrent|}}{g'_{\tuneparam}(\lambda^{|\estCurrent|})}_{\ell} \; \cdot  |\est_{\ell}| + \text{constants}(b), 
\end{equation}
where $V^{|\estCurrent|}, \lambda^{|\estCurrent|}$ are eigenvectors/values of $\mathcal{L}(|\estCurrent|)$.

\begin{proposition} \label{prop:spect_pen_maj}
$Q(\cdot | \estCurrent) $ is convex and majorizes $\fclsFunc{\tuneparam}(\cdot)$ at $\estCurrent$.
The quantity $ \lapcoeff{V^{|\estCurrent|}}{g'_{\tuneparam}(\lambda^{|\estCurrent|})}$ does not depend on the choice of eigenvectors, $V^{|\estCurrent|}$.
Majorizing at $0$ leaves a Lasso i.e. $Q(\est | 0)  =  g'_{\tuneparam}(0) ||\est||_1$.
Additionally, any fixed point of the below MM algorithm is a stationary point of Problem \eqref{prob:con_comp_pen}.
\end{proposition}

Equipped with our weighted $L_1$ surrogate function we can use the following MM algorithm to hunt for stationary points of the penalized Problem \eqref{prob:con_comp_pen}.
Similar \textit{reweighted $L_1$} algorithms appear elsewhere in the literature \citep{krishnapuram2005sparse, zou2008one,  gasso2009recovering, candes2008enhancing, fan2014strong}. 

\begin{algorithm}[H] \label{algo:lla}
\DontPrintSemicolon
  
\KwInput{Tuning parameter value $\tuneparam \ge 0$ and initializer $\est\spa{0} \in \mathbb{R}^D$}
\For{s=0, 1, 2, \dots}
{

$V\spa{s}, \lambda\spa{s} \leftarrow$ eigenvectors and eigenvalues of $\mathcal{L}\left(|\est\spa{s}| \right)$ 

$ M\spa{s} \leftarrow \lapcoeff{V\spa{s} }{g_{\tuneparam}'(\lambda\spa{s}) }$ \tcp*{Obtain weights for Lasso surrogate using \eqref{eq:lap_coeff_formula}}

\begin{equation}  \label{prob:weighted_lasso_maj}
\est\spa{s + 1} =  \underset{\est \in \mathbb{R}^D}{\textup{argmin}} \;\; \ell(\est)  +  \frac{1}{2} |\est|^T M\spa{s} 
\end{equation}

}
\caption{Local linear approximation algorithm for \penabr\ penalized Problem \eqref{prob:con_comp_pen}}
\end{algorithm}

Each loop of Algorithm \ref{algo:lla} requires computing an eigen-decomposition and a solution to a (weighted) Lasso penalized problem.
The theory in Section \ref{s:theory} shows we may only need a few (e.g. 2) LLA iterations.
Additional computational details of Algorithm \ref{algo:lla} are discussed in Appendix \ref{a:tune_param_ubd}. 

A key insight into this algorithm is the following.
Suppose a point $\estCurrent$ is such that $\mathcal{A}(\estCurrent)$ has exactly $K$ connected components.
Further suppose the $K+1$st smallest eigenvalue is large enough and $g_{\tuneparam}$ is such that $g'_{\tuneparam}(\lambda_{(K+1)}(\mathcal{L}(|\estCurrent|))) = 0$ (think of penalties such as SCAD that are eventually flat).
In this case we can check that up to additive constants
$$
Q(\est | \estCurrent) =
\frac{1}{2} g'_{\tuneparam}(0) \sum_{ (ij) \in [D]}  \left(\frac{1}{|C(i)|} + \frac{1}{|C(j)|}\right)  \cdot \ind{i \text{ and } j \text{ are in different connected components} }  \cdot |\est_{(ij)}| 
$$
where $|C(i)|$ denotes the number of vertices in the connected component that node $i$ belongs to.
In other words $Q(\est | \estCurrent)$ puts a Lasso penalty on all the edges that go between connected components in $\mathcal{A}(\estCurrent)$, but does not penalize any edge within a in connected component.
Thus if Algorithm \ref{algo:lla} reaches a point with $K$ connected components, then the \penabr\ penalty may force the next iteration to have the same connected components (see Theorem  \ref{thm:orc_is_stat_point}).
Similar behavior occurs when $\estCurrent$ is close to a graph that has $K$ connected components (see Theorem \ref{thm:init_from_close_give_lasso_orc}).


\section{Oracle properties of the \penabr\ penalty} \label{s:theory}

This section studies two questions:
\begin{enumerate}

\item When is the \textit{block oracle solution} (made precise below) a stationary point of the \penabr\ penalized Problem \eqref{prob:con_comp_pen}?

\item When does the LLA algorithm find the block oracle solution after exactly two steps?

\end{enumerate}

Section \ref{ss:theory__stat_pt} answers the first question, thus providing theoretical guarantees for a theoretical estimator that we are by no means guaranteed to find.
Section \ref{ss:theory__two_step_lla} answers the second question, thus providing theoretical guarantees for estimators actually used in practice.
Recall that similar guarantees for an entrywise folded concave concave penalty require only one LLA step \citep{fan2014strong}.
Sections \ref{ss:theory__stat_pt} and \ref{ss:theory__two_step_lla} provide deterministic conditions under which their claims follow.
Applying these results requires verifying these conditions hold with high-probability for particular statistical models, which is the topic of the following Sections \ref{s:bd_shrink} and \ref{s:regression_examples}.

\subsection{When the oracle is a stationary point} \label{ss:theory__stat_pt}

Our goal is to estimate some ``target" parameter, $\esttarg \in \mathbb{R}^D$ where $D = {d \choose 2}$.
For example, $\esttarg$ might be the upper triangular entries of a correlation matrix for a $d$ dimensional random vector.

The results below depend on the following graph quantities.
Let $\binGraphTarg$ be the  binary graph  whose adjacency matrix is given by the support graph, $\matsupp{\mathcal{A}(\esttarg)}$.
Let $\suppTarg$ be the connected component support set i.e. $ \suppSetTarg = \{(i,j) : \blocksupp{\mathcal{A}(\esttarg)}_{ij} = 1, i < j\}$; this is the support set we will attempt to estimate.
Also let $\nccTarg$ be the number of connected components of $\binGraphTarg$,  let $\maxccTarg$ be the largest number of nodes in a connected component and let $\minccTarg$ be the smallest size of a connected component that is not an isolated vertex.

A \textit{block oracle} estimator, $\estorc$, is a solution to the following constrained problem
\begin{equation} \label{prob:block_oracle}
\begin{aligned}
& \underset{\est \in \mathbb{R}^{D} }{\textup{minimize}}  & & \ell(\est) \\ 
& \text{subject to } & & \est_{\suppSetTargPow{C}} = 0.
\end{aligned}
\end{equation}
In other words, a block oracle knows the true between block edges for the target parameter.
Note this estimator is unaware of any within block edges that are zero.
For simplicity we assume Problem \eqref{prob:block_oracle} has a unique solution, $\estorc$, which satisfies the first order necessary conditions $\nabla_{\suppSetTarg} \ell(\estorc) = 0$.

We next define two sets specifying the values of $\est$ that satisfy regularity conditions.
Let
\begin{equation}\label{eq:sets_block_orc}
\setNiceGrad{t} :=  \left\{ \est \text{ s.t. } || \nabla_{\suppSetTargPow{C}} \ell(|\est|) ||_{\text{max}}  < t \right\}, 
\end{equation}
be the set of points whose gradients are small for the between block edges $\suppSetTargPow{C}$.
Let $\ccTarg{1}, \dots, \ccTarg{\nccNzTarg} \subseteq \suppSetTarg$ be the indices of the $\nccNzTarg \le \nccTarg$ non-zero connected components of $\binGraphTarg$ (recall isolated vertices correspond to zero rows of $\mathcal{A}(\esttarg)$).
Let
\begin{equation}\label{eq:each_cc_has_big_ssm_eval}
\setBigEvalComp{t} :=  \left\{ \est \text{ s.t. } \min_{k \in [\nccNzTarg]}   \lapSmEval{2}{|\est_{ \ccTarg{k}}|} \ge t \right\},
\end{equation}
be the set of points with a large spectral gap.

The theory in this section applies to a class of concave penalties defined in \citep{fan2014strong}.
\begin{definition} \label{def:scad_like_pen_func}
Let $a_0 \ge a_1 > 0$ and $b_2 > b_1 > 0$. 
A SCAD-like 
concave penalty function, $g_{\tuneparam}: \mathbb{R}_+ \to \mathbb{R}$, satisfies the following
\begin{enumerate}
\item $g_{\tuneparam}(t)$ is increasing and concave for $t \in [0, \infty)$ and $g_{\tuneparam}(0) = 0$,

\item $g_{\tuneparam}(t)$ is differentiable for $t \in (0, \infty]$ with $g_{\tuneparam}'(0) := g_{\tuneparam}'(0+) = a_0 \tuneparam$, 

\item $g_{\tuneparam}'(t) \ge a_1 \tuneparam$ for $t \in (0, b_1 \tuneparam]$, 

\item $g_{\tuneparam}'(t) =0 $ for $t \in [b_2 \tuneparam, \infty)$. 

\end{enumerate}
\end{definition}
\noindent For example, if $g_{\tuneparam}(\cdot)$ is the SCAD penalty where
$$
g_{\tuneparam}'(t) = \tuneparam \ind{t \le \tuneparam} + \frac{[a \tuneparam - t]_+}{a - 1} \ind{t > \tuneparam}, \qquad \text{ for some } a > 2
$$ 
then $a_0= a_1 = 1$, $b_1=1$, $b_2= a$.
Other penalties satisfying this definition include the MCP and hard-thresholding penalty, see\footnote{Beware our notation is slightly different.} \citep{fan2014strong}.

Throughout our analysis we make the following basic assumption\footnote{We can relax this uniqueness assumption by checking the regularity conditions hold for \textit{any} solution to Problem \eqref{prob:block_oracle}.}  
\begin{assumption} \label{assu:basic}
The loss function $\ell(\cdot)$ in Problem \eqref{prob:con_comp_pen} is convex,
$g_{\tuneparam}(\cdot)$ satisfies Definition \ref{def:scad_like_pen_func},
$\estorc$ is the unique solution to Problem \eqref{prob:block_oracle},
and $\mathcal{A}(\esttarg)$ has $\nccTarg$ connected components as above.
\end{assumption}

\begin{theorem} \label{thm:orc_is_stat_point}
Suppose Assumption \ref{assu:basic} is satisfied.
If $\estzero \in \mathbb{R}^{D}$ is any point whose block support is a subset of $\esttarg$'s block support, i.e. $ \blocksupp{\mathcal{A}(\estzero)} \subseteq  \blocksupp{\mathcal{A}(\esttarg)}$, and satisfies
\begin{equation} \label{eq:orc_is_stat_point__big_eval_compwise}
\estzero \in \setBigEvalComp{b_2 \tuneparam}
\end{equation}
and the block oracle satisfies
\begin{equation} \label{eq:orc_has_nice_grad} 
\estorc \in \setNiceGrad{\frac{a_0 \tuneparam}{\maxccTarg}}
\end{equation}
then taking one LLA step from $\estzero$ results in $\estorc$.
In particular, if 
\begin{equation} \label{eq:orc_is_stat_point__orc_in_set_for_stat_point}
\estorc \in  \setBigEvalComp{b_2 \tuneparam} \cap \setNiceGrad{\frac{a_0 \tuneparam}{\maxccTarg}},
\end{equation}
then $\estorc$ is fixed point of the LLA algorithm and therefore is a stationary point of Problem \eqref{prob:con_comp_pen}.
\end{theorem}

Let us unpack the conditions of this theorem.
The block support assumption on $\estorc$ guarantees the connected components of $\mathcal{A}(\estorc)$ are a subset of the target's connected components.
Combining this with the assumption $\estorc \in \setBigEvalComp{b_2 \tuneparam}$ guarantees these connected components are \textit{exactly} the target's connected components.\footnote{E.g. this rules out the case that one connected component is zeroed out while another connected component is split into two components.} 
This assumption -- in conjunction with the spectral properties of the Laplacian and form of $g_{\tuneparam}$ -- guarantee the Laplacian coefficient,  $\lapcoeff{V^{\text{orc}} }{w^{\text{orc}} }$,  computed at the oracle puts zero penalty on the within block edges and a large penalty on the between block edges (recall Section \ref{ss:lla}). %
Finally, the condition $\estorc \in  \setNiceGrad{a_0 \tuneparam}$ guarantees the between block edges are killable by the latter Lasso terms.

Next we give a sufficient condition to verify elements of $\setBigEvalComp{t}$.
Let
\begin{equation}
\targgap := \lapSmEval{\nccTarg + 1}{|\esttarg|} = \min_{k \in [\nccNzTarg]}  \lapSmEval{2}{|\esttarg_{ \ccTarg{k}}|}.
\end{equation}
be the \textit{target spectral gap}, which quantifies the minimal signal strength of the target graph.\footnote{This quantity is the analog of $||\esttarg_{\suppSizeTarg}||_{\text{min}}$ from \citep{fan2014strong}.} 
This quantity depends on the magnitude of the (smallest) entries of $\esttarg$ as well as on the topology of the binary graph $\binGraphTarg$.
For binary graphs, finding lower bounds on the second smallest Laplacian eigenvalue -- also called the \textit{algebraic connectivity} or \textit{Fiedler value} -- is an active area of research \citep{fiedler1973algebraic, lu2007lower, de2007old, rad2011lower}.
For some important special cases we can give lower bounds on $\targgap$.
When the non-zero connected components of $\binGraphTarg$ are fully connected we have\footnote{We have used the elementary Fact \ref{fact:evals_increase_for_pos_mats} stated in Appendix \ref{a:proofs_lla}.} 
\begin{equation} \label{eq:targ_gap_from_fully_conn}
\targgap \ge \minccTarg  \minSuppMagTarg,
\end{equation}
which shows the more nodes in a connected component, the larger the signal strength.
We can relax the fully connected assumption for \eqref{eq:targ_gap_from_fully_conn}  e.g. zeroing out several edges should not dramatically change $\targgap$. 
See Appendix \ref{as:graph_spectra} for additional details. 

We next give the non-asymptotic probability that $\estorc$ is a stationary point.
\footnote{We use Weyl's inequality and the fact that the operator norm of a block diagonal matrix is the largest operator norm of the blocks thus $\est_{\suppSetTargPow{C}} = 0 \text{ and } \opn{\mathcal{L}(|\est| - |\esttarg|)} \le \targgap - t \implies \est \in \setBigEvalComp{t}$} 
Let
\begin{equation}
\probNiceGradStatPt  :=  \prob{ || \nabla_{\suppSetTargPow{C}} \ell(\estorc) ||_{\text{max}}  >  \frac{a_0 \tuneparam}{\maxccTarg}  },
\qquad
\probSmallResidStatPt   :=  \prob{\opn{\mathcal{L}(|\estorc| - |\esttarg |)} \ge \targgap - b_2 \tuneparam}.
\end{equation}

\begin{corollary} \label{cor:orc_is_stat_pt_whp}
Suppose Assumption \ref{assu:basic} is satisfied.
Then $\estorc$  is a stationary point of Problem \eqref{prob:con_comp_pen} with probability at least $1 - \probNiceGradStatPt - \probSmallResidStatPt$.
\end{corollary}

\subsection{When the LLA algorithm converges to the block oracle in two steps}\label{ss:theory__two_step_lla}

We first make a few constructions that capture regularity conditions needed below.
The following set captures the set of points with a large Laplacian spectral gap 
\begin{equation}\label{eq:big_eval_gap}
\setEvalGap{s, t} :=  \left\{ \est \text{ s.t. } \lapSmEval{\nccTarg}{|\est|}  \le s \text{ and }  \lapSmEval{\nccTarg+1}{|\est|} \ge t  \right\}.
\end{equation}

For a weight vector  $w \in \mathbb{R}^D_+$ with $w_{\suppSetTargPow{C}} = 0$ consider the following \textit{Lasso block oracle} problem
\begin{equation} \label{prob:weighted_lasso_oracle}
\begin{aligned}
& \underset{\est \in \mathbb{R}^{D} }{\textup{minimize}}  & & \ell(\est) +  \frac{1}{2} \sum_{(ij) \in \suppSetTarg} w_{(ij)} |\est_{(ij)}| \\ 
& \text{subject to } & & \est_{\suppSetTargPow{C}} = 0.
\end{aligned}
\end{equation}
This problem reduces to the block oracle Problem \eqref{prob:block_oracle} when $w = 0$.
The following are two sets of block oracle solutions with ``small" Lasso penalties,
\begin{equation} \label{eq:def_lassoOrcSetOp}
\begin{aligned}
\lassoOrcSetOp{\residDiffGapRatio, \tuneparam} := \Big\{ \est \in  \mathbb{R}^{D}|
& \est \text{ is a minimizer of } \eqref{prob:weighted_lasso_oracle}  \text{ for } w \in \mathbb{R}^{D}_+, \text{ s.t. } w_{\suppSetTargPow{C}} = 0 \text{ and }  ||w||_{\text{max}} \le 2^5 a_0 \tuneparam \residDiffGapRatio^2
 \Big\}. 
\end{aligned}
\end{equation}
and
\begin{equation} \label{eq:def_lassoOrcSetFrob}
\begin{aligned}
\lassoOrcSetFrob{\residDiffGapRatio, \tuneparam} := \Big\{ \est \in  \mathbb{R}^{D}|
& \est \text{ is a minimizer of } \eqref{prob:weighted_lasso_oracle}  \text{ for } w \in \mathbb{R}^{D}_+ , \text{ s.t. } w_{\suppSetTargPow{C}} = 0 \text{ and}\\
&  ||w||_{\text{max}} \le 2^5 a_0 \tuneparam \residDiffGapRatio^2,
\; ||w ||_2 \le2^{9/2} \maxccTargPow{1/2} a_0 \tuneparam \residDiffGapRatio^2,
\text{ and }  ||w ||_1 \le 2^4 \maxccTarg a_0 \tuneparam   \residDiffGapRatio^2
 \Big\}. 
\end{aligned}
\end{equation}

\begin{theorem} \label{thm:init_from_close_give_lasso_orc} 
Suppose Assumption \ref{assu:basic} is satisfied and Lasso oracle Problem \eqref{prob:weighted_lasso_oracle} has a unique solution for any weight vector.
Fix $\residDiffGapRatio \ge 0$.
Let $\estone, \esttwo$ be the result of taking the first and second LLA steps from $\estzero$.
If
\begin{equation} \label{eq:init_from_close_give_lasso_orc__big_gap}
\estzero \in \setEvalGap{b_1 \tuneparam, b_2 \tuneparam} 
\end{equation}
\begin{equation}\label{eq:init_from_close_give_lasso_orc__small_ratio__op} 
\frac{\opn{\mathcal{L}(|\estzero| - |\esttarg|)}}{\targgap}  \le 
 \residDiffGapRatio \wedge  \frac{1}{2^6} \sqrt{\frac{a_1}{a_0 \maxccTarg}},
\end{equation}
\begin{equation}  \label{eq:init_from_close_give_lasso_orc__nice_grad}
\lassoOrcSetOp{\residDiffGapRatio, \tuneparam} \subseteq \setNiceGrad{\frac{a_1 \tuneparam}{4 \maxccTarg}},
\end{equation}
then  $\estone \in \lassoOrcSetOp{\residDiffGapRatio, \tuneparam}$.
If additionally
\begin{equation}\label{eq:init_from_close_give_lasso_orc__lasso_orc_uniformly_big_eval}
\lassoOrcSetOp{\residDiffGapRatio, \tuneparam}  \subseteq \setBigEvalComp{b_2 \tuneparam}
\end{equation}
then $\esttwo = \estorc$. Furthermore, taking an additional LLA step from $\esttwo$ results in the same estimate i.e. the LLA algorithm has converged to the block oracle estimate.

Furthermore, if \eqref{eq:init_from_close_give_lasso_orc__small_ratio__op} is replaced with
\begin{equation}\label{eq:init_from_close_give_lasso_orc__small_ratio__frob}
\min\left(||\mathcal{L}(|\estzero| - |\esttarg|)||_F, \nccTargPow{1/2} \opn{\mathcal{L}(|\estzero| - |\esttarg|)} \right)
\le
\targgap  \cdot \left( \residDiffGapRatio \wedge  \frac{1}{2^6} \sqrt{\frac{a_1}{a_0 \maxccTarg}} \right),
\end{equation}
then the above statements hold with $\lassoOrcSetFrob{\residDiffGapRatio, \tuneparam} $ in place of $\lassoOrcSetOp{\residDiffGapRatio, \tuneparam}$.

\end{theorem}

\begin{figure}[H]
 \centering
\centering
\includegraphics[width=\linewidth, height=.65\linewidth]{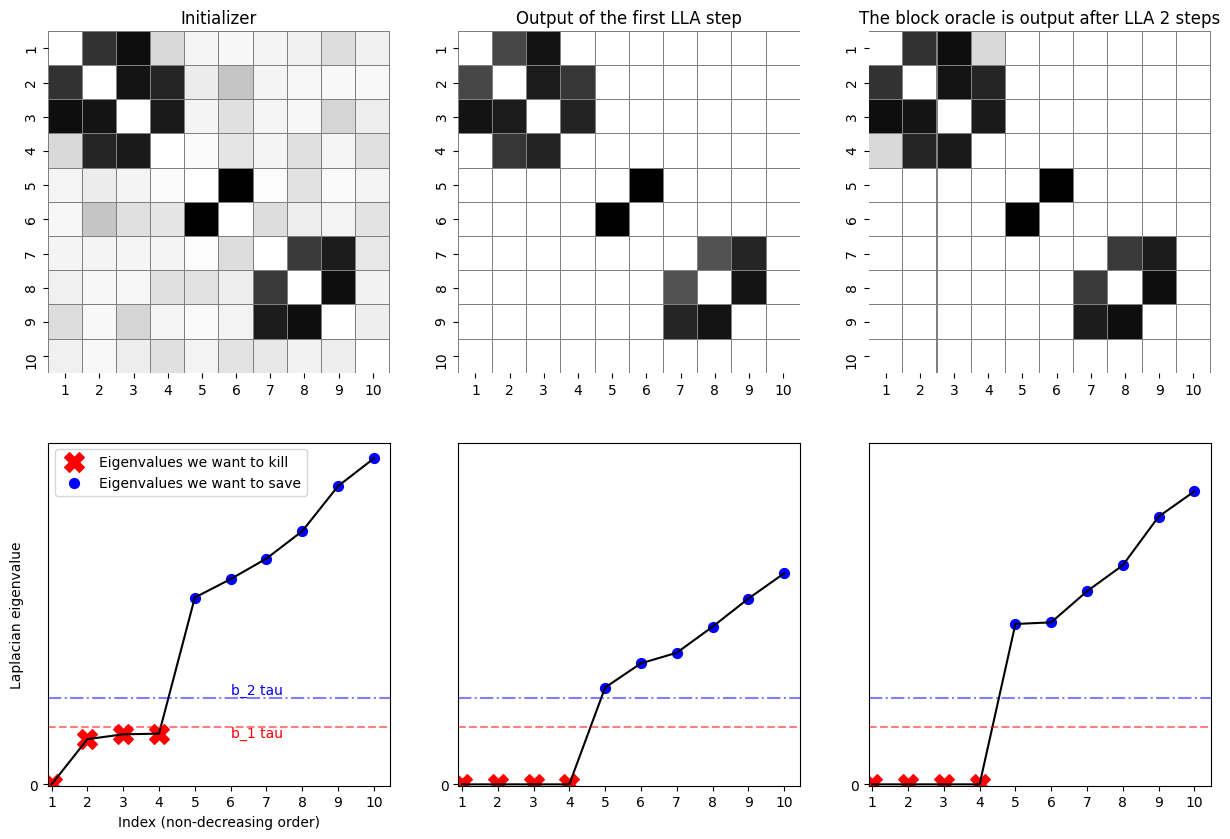}
\caption{
Two step convergence of the LLA algorithm where the target graph is given in Figure \ref{fig:bd_sym_ex}.
The columns show the initializer and the output of the first two LLA steps.
The first row show the matrices $\mathcal{A}(|\est\spa{s}|)$ at each step, $s=0, 1, 2$.
The second row shows the spectrum of the Laplacian, $\mathcal{L}(|\est\spa{s}|)$.
In Theorem \ref{thm:init_from_close_give_lasso_orc} the estimate, $\estone$, after the first LLA step is close to the block oracle i.e. all between block edges are zero but the within block edges are not exactly the block oracle.
The estimate, $\esttwo$, after the second LLA step is exactly the block oracle.
Note some edges in the block oracle may be incorrectly identified as non-zero e.g. edge $(1, 4)$ in the third column.
}
\label{fig:LLA_algo_two_steps_grid}
\end{figure}

Figure \ref{fig:LLA_algo_two_steps_grid} gives a visual depiction of the two step convergence in Theorem \ref{thm:init_from_close_give_lasso_orc}.
Let us unpack the conditions of this theorem.
Condition \eqref{eq:init_from_close_give_lasso_orc__big_gap} ensures $g_{\tuneparam}$ only penalizes the smallest $\nccTarg$ Laplacian eigenvalues of the initializer and not the larger ones.
Condition \eqref{eq:init_from_close_give_lasso_orc__small_ratio__op} says the initializer must be close enough to the target parameter such their Laplacian eigenvectors look similar.
Condition \eqref{eq:init_from_close_give_lasso_orc__nice_grad} guarantees the between block edges are killable by the Lasso penalty that comes from the Laplacian coefficient of the initializer.
Under these assumptions, taking one LLA step from $\estzero$ ends up in $\lassoOrcSet{\residDiffGapRatio, \tuneparam}$. 
Assumption \eqref{eq:init_from_close_give_lasso_orc__lasso_orc_uniformly_big_eval}  guarantees taking an LLA step from any element of $\lassoOrcSet{\residDiffGapRatio, \tuneparam}$ lands at the block oracle (by Theorem \ref{thm:orc_is_stat_point}).
The parameter $\residDiffGapRatio$ is introduced because making the set $ \lassoOrcSet{\residDiffGapRatio, \tuneparam}$ small enough so that the conditions hold with high probability requires a model specific analysis.

We next give a non-asymptotic probability bound that the LLA algorithm converges to the block oracle after two steps.
Let
\begin{equation} \label{eq:lla_probs}
\begin{aligned}
 \initRatioConsts & := \min \left(\residDiffGapRatio,  \frac{1}{2^6} \sqrt{\frac{a_1}{a_0 \maxccTarg}}, \frac{b_1}{b_1 + b_2}\right) (b_1 + b_2)    \\
\probGoodInitLLAOp & := \prob{
 \opn{\mathcal{L}(|\estinit| - |\esttarg|)} \ge \initRatioConsts \tuneparam } \\
\probNiceGradLLAOp & :=  \prob{ \sup_{\est \in \lassoOrcSetOp{\residDiffGapRatio, \tuneparam}} || \nabla_{\suppSetTargPow{C}} \ell(\est) ||_{\text{max}}  >  \frac{a_1 \tuneparam}{4 \maxccTarg}  }\\
\probSmallResidLLAOp & :=  \prob{  \sup_{\est \in \lassoOrcSetOp{\residDiffGapRatio, \tuneparam}}  \opn{\mathcal{L}(|\est| - |\esttarg |)} \ge  b_1 \tuneparam }.
\end{aligned}
\end{equation}
We also define $\probGoodInitLLAFrob$ similarly to $\probGoodInitLLAOp$, but with the left hand side of \eqref{eq:init_from_close_give_lasso_orc__small_ratio__frob} in place of the operator norm.
Also let $\probNiceGradLLAFrob$ be the quantity $\probSmallResidLLAOp$ but with   $ \lassoOrcSetFrob{\residDiffGapRatio, \tuneparam}$ in place if $ \lassoOrcSetOp{\residDiffGapRatio, \tuneparam}$ (similarly for $\probSmallResidLLAFrob$).

\begin{corollary} \label{cor:suff_cond_whp_for_init_from_close_give_lasso_orc}
Suppose Assumption \ref{assu:basic} is satisfied, the Lasso oracle Problem \eqref{prob:weighted_lasso_oracle} has a unique solution for any weight vector and
\begin{equation} \label{eq:targ_gap_at_least_tune_param}
\targgap \ge (b_1 + b_2) \tuneparam.
\end{equation}
Then the LLA algorithm initialized from $\estinit$ will converge to $\estorc$ in two steps with probability at least $1 - \probGoodInitLLAOp - \probNiceGradLLAOp - \probSmallResidLLAOp$.
This statement also holds with probability at least $1 - \probGoodInitLLAFrob - \probNiceGradLLAFrob - \probSmallResidLLAFrob$.
\end{corollary}

To apply this corollary have to verify regularity conditions hold for every element of $\lassoOrcSet{\residDiffGapRatio, \tuneparam}$.
If the loss function, $\ell(\cdot)$ is strongly convex when restricted to $\suppSetTarg$, then we can reduce these conditions to statements about $\estorc$ 
(see Corollary \ref{cor:convex_loss_bounds_lasso_orc_dist_to_orc} in Appendix \ref{as:lin_reg__prelim}) 
, which can then be controlled using standard concentration inequalities.
The logistic regression example in Section \ref{ss:theo_ex_logisitc} uses a more involved fixed point argument due to the lack of strong convexity.


\section{Block diagonal shrinkage} \label{s:bd_shrink}

This section shows how the theory from Section \ref{s:theory} applies when $\ell(\cdot)$ is the squared loss.
Suppose  $\estok$ is an ``ok" estimate (e.g. the sample covariance matrix) for a target parameter $\esttarg$ (e.g. the population covariance matrix).
We can obtain a better estimator by shrinking $\estok$ towards a block diagonal matrix by computing proximal operator of the \penabr\ penalty,
\begin{equation}  \label{prob:prox_fclsp}
\begin{aligned}
& \underset{\est \in \mathbb{R}^{{d \choose 2}}}{\textup{minimize}}  & & \frac{1}{2} ||\estok - \est||_2^2 + \fcls{\tuneparam}{\est}.
\end{aligned}
\end{equation}
This ``shrink an ok estimator" idea is used for classical sparse estimation problems such as sparse covariance estimation \citep{bickel2008covariance, rothman2009generalized}.
Instead of solving this problem exactly, we can use the LLA algorithm to approximate its solution.
Below we show that as long as  $\estok$ is an ``ok" estimate of $\esttarg$ and a ``good enough" initializer, $\estinit$, is used then the LLA algorithm obeys the strong oracle property after two steps.

The results of this section are directly applicable to many statistical models including covariance estimation  \citep{ravikumar2011high}, correlation matrix estimation for non-paranormal copula models \citep{liu2012high}, and multi-layer community detection.
While our results can be applied immediately to these settings, we will not state these extensions explicitly for the sake of brevity.

Throughout this section we will assume Assumption \ref{assu:basic} is satisfied.
The results below will depend on the target parameter $\esttarg$ via: the spectral gap $\targgap$, the block sparsity support size $\suppSizeTarg := |\suppSetTarg|$, the number of connected components, $\nccTarg$, and the number of nodes in the largest connected component $\maxccTarg$.
We will also need the following quantities based on the binary target graph $\binGraphTarg$: the number of non-zero connected components, $\nccNzTarg$ the number of non-isolated vertices $\NumNonIso$, and the largest node degree, $\maxDegTarg$.

Recall $c, C, c_i$ denote absolute constants that may change from line to line.


\subsection{Stationary point for independent, unbiased, sub-Gaussian residuals}\label{ss:bd_shrink__stat_pt}

In this section we work under the following assumption.
\begin{assumption} \label{assumpt:sub_g_indep_unbiased}
Assume the entries of  $\estok$ are independent, $\expect{\estok} = \esttarg$ and have  sub-Gaussian tails with variance proxy $\frac{\sigma^2}{n}$.
That is $\prob{|\estok_{\ell} - \esttarg_{\ell}| > t} \le 2 \exp \left(\frac{- n t^2}{2 \sigma^2} \right)$ for each $\ell \in [D]$.
\end{assumption}

While this setting is meant to illustrate the theory in the nicest statistical setting it is also motivated by the following multi-layer, weighted network community detection model (similar models are studied in \cite{levin2019recovering}). 
Suppose we observe $n$ independent networks with real valued edges $\{\mathcal{A}(x_i) \}_{i=1}^n$, where each $x_i \in \mathbb{R}^{D}$.
Assume an additive noise model $x_i = \esttarg + e_i$ where $e_i \in \mathbb{R}^D$ $i=1, \dots, n$ are independent, mean zero noise vectors whose entries follow sub-Gaussian distributions with variance proxy $\sigma^2$.
Let $\estok = \frac{1}{n} \sum_{i=1}^n x_i$. Note $\estok = \esttarg   + \epsilon$ where the entries of $\epsilon$ are independent, mean 0 and follow sub-Gaussian distributions with variance proxy $\frac{\sigma^2}{n}$.

\begin{theorem} \label{thm:bd_shrink__stat_pt__nice}  
Suppose Assumption \ref{assumpt:sub_g_indep_unbiased} holds, $\targgap \ge 2 b_2\tuneparam$ and
\begin{equation} \label{eq:bd_shrink__indep_unbiased__stat_pt__targ_gap__assumpt}
4  b_2 \tuneparam \ge  
c_1 \maxccTarg \left(\frac{\sigma}{\sqrt{n}} + ||\esttarg||_{\text{max}} \right) \exp \left( -\frac{c_2 n  \minSuppMagTarg^2}{\sigma^2}\right) + 
c_3 \sigma \sqrt{\frac{\maxccTarg}{n}}.
\end{equation}
Then $\estorc$ is a stationary point of Problem \eqref{prob:prox_fclsp} with probability $1 - \probNiceGradStatPt - \probSmallResidStatPt$, where
\begin{equation}
\probSmallResidStatPt 
= \NumNonIso \exp\left(-\frac{ n c_4 b_2^2 \tuneparam^2 }{8 \maxccTarg \sigma^2} \right) 
+ \nccNzTarg 4 \exp \left(- \frac{n c_5 b_2^2 \tuneparam^2 }{16 \sigma^2 }  \right)
\end{equation}
and
$$
\probNiceGradStatPt = 2 (D - \suppSizeTarg) \exp \left( - \frac{n a_0^2 \tuneparam^2}{2 \maxccTargPow{2} \sigma^2} \right).
$$
\end{theorem}

Condition \eqref{eq:bd_shrink__indep_unbiased__stat_pt__targ_gap__assumpt} can easily be weakened (e.g. allowing within block zeros), but we state it as above for the sake of exposition.\footnote{This assumption requires $\minSuppMagTarg > 0$ (i.e. the target graph should be fully connected) and  $||\esttarg||_{\text{max}}$ to be not too large. Both of these requirements can be significantly weakened by modifying \eqref{eq:fuzzy_sleepy_cat} in the proof of Theorem \ref{thm:bd_shrink__stat_pt__nice}.} 
In the setting of this theorem we see there exists a value of $\tuneparam$ such that $\esttarg$ is a stationary point with overwhelming probability if $\targgap \ge c \sigma \sqrt{\max\left(\maxccTarg \log \NumNonIso, \maxccTargPow{2} \log(D - \suppSizeTarg) \right)}$. 
\begin{remark}
Suppose the non-empty connected components of the target graph all have $\qTarg$ nodes; for this topology $\targgap \ge \qTarg \minSuppMagTarg$.
We will refer to this as the \textit{basic target graph} in the sequel.
\end{remark}
For the basic target graph, the spectral gap requirement translates to $ \minSuppMagTarg \ge c \sigma \sqrt{\log(D - \suppSizeTarg)}$.
This is the same condition for the hard thresholding estimator to succeed with overwhelming probability,
see Appendix \ref{a:proofs__bd_shrink}. 
Next consider the noisy within block edge setting where the elements of $\suppSetTarg$ have sub-Gaussian parameter $\omega >> \sigma$.
In this case we can check the target gap requirement becomes $\minSuppMagTarg \ge c \sqrt{  \max \left( \frac{\omega^2}{\qTarg}  \log \NumNonIso,  \sigma^2 \log(D - \suppSizeTarg)  \right)}$.
Hard-thresholding on the other hand requires $\minSuppMagTarg \ge c \sqrt{  \max \left( \omega^2 \log \NumNonIso,  \sigma^2 \log(D - \suppSizeTarg)  \right)}$. In this case $\omega$ can be up to $\sqrt{\qTarg}$ times larger than $\sigma$ before these two conditions are equal.

\subsection{Two step convergence}\label{ss:bd_shrink__two_step}

This section considers more general conditions on $\estok$ by relaxing Assumption \ref{assumpt:sub_g_indep_unbiased} e.g. we do not assume unbiasedness or independence.
We measure the quality of $\estok$ via an entrywise residual tail condition by borrowing the following definition from \citep{ravikumar2011high}.
\begin{definition} \label{def:tail_error}
Let $f: \mathbb{N} \times \mathbb{R}_+ \to \mathbb{R}_+$ be monotonically decreasing in both arguments and let $\tailconst \in [0, \infty)$.
We say an estimator $\estok \in \mathbb{R}^{D}$ satisfies a $\tailclass{f}{\tailconst}$ residual tail condition if for each $\ell \in [D]$,
\begin{equation} \label{eq:resid_tail}
\prob{ |\estok_{\ell} - \esttarg_{\ell}| > t } \le f(n, t), \qquad \text{for all } t \in (0, 1/\tailconst],
\end{equation}
where $1/0 := + \infty$ by convention (i.e. if $\tailconst = 0$ the bound is valid for all $t > 0$).
We assume $f(n, t) = 1$ for $t \in (0, 1/\tailconst]^C$.
\end{definition}
For example, in Lemma 1 of \citep{ravikumar2011high} the empirical covariance matrix from $n$ samples satisfies a $\tailclass{f}{\tailconst}$ residual tail condition for some $\tailconst$ and $f(n, t) = 4 \exp\left(- c n t^2 \right)$ for some constant $c$.

\begin{theorem} \label{thm:bd_shrink__two_step_general}
Assume $\estok$ is an estimator satisfying a $\tailclass{f}{\tailconst}$ residual tail condition and $\targgap \ge (b_1 + b_2) \tuneparam$.
Suppose $\residDiffGapRatio$ satisfies 
\begin{equation} \label{eq:thoe_ex__bd_shrink__reg_sat_whp__init_tol}
 \residDiffGapRatio  \le  \sqrt{\frac{b_1}{2^7 a_0 \maxccTarg} }
\end{equation}
and $ \max(\frac{a_1 \tuneparam}{4}, \frac{b_1 \tuneparam}{4 \maxccTarg}) \le \frac{1}{\tailconst}$.
Then the LLA algorithm initialized by $\estinit$ converges to $\estorc$ after two steps with probability at least $1 -\probGoodInitLLAOp - \probNiceGradLLAOp -\probSmallResidLLAOp$, where
\begin{equation}\label{thm:theo_ex_nbd_with_init_png}
\probNiceGradLLAOp :=   (D - \suppSizeTarg) f\left(n, \frac{a_1 \tuneparam}{4 \maxccTarg} \right),
\end{equation}
and
\begin{equation} \label{thm:theo_ex_nbd_with_init_pbg}
\probSmallResidLLAOp :=  \suppSizeTarg f\left(n, \frac{b_1 \tuneparam}{4 \maxccTarg} \right).
 \end{equation}

\end{theorem}

Next we turn our attention to the initializer.
For very high dimensional settings when $d >> n$, $\estok$ may not be a sufficiently good initializer.
Fortunately by thresholding $\estok$ we can obtain a good enough initializer for high-dimensional settings.
Note the following definition, borrowed from  \citep{rothman2009generalized}, encompasses many common thresholding operators including  hard-thresholding, soft-thresholding, etc.

\begin{definition} \label{def:gen_thresh}
A generalized thresholding operator is a function $\threshop{\cdot}{\ThreshTuneParam}: \mathbb{R} \to \mathbb{R}$ satisfying
$$
1.\; |\threshop{z}{\ThreshTuneParam}| \le |z| \qquad
2.\; \threshop{z}{\ThreshTuneParam} = 0, \text{ if  }|z| \le \tuneparam \qquad
3.\; |\threshop{z}{\ThreshTuneParam} - z| \le \tuneparam,
$$
for all $z \in \mathbb{R}$ and $\ThreshTuneParam \ge 0$.
\end{definition}

The following corollary puts Theorem  \ref{thm:bd_shrink__two_step_general} together with three choices of the initializer; $\estok$, $\threshop{\estok}{\ThreshTuneParam}$, and $0$.
Recall that if we initialize from 0, the first step of the LLA algorithm is equivalent to a Lasso penalty with parameter $a_0 \tuneparam$ (by Proposition \ref{prop:spect_pen_maj} and Definition \ref{def:scad_like_pen_func}).

\begin{corollary} \label{cor:bd_shrink__initializer__general} 
Suppose the assumptions of Theorem \ref{thm:bd_shrink__two_step_general} hold.
If the LLA algorithm is initialized by $\estok$ then
\begin{equation} \label{eq:bd_shrink_est_ok_init}
\probGoodInitLLAOp =  D f \left(n, \frac{\initRatioConsts \tuneparam}{2 d} \right).
\end{equation}

Let $\threshop{\cdot}{\ThreshTuneParam}$ be any generalized thresholding operator satisfying Definition \ref{def:gen_thresh}.
Assume  $\residDiffGapRatio$ satisfies \eqref{eq:thoe_ex__bd_shrink__reg_sat_whp__init_tol} and recall $\initRatioConsts$ from \eqref{eq:lla_probs}.
If
\begin{equation} \label{eq:thoe_ex__bd_shrink__initializer__param_init}
\ThreshTuneParam \le 
\frac{ \initRatioConsts \tuneparam}{12 \maxDegTarg}   
\end{equation}
and the LLA algorithm is initialized by $\threshop{\estok}{\ThreshTuneParam}$ then
\begin{equation} \label{eq:bd_shrink_gen_thresh_init}
\probGoodInitLLAOp =
2 \nnzSizeTarg f\left(n, \ThreshTuneParam  \right) + 2 D f\left(n, \frac{1}{2} \ThreshTuneParam \right).
\end{equation}

Furthermore, if 
$a_0 \le  \frac{ \initRatioConsts }{ 12 \maxDegTarg } $ and the LLA algorithm is initialized at 0 then
\begin{equation} \label{eq:thoe_ex__bd_shrink__initializer__prob_init_inot_from_one}
\probGoodInitLLAOp =  2 \nnzSizeTarg f\left(n, \tuneparam  \right) + 2D f\left(n, \frac{1}{2} \tuneparam  \right).
\end{equation}
\end{corollary}

A similar statement using $\probGoodInitLLAFrob$ is possible (see the proof of Corollary \ref{cor:bd_shrink__initializer__general}).
In the independent, sub-Gaussian setting of the previous section we can improve upon the direct applications of Theorem \ref{thm:bd_shrink__two_step_general} and Corollary \ref{cor:bd_shrink__initializer__general}.
\begin{remark} \label{rem:two_step_reates_under_indep_unbaised}
If $\estok$ satisfies Assumption \ref{assumpt:sub_g_indep_unbiased} and $8 b_1 \tuneparam$ is at least the right hand size of \eqref{eq:bd_shrink__indep_unbiased__stat_pt__targ_gap__assumpt} then
$$
\probSmallResidLLAOp
=  \NumNonIso \exp\left(-\frac{ n  c_1 b_1^2 \tuneparam^2 }{\maxccTarg  \sigma^2} \right) + 4  \nccNzTarg \exp \left(- \frac{n c_2 b_1^2 \tuneparam^2}{ \sigma^2 } \right).
$$
Similarly, if the LLA algorithm is initialized by $\estok$ and  $8 \initRatioConsts \tuneparam$ is at least \eqref{eq:bd_shrink__indep_unbiased__stat_pt__targ_gap__assumpt} then
$$
\probGoodInitLLAOp 
=  \NumNonIso \exp\left(-\frac{ n c_3 \initRatioConsts^2 \tuneparam^2 }{\maxccTarg \sigma^2} \right) + 4  \nccNzTarg \exp \left(- \frac{n c_4 \initRatioConsts^2 \tuneparam^2}{ \sigma^2 } \right) 
+  2 d \exp \left( - \frac{n c_5 \initRatioConsts^2 \tuneparam^2}{ d \sigma^2} \right).
$$
These two claims follow from technical lemmas provided in the appendix.
\end{remark}

Next we put these results into context by working under the independent, sub-Gaussian setting of Assumption \ref{assumpt:sub_g_indep_unbiased}.  
In this case, if we use $\estok$ as the initializer we need  $\targgap \ge  c \sigma  \sqrt{\maxccTarg d \log d}$
for there to be a value of $\tuneparam$ such that the LLA algorithm converges to the oracle in two steps with overwhelming probability.
On the other hand if we initialize with an appropriately chosen thresholded initializer we need 
$\targgap \ge  c \sigma  \sqrt{\maxccTarg \maxDegTargPow{2}  \log d }$.
Without the independence/unbiasedness of Assumption \ref{assumpt:sub_g_indep_unbiased} these become slightly worse. 
If we assume the basic target graph topology from the previous section then these conditions translate to  
$\minSuppMagTarg \ge  c \sigma  \sqrt{\frac{d \log d}{\qTarg} }$
and
$\minSuppMagTarg \ge  c \sigma  \sqrt{\qTarg\log d }$
respectively.\footnote{For this graph topology note $\frac{d}{\qTarg} = \nccNzTarg + \frac{\NumIso}{\qTarg}$ is linear in $\nccNzTarg$, not in $\qTarg$.} 
These rates are worse than those required for stationary point conditions and hard thresholding estimator discussed in the previous section.
The empirical evidence from Section \ref{s:sim} suggests these rates can be improved.

\section{Block sparse regression} \label{s:regression_examples}

This section studies the \penabr\ penalty LLA algorithm for linear and logistic regression.
As in the previous section we assume Assumption \ref{assu:basic} is satisfied throughout this section.
The setup in this section closely follows the examples of \cite{fan2014strong}.

\subsection{Linear regression} \label{ss:theo_ex_lin_reg}

Consider the least squares objective function with the \penabr\ penalty, 
\begin{equation}  \label{prob:lin_reg_fclsp}
\begin{aligned}
& \underset{\est \in \mathbb{R}^{{d \choose 2}}}{\textup{minimize}}  & & \frac{1}{2n} ||y - X \est||_2^2 + \fcls{\tuneparam}{\est}\\ 
\end{aligned}
\end{equation}
where $y \in \mathbb{R}^{n}$ is the response vector and $X \in \mathbb{R}^{n \times D}$ is the covariate matrix.

Let $\esttarg \in \mathbb{R}^D$ be the true parameter vector in the linear regression model $y = X \esttarg + \epsilon$. 
We assume then entries of $\epsilon \in \mathbb{R}^n$ are $i.i.d$ mean zero $\sigma$ sub-Gaussian.
For linear regression the oracle estimator has a closed form,
$$
\estorc = (\estorc_{\suppSetTarg}, 0), \text{ where } \estorc_{\suppSetTarg} = (X^T_{\suppSetTarg} X_{\suppSetTarg})^{-1} X^T_{\suppSetTarg} y.
$$

We will need the following quantities: $\minevalxsshort:=  \minevalxsexpl$, $\maxevalxshort :=\maxevalxexpl$, $\maxColNorm :=  \frac{1}{n}\max_{j \in [D]} ||X_j||_2^2$ (which is often normalized to be 1). 

\begin{theorem} \label{thm:thoe_ex__lin_reg__reg_sat_whp}
Assume  $\targgap \ge (b_1 + b_2) \tuneparam$ and let
\begin{equation} \label{eq:theo_ex_lin_reg_init_tol}
\residDiffGapRatio \le
\sqrt{\frac{a_1 \minevalxsshort }{a_0  \sqrt{2^{11} \maxccTarg \maxColNorm \maxevalxshort}}}
\wedge
\sqrt{\frac{b_1\minevalxsshort }{64 a_0 \maxccTarg}}.
\end{equation}
Then the LLA algorithm initialized by $\estinit$ converges to $\estorc$ in two steps with probability at least $1 -\probGoodInitLLAFrob - \probNiceGradLLAFrob - \probSmallResidLLAFrob$, where
\begin{equation}\label{eq:thoe_ex__lin_reg__reg_sat_whp__prob_nice_grad}
\probNiceGradLLAFrob =  2 (D - \suppSizeTarg) \exp\left(\frac{-n a_1^2 \tuneparam^2}{32 \maxccTargPow{2} \sigma^2 \maxColNorm} \right)
\end{equation}
and
\begin{equation} \label{eq:thoe_ex__lin_reg__reg_sat_whp__prob_small_resid}
\probSmallResidLLAFrob = 
2 \suppSizeTarg \exp \left(  \frac{-n \minevalxsshort b_1^2 \tuneparam^2 }{32 \maxccTargPow{2} \sigma^2}\right).
\end{equation}
\end{theorem}

We complement Theorem \ref{thm:thoe_ex__lin_reg__reg_sat_whp} by showing that the LLA algorithm can be initialized by a Lasso estimate as in  \cite{fan2014strong}.
We use Theorem 7.13 of \cite{wainwright2019high} to bound the Lasso error under the assumption that $X$ satisfies the $(\RestrEvalLinReg, 3)$ \textit{restricted eigenvalue condition} stated in (7.22) of \cite{wainwright2019high}.
Under the sub-Gaussian error assumptions of this section the Lasso solution,  $\widehat{\est}^{\text{lasso}}$, with penalty parameter $ \ThreshTuneParam$ satisfies
\begin{equation*}
||\widehat{\est}^{\text{lasso}} - \esttarg||_1 \le \frac{12  \nnzSizeTarg \ThreshTuneParam  }{\RestrEvalLinReg}  
\end{equation*}
with probability at least $1 - 2 \exp(- \frac{ n  \ThreshTuneParam^2}{8 \sigma^2 \maxColNorm^2} )$. 

\begin{corollary} \label{cor:theo_ex_lin_reg_lasso_init}
Suppose the assumptions of Theorem \ref{thm:thoe_ex__lin_reg__reg_sat_whp} hold and
$$
\frac{12  \nnzSizeTarg\ThreshTuneParam  }{\RestrEvalLinReg}  
\le  \frac{\initRatioConsts \tuneparam}{4},
$$
where $\residDiffGapRatio$ is given in \eqref{eq:theo_ex_lin_reg_init_tol} and $\initRatioConsts$ is given in \eqref{eq:lla_probs}.
Then the two step convergence claim of Theorem \ref{thm:thoe_ex__lin_reg__reg_sat_whp} holds when the LLA algorithm is initialized by $\widehat{\est}^{\text{lasso}}$  where
$$
\probGoodInitLLAFrob = 2 \exp(- \frac{ n \ThreshTuneParam^2}{8 \sigma^2 \maxColNorm^2} ).
$$
Furthermore, if $\frac{12  \nnzSizeTarg a_0  }{\RestrEvalLinReg} \le  \frac{\initRatioConsts }{4} $ then the LLA algorithm initialized by 0 converges to $\estorc$ after three steps with probability at least $1 - 2 \exp(- \frac{ n  \tuneparam^2}{8 \sigma^2 \maxColNorm^2} ) - \probNiceGradLLAFrob - \probSmallResidLLAFrob$.
\end{corollary}

Assume the target graph has the basic topology discussed in the previous section.
In this case there exists a value of $\tuneparam$ such that the LLA algorithm initialized from an appropriately selected Lasso initializer converges to $\esttarg$ in two steps with overwhelming probability if $\minSuppMagTarg \ge c \sigma \nccNzTarg \qTargPow{3/2}  \sqrt{ \log d}$.
Section \ref{s:sim} suggests this rate may be improvable.

%
%
%
%
%
%
%
%
%
%

\subsection{Logistic regression} \label{ss:theo_ex_logisitc}

This section considers logistic regression when the true coefficient vector is block sparse.
Suppose we observe a covariate matrix $X \in \mathbb{R}^{n \times D}$ and a random binary response vector $y \in \{0, 1\}^n$.
We assume the standard logistic regression setup 
\begin{equation} \label{eq:logistic_setup}
y_i | X(i, :) \sim \text{Bernoulli} \left( \frac{e^{X(i, :)^T \esttarg }}{1 + e^{X(i, :)^T \esttarg }} \right), i=1, \dots, n \text{, independently},
\end{equation}
where $\esttarg \in \mathbb{R}^D$ is the true coefficient vector.
The \penabr\ penalized logistic regression problem is then given by
\begin{equation}  \label{prob:logistic_reg_fclsp}
\begin{aligned}
& \underset{\est \in \mathbb{R}^{D}}{\textup{minimize}}  & & \frac{1}{n} \sum_{i=1}  \left(-y_i X(i, :)^T \est + \LogisitcLink{X(i, :)^T \est}  \right)+ \fcls{\tuneparam}{\est},
\end{aligned}
\end{equation}
where $\LogisitcLink{t} = \log(1 + e^t)$ is the canonical link function.

To ease notation we make a few definitions.
Let $\LogisticMaxVarNorm := \max_{j \in [D]} \frac{ ||X_{j}||_2^2}{n}$ and $\LogisticMaxX := ||X||_{\text{max}}$.
Let  $\LogisitcH{\cdot}: \mathbb{R}^{D} \to \mathbb{R}^n$ be the function given by  $\LogisitcH{\est}_i = \LogisitcLinkGG{X(i, :)^T \est}, i = 1, \dots, n$.
Let $\LogisticEvalMax :=  \max_{j \in [D] } \lambda_{\text{max}} \left(\frac{1}{n} X_{\suppSetTarg}^T \text{diag}(|X_j|) X_{\suppSetTarg} \right)$, 
$\LogisticXsHXscOne :=   || \frac{1}{n}    X_{\suppSetTargPow{C}}^T \text{diag}(\LogisitcH{\esttarg})) X_{\suppSetTarg}  ||_{\text{max}}$,
and \\
$\LogisticXshXsInv := \opnorm{ \left( \frac{1}{n} X_{\suppSetTarg}^T \text{diag}(\LogisitcHS{\esttarg}) X_{\suppSetTarg} \right)^{-1}}{\infty}{\infty} $.

%
\begin{theorem} \label{thm:thoe_ex__logistic__reg_sat_whp}
Assume  $\targgap \ge (b_1 + b_2) \tuneparam$, the columns of $X_{\suppSetTarg}$ are linearly independent and let
\begin{equation} \label{eq:theo_ex_logisitc_init_tol}
\residDiffGapRatio \le
\sqrt{ \frac{1}{2^7 a_0 \LogisticXshXsInv}  \min\left(
\frac{b_1}{2 \maxccTarg},
\frac{a_1}{ 8 \maxccTarg \left( \suppSizeTarg \LogisticXsHXscOne + \frac{\sqrt{\suppSizeTarg} \LogisticEvalMax}{8} \right)}
\right)} 
\end{equation}
and assume
\begin{equation} \label{eq:theo_ex_logistic_tune_param_ubd}
\tuneparam \le 
\frac{4}{\suppSizeTarg \LogisticEvalMax \LogisticXshXsInv }  
\min \left(
\frac{\maxccTarg}{2b_1},
\frac{8 \maxccTarg \left(\suppSizeTarg \LogisticXsHXscOne + \frac{\sqrt{\suppSizeTarg} \LogisticEvalMax}{8} \right)}{a_1}
\right).
\end{equation}

Then the LLA algorithm initialized by $\estinit$ converges to $\estorc$ in two steps with probability at least $1 -\probGoodInitLLAOp - \probNiceGradLLAOp  - \probSmallResidLLAOp$, where
\begin{equation}\label{eq:theo_ex_logistic__reg_sat__prob_nice_grad}
\probNiceGradLLAOp =  
2 \suppSizeTarg \exp\left(- \frac{n}{\LogisticMaxVarNorm} \left(\frac{a_1}{ 8 \maxccTarg \left(\suppSizeTarg \LogisticXsHXscOne + \frac{\sqrt{\supptarg} \LogisticEvalMax}{8} \right)}\right)^2 \tuneparam^2 \right)
+ 
(D - \supptarg) \exp\left(\frac{-na_1^2 \tuneparam^2}{64 \maxccTargPow{2} \LogisticMaxVarNorm} \right)
\end{equation}
and
\begin{equation} \label{eq:theo_ex_logistic__reg_sat__prob_small_resid}
\probSmallResidLLAOp =  
2 \supptarg \exp\left(- \frac{n}{\LogisticMaxVarNorm} \frac{b_1^2 \tuneparam^2}{4\maxccTargPow{2}} \right)
\end{equation}

\end{theorem}

As in the previous section, Lasso penalized logistic regression gives a ``good enough" initial estimator assuming an appropriate restricted eigenvalue condition.
In the following theorem $\widehat{\est}^{\text{lasso}}$ is the Lasso solution with tuning parameter $\ThreshTuneParam$ and we assume a restricted eigenvalue condition \eqref{eq:logistic_restr_eval} is satisfied with parameter $\LogisticRestrEval$.
See Appendix \ref{as:logistic_reg__prelim} for details. 

\begin{corollary} \label{cor:theo_ex_logisitc_lasso_init}
Suppose the assumptions of Theorem \ref{thm:thoe_ex__logistic__reg_sat_whp} hold and
$$
\frac{20 \nnzSizeTarg}{\LogisticRestrEval}\ThreshTuneParam 
\le  \frac{\initRatioConsts \tuneparam}{4} 
$$
where $\residDiffGapRatio$ is given in \eqref{eq:theo_ex_logisitc_init_tol}.
Then the two step convergence claim of Theorem \ref{thm:thoe_ex__logistic__reg_sat_whp} holds when the LLA algorithm is initialized by $\widehat{\est}^{\text{lasso}}$  where
$$
\probNiceGradLLAOp = \probNiceGradLLAFrob = 2 D \exp(- \frac{ n}{2\LogisticMaxVarNorm } \ThreshTuneParam^2 ).
$$

Furthermore if
$
\frac{20 \nnzSizeTarg a_0 }{\LogisticRestrEval} 
\le  \frac{\initRatioConsts }{4} 
$
then the LLA algorithm initialized by 0 converges to $\estorc$ after three steps with probability at least $1 - 2 D \exp(- \frac{ n}{2\LogisticMaxVarNorm }  \tuneparam^2 ) - \probNiceGradLLAOp - \probSmallResidLLAOp$.
\end{corollary}


%
%
%
%
%
%
%
%


\section{Simulations} \label{s:sim}

This section studies several \penabr\ penalized models discussed in previous sections through simulations.
The block diagonal shrinkage algorithm from Section \ref{s:bd_shrink} is examined for a Gaussian sequence model (Section \ref{ss:sim__gaussian_seq}) and for covariance estimation (Section \ref{ss:sim__covar}).
The block sparse linear and logistic regression models from Section \ref{s:regression_examples} are examined in Sections \ref{ss:sim__lin_reg} and \ref{ss:sim__log_reg} respectively.

For each model we look at three target graph topologies with two fully connected components consisting of: 5 , 10,  and 25 nodes.
In other words $\esttarg$ has $D= {2 \cdot 5 \choose 2} = 45$ dimensions for the first graph, ${2 \cdot 10 \choose 2} = 190$ for the second, and $ { 2 \cdot  25 \choose 2} = 1225$ for the third.\footnote{For these graphs $\esttarg$ has $2 \times {5 \choose 2} = 20$, $2 \times {10 \choose 2} = 90$, and $ 2 \times  {25 \choose 2} = 600$ non-zero elements respectively.} 

This section focuses on how close the \penabr\ penalized estimator -- and its entrywise sparse competitors -- comes the block oracle, $\estorc$.
For each model and graph topology the competing estimators are fit for a range of number of samples and each experiment is repeated $20$ times.
Throughout this section $g_{\tuneparam}$ is set to be the SCAD function\footnote{We found empirically that smaller values of $a$ worked better for the \penabr\ penalty contrary to the the standard default of $a=3.7$.} with $a=2.1$. 

Each of the algorithms in this section require selecting a tuning parameter e.g. $\tuneparam$ for the \penabr\ penalty.
To simplify the comparison of competing models we  ``cheat" by picking the best possible tuning parameter for the models presented below.
In other words, we fit each model for a range of tuning parameter values then pick the tuning parameter that gives the best performance (i.e. the one whose estimate is nearest the block oracle).
The tuning parameter range for the \penabr\ models is specified using the discussion in Appendix \ref{a:tune_param_ubd}. 
Choosing a good initializer is important for obtaining good performance.
For each model below we tune the initializers using 10-fold cross-validitaion i.e. we do not ``cheat" for the initializers.
For tuning the \penabr\ penalty in practice we suggest cross-validation.

The code to reproduce the results is provided in the github repository linked to in Section \ref{s:intro}.
This software makes use of several standard python libraries \citep{matplotlib2007hunter, scikit2011pedrogosa, harris2020array, scipy2020virtanen, seaborn2021waskom}, an accelerated coordinate descent package for linear regression \citep{massias2020dual, bertrand2021anderson} and a recently developed package for penalized GLMs \citep{carmichael2021yet}.

\subsection{Gaussian sequence model} \label{ss:sim__gaussian_seq} 

This section examines the block diagonal shrinkage algorithm from Section \ref{s:bd_shrink} in the context of a Gaussian sequence model.
For this model we observe $x_1, \dots, x_n \in \mathbb{R}^{D}$, $x_i = \esttarg + \sigma \epsilon_i$ where each $\epsilon_i \sim N(0, I_D)$ independently.
Each non-zero entry of $\esttarg$ is set to 1 and we set\footnote{This value is selected such that $||\estorc - \esttarg||_2 \approx ||\esttarg||_2$ with $n=10$ samples, which indicates we are in a very challenging signal to noise regime.} 
$\sigma = 3.25$.
For the objective of Problem \eqref{prob:prox_fclsp} we use $\estok = \overline{x}$.

We examine the LLA algorithm with three different initialization strategies:  1) hard-thresholding tuned with 10 fold cross-validation 2) from the empirical means, $\estok$, and 3) from 0.
We also examine the LLA after two steps (three for the 0 initialization) and after convergence to a stationary point.
As a baseline for comparison we contrast the block diagonal shrinkage estimate to the entrywise hard-thresholding estimate and the empirical means.
Recall the \penabr\ and baseline hard-thresholded estimates are tuned via ``cheating" by selecting the tuning parameter whose estimate is nearest to $\estorc$.

\begin{figure}[H]
 \centering
\begin{subfigure}[t]{0.3\textwidth}
\centering
\includegraphics[width=1.1\linewidth, height=1.5\linewidth]{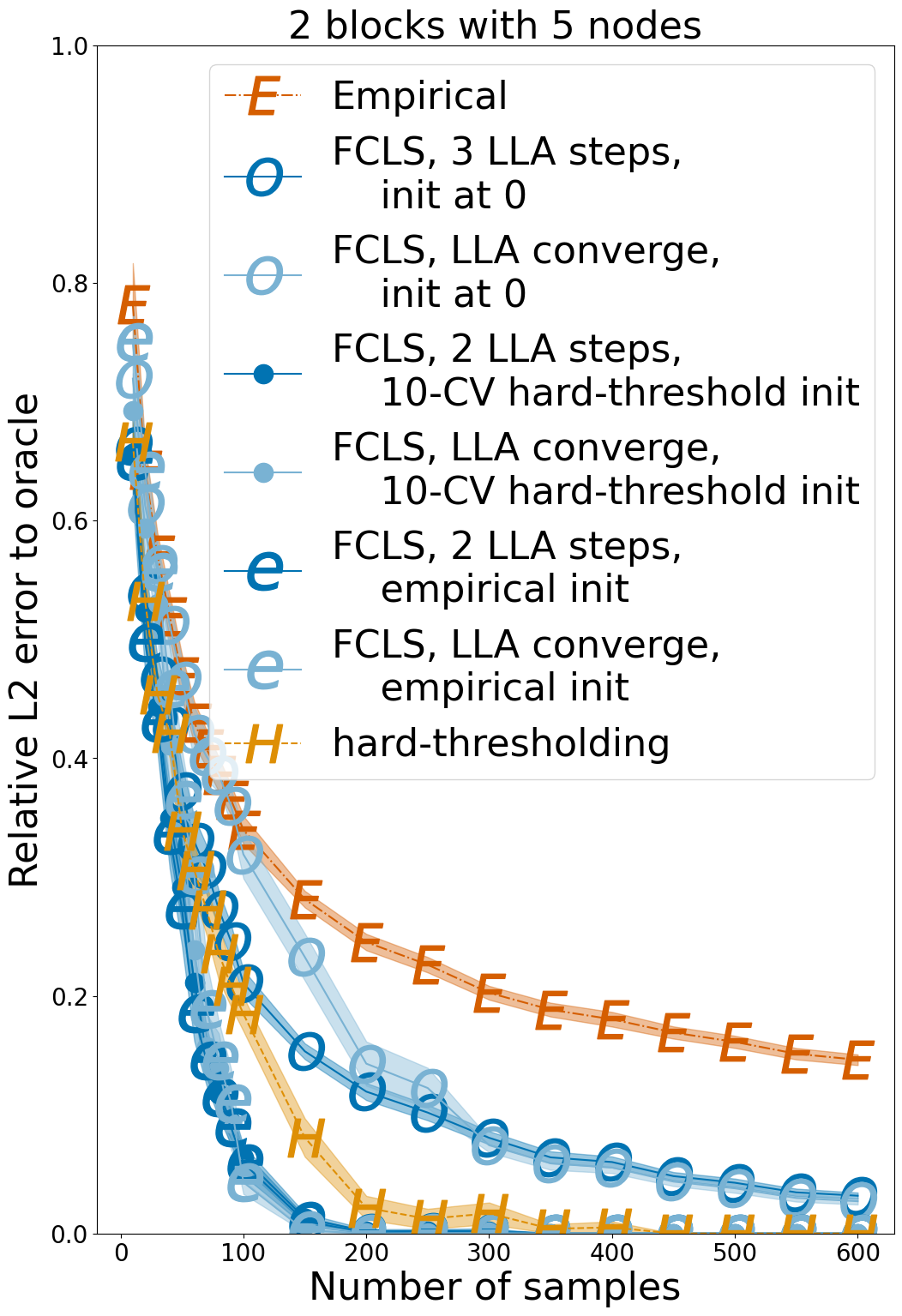}
\caption{
Two $5 \times 5$ blocks
}
\label{fig:means_est__bsize=5_2__vs__oracle__L2_rel}
\end{subfigure}
\hfill
\begin{subfigure}[t]{0.3\textwidth}
\centering
\includegraphics[width=1.1\linewidth, height=1.5\linewidth]{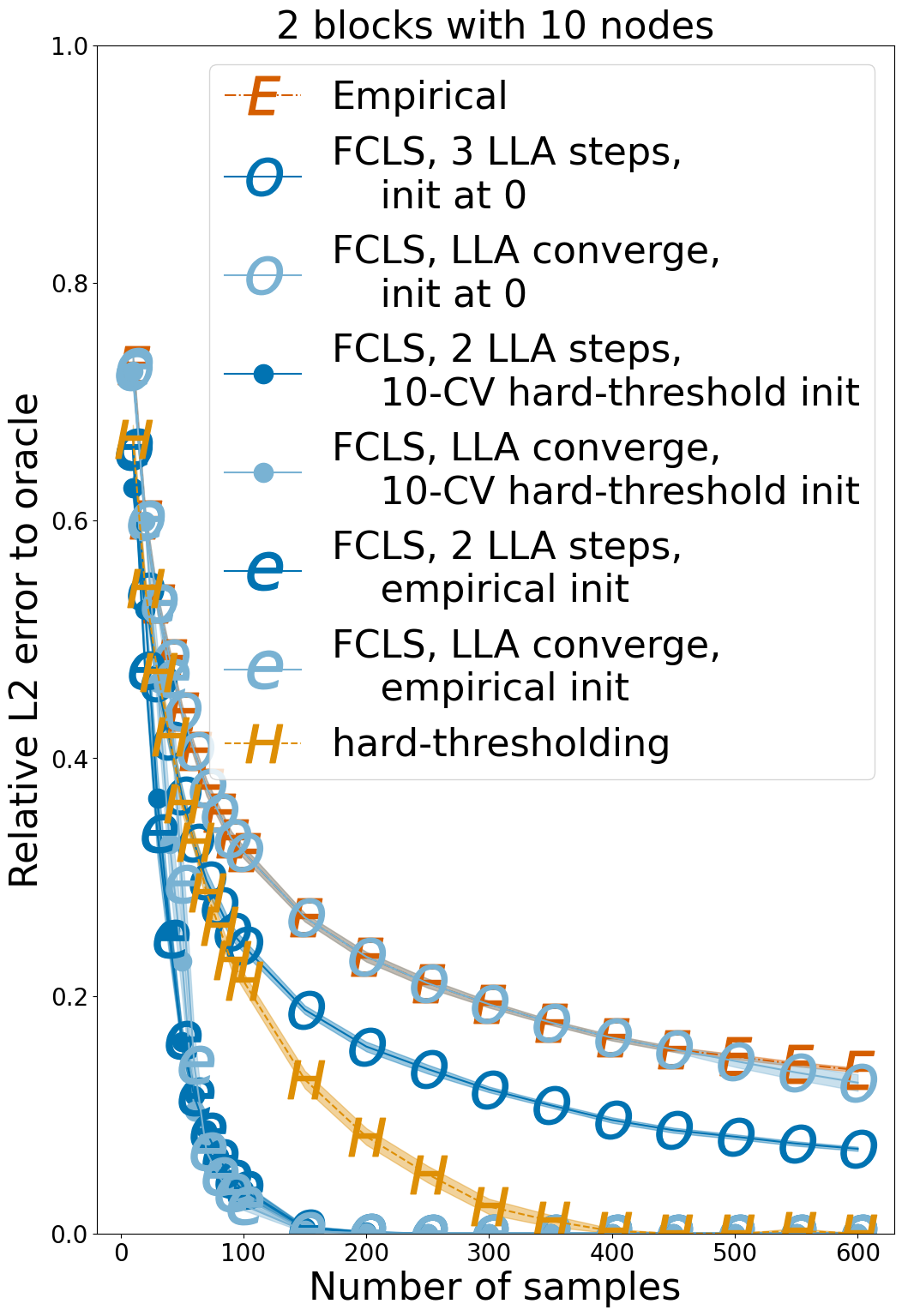}
\caption{
Two $10 \times 10$ blocks
}
\label{fig:means_est__bsize=10_2__vs__oracle__L2_rel}
\end{subfigure}
\hfill
\begin{subfigure}[t]{0.3\textwidth}
\centering
\includegraphics[width=1.1\linewidth, height=1.5\linewidth]{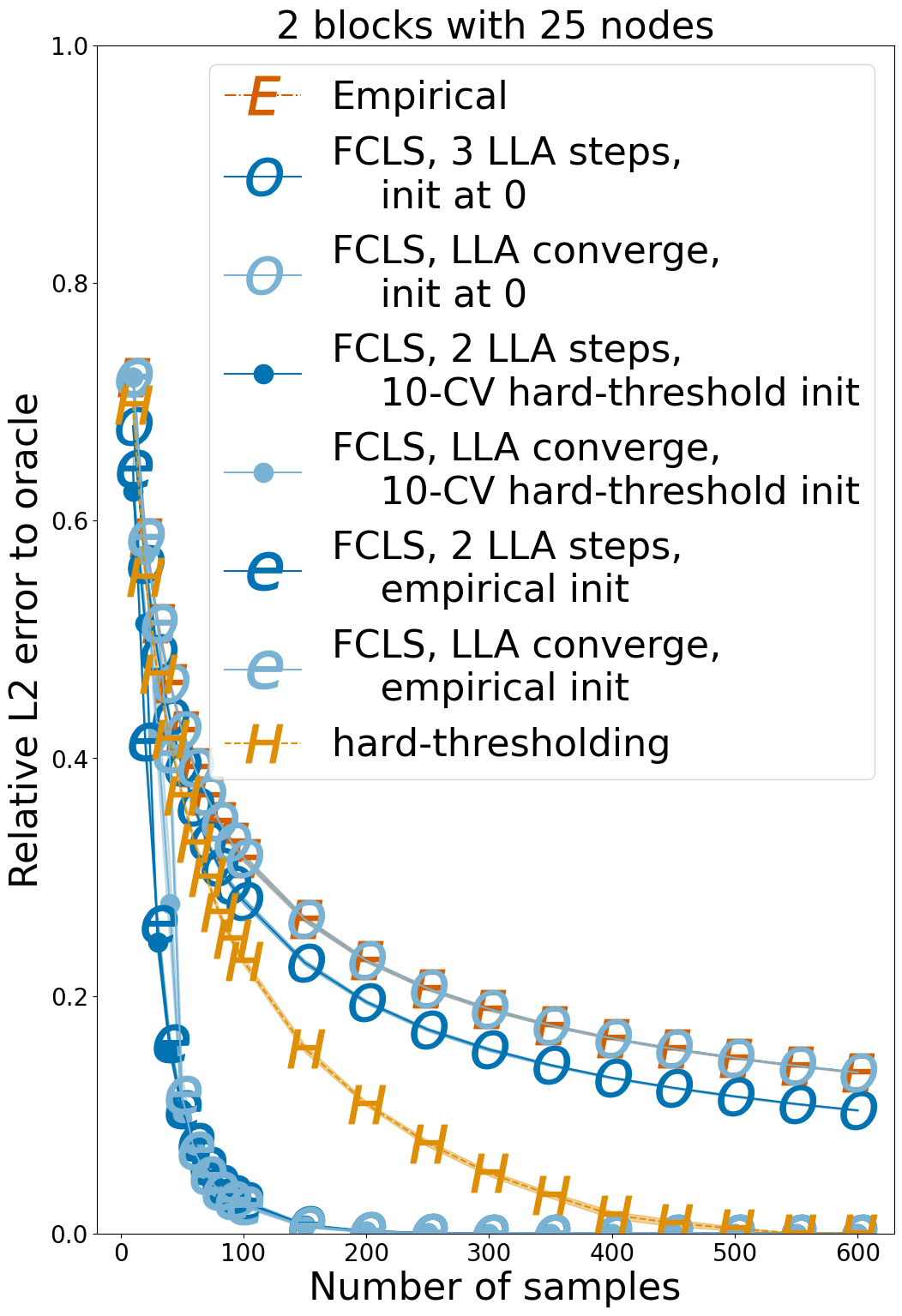}
\caption{
Two $25 \times 25$ blocks
}
\label{fig:means_est__bsize=25_2__vs__oracle__L2_rel}
\end{subfigure}
\caption{
Results for Gaussian sequence model.
The $x$-axis shows the number of samples and the $y$-axis shows the relative $L_2$ norm error of the estimated coefficient to the oracle coefficient i.e. $||\widehat{\est} - \estorc||_2 / ||\estorc||_2$. 
The shaded area shows the standard error after 20 Monte-Carlo repetitions.
The \penabr\ penalized models are shown in blue with either a dot (hard-threshold initialization), 0 (0 initialization) or e (empirical initialization).
}
\label{fig:means_est__oracle__L2_rel}
\end{figure}

Figure \ref{fig:means_est__oracle__L2_rel} examines the models across a range of number of samples for each of the three graph topologies.
This figure shows the FCLS penalty initialized from either 10-CV hard-thresholding or from the empirical means performs the best.
The performance gap between block diagonal shrinkage and hard-thresholding appears to widen for target graphs with larger connected components.

Figure \ref{fig:means_est__oracle__L2_rel} shows 10-CV hard-thresholding and empirical means initialization strategies perform similarly while the 0 initialization strategy performs worse.
For the first two initializers, running the LLA algorithm to convergence does not provide much benefit over just two steps. 
These observations hold for the models in the following sections so we omit them from the following figures.

\subsection{Covariance estimation} \label{ss:sim__covar} 

This section examines the block diagonal shrinkage algorithm from Section \ref{s:bd_shrink} in the context of covariance estimation.
In this case we observe  $x_1, \dots, x_n \in \mathbb{R}^{D}$ independently where  $x_i \sim N(0, \Sigma^*)$ and $\Sigma^* = I_{D} +\mathcal{A}(\esttarg)$. 
The non-zero entries of $\esttarg$ are set to 0.3. 
Here $\estok$ is set to the upper triangular elements of the empirical covariance matrix $\widehat{\Sigma}$. 
As in the previous section the LLA algorithm is initialized with a hard-thresholding operation tuned with 10-fold cross-validation.

\begin{figure}[H]
 \centering
\begin{subfigure}[t]{0.3\textwidth}
\centering
\includegraphics[width=\linewidth]{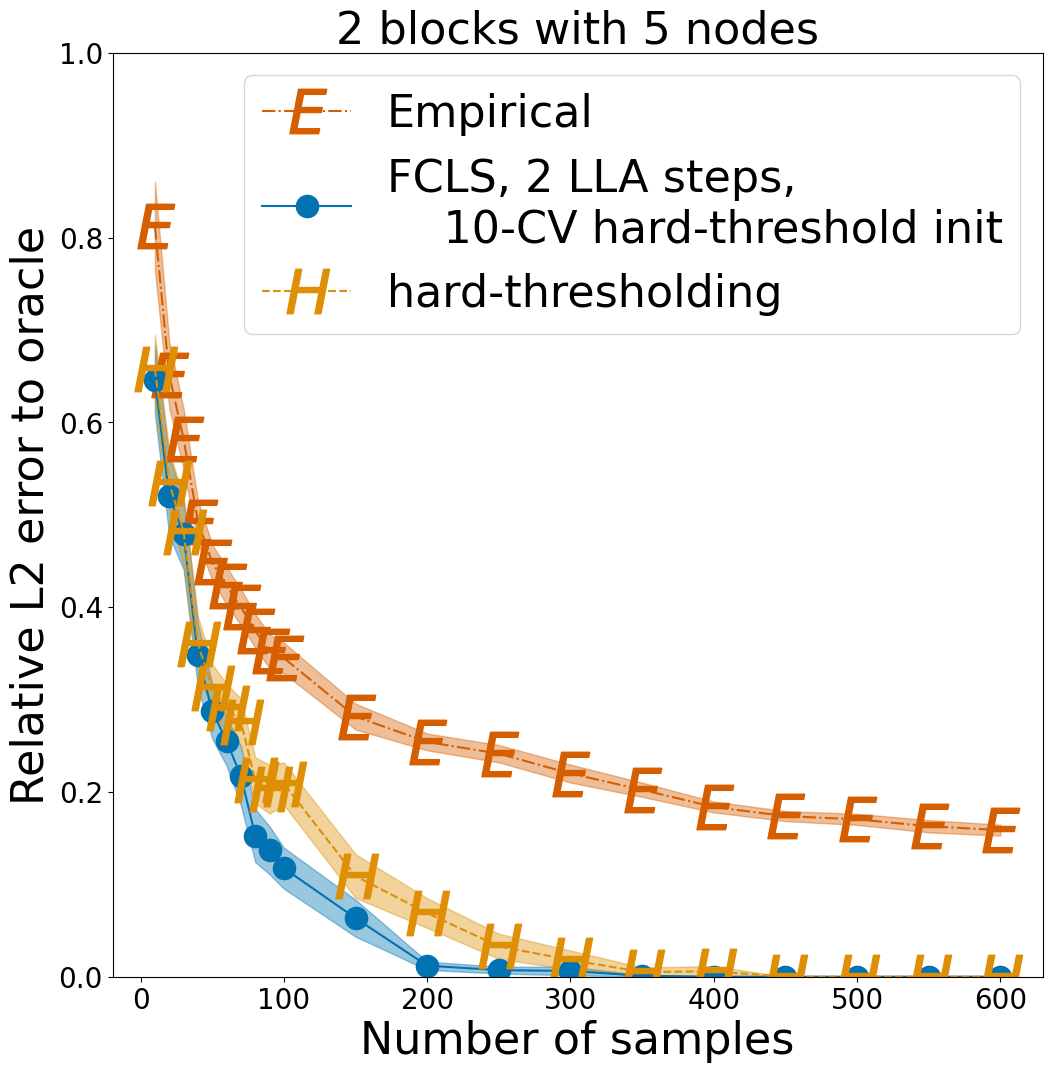}
\caption{
Two $5 \times 5$ blocks
}
\label{fig:covar__bsize=5_2__vs__oracle__L2_rel}
\end{subfigure}
\hfill
\begin{subfigure}[t]{0.3\textwidth}
\centering
\includegraphics[width=\linewidth]{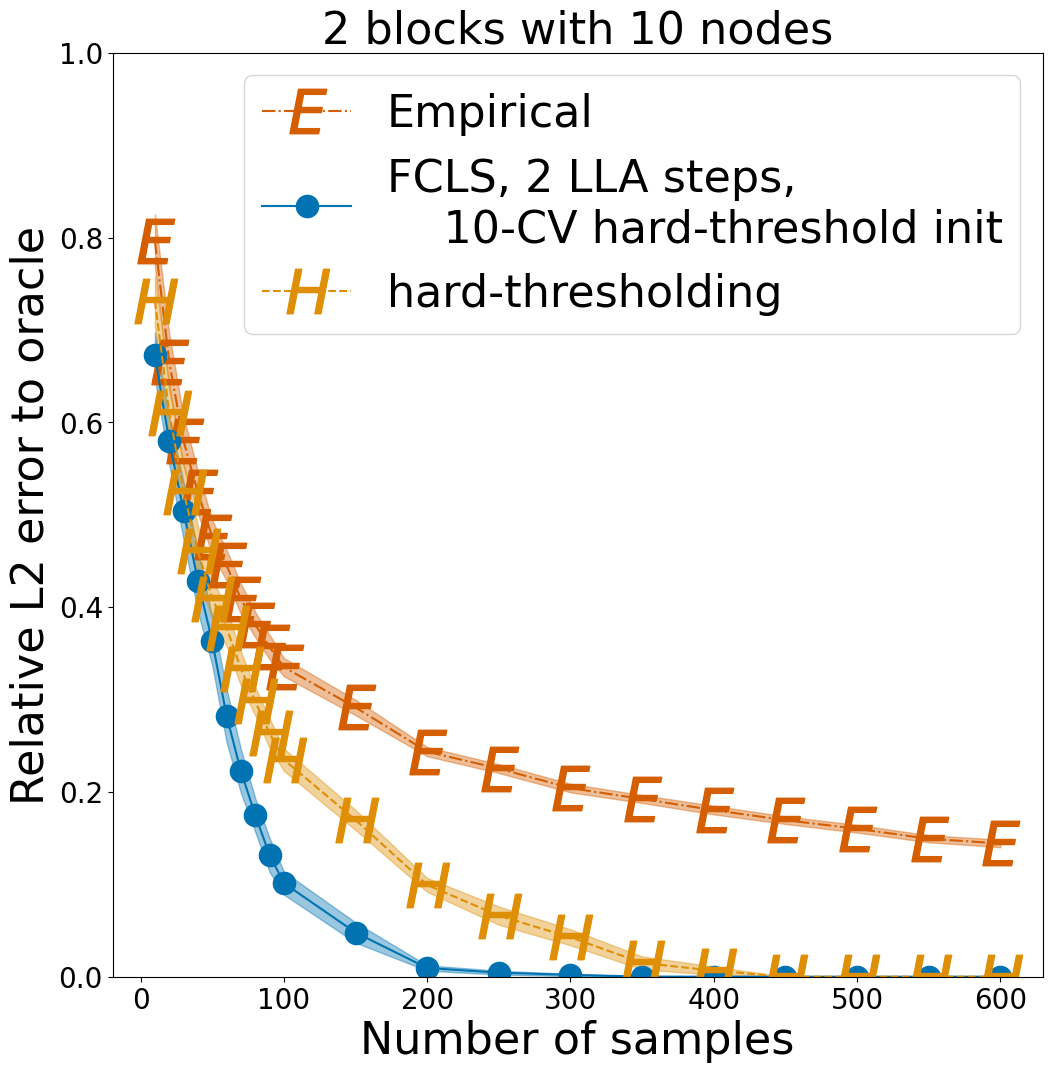}
\caption{
Two $10 \times 10$ blocks
}
\label{fig:covar__bsize=10_2__vs__oracle__L2_rel}
\end{subfigure}
\hfill
\begin{subfigure}[t]{0.3\textwidth}
\centering
\includegraphics[width=\linewidth]{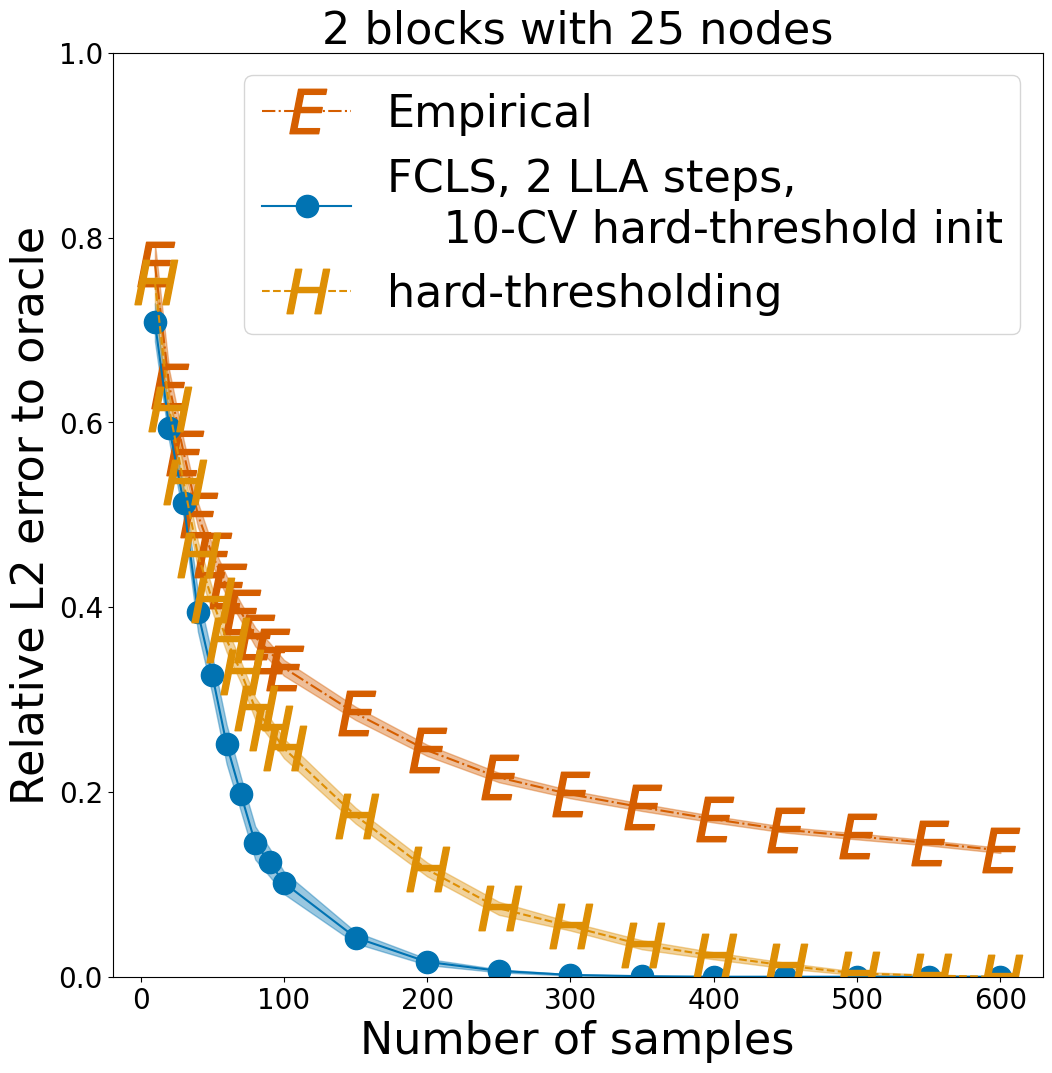}
\caption{
Two $25 \times 25$ blocks
}
\label{fig:covar__bsize=25_2__vs__oracle__L2_rel}
\end{subfigure}
\caption{
Results for covariance estimation.
The \penabr\ penalized model is shown with blue dots.
}
\label{fig:covar__oracle__L2_rel}
\end{figure}

Figure \ref{fig:covar__oracle__L2_rel} examines the models across a range of number of samples for each of the three graph topologies.
As in the previous section the block diagonal shrinkage estimate has the best performance.
Similarly, the gap between block diagonal shrinkage and hard-thresholding appears to increase with larger block sizes.

\subsection{Linear regression} \label{ss:sim__lin_reg} 

This section examines the block sparse linear regression example from Section \ref{ss:theo_ex_lin_reg}.
In this case we observe  $(y_1, x_1), \dots, (y_i, x_n)$ where $x_i \sim N(0, I_D)$ and $y_i = x_i^T \esttarg + \epsilon_i$ for $\epsilon_i \sim N(0, 1)$ independently.
The non-zero entries of $\esttarg$ are randomly set to $\pm 1$.

The LLA algorithm for the \penabr\ penalized model is initialized using a Lasso estimator that is tuned via 10-fold cross-validation.
As a baseline for comparison we compute the Lasso estimate and the entrywise SCAD penalized estimator with $a=3.7$ following \citep{fan2001variable}. 
This entrywise SCAD penalty is fit with one step of the LLA algorithm as in \cite{fan2014strong} and initialized with the same Lasso estimate.
Recall both the Lasso and entrywise SCAD estimators are tuning via ``cheating".

\begin{figure}[H]
 \centering
\begin{subfigure}[t]{0.3\textwidth}
\centering
\includegraphics[width=\linewidth]{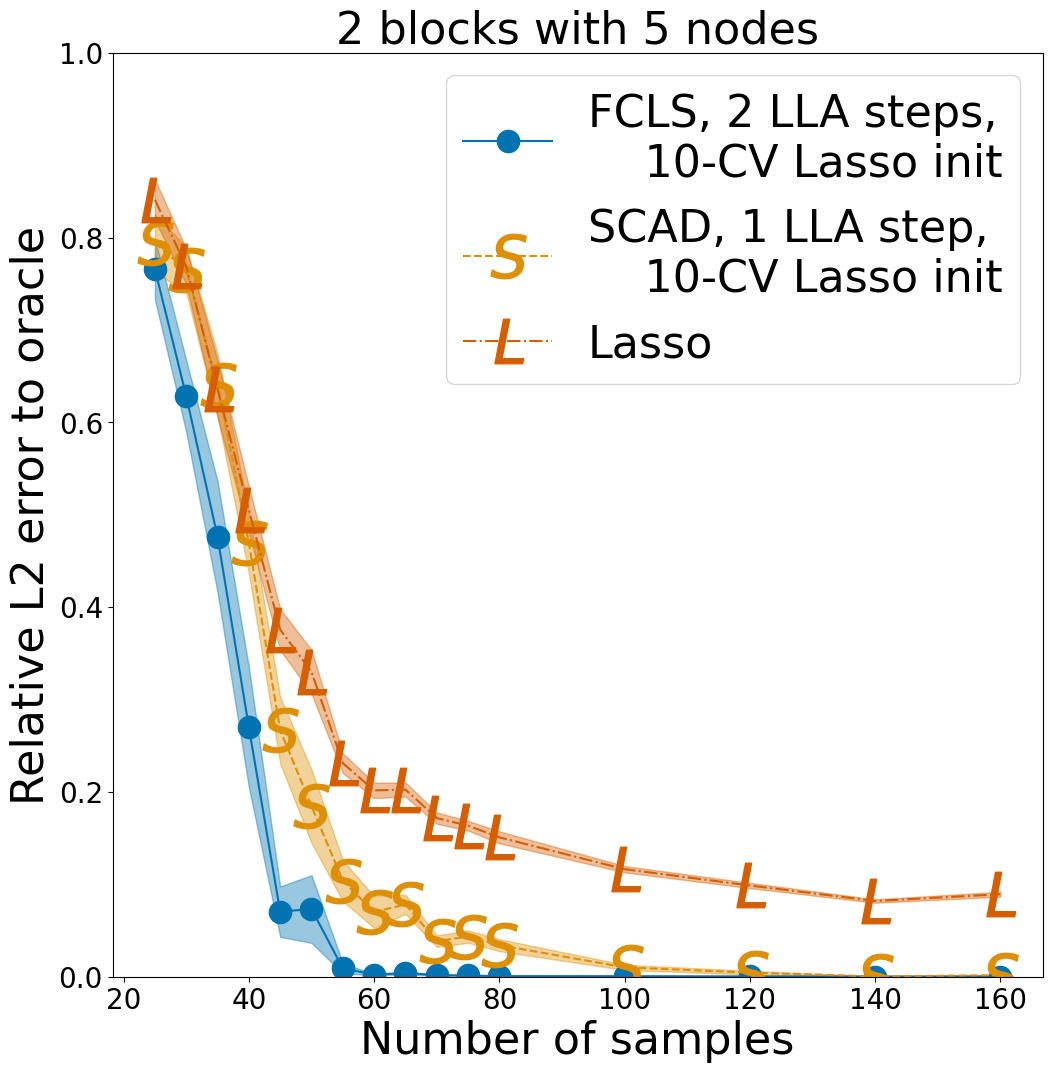}
\caption{
Two $5 \times 5$ blocks
}
\label{fig:lin_reg__bsize=5_2__vs__oracle__L2_rel}
\end{subfigure}
\hfill
\begin{subfigure}[t]{0.3\textwidth}
\centering
\includegraphics[width=\linewidth]{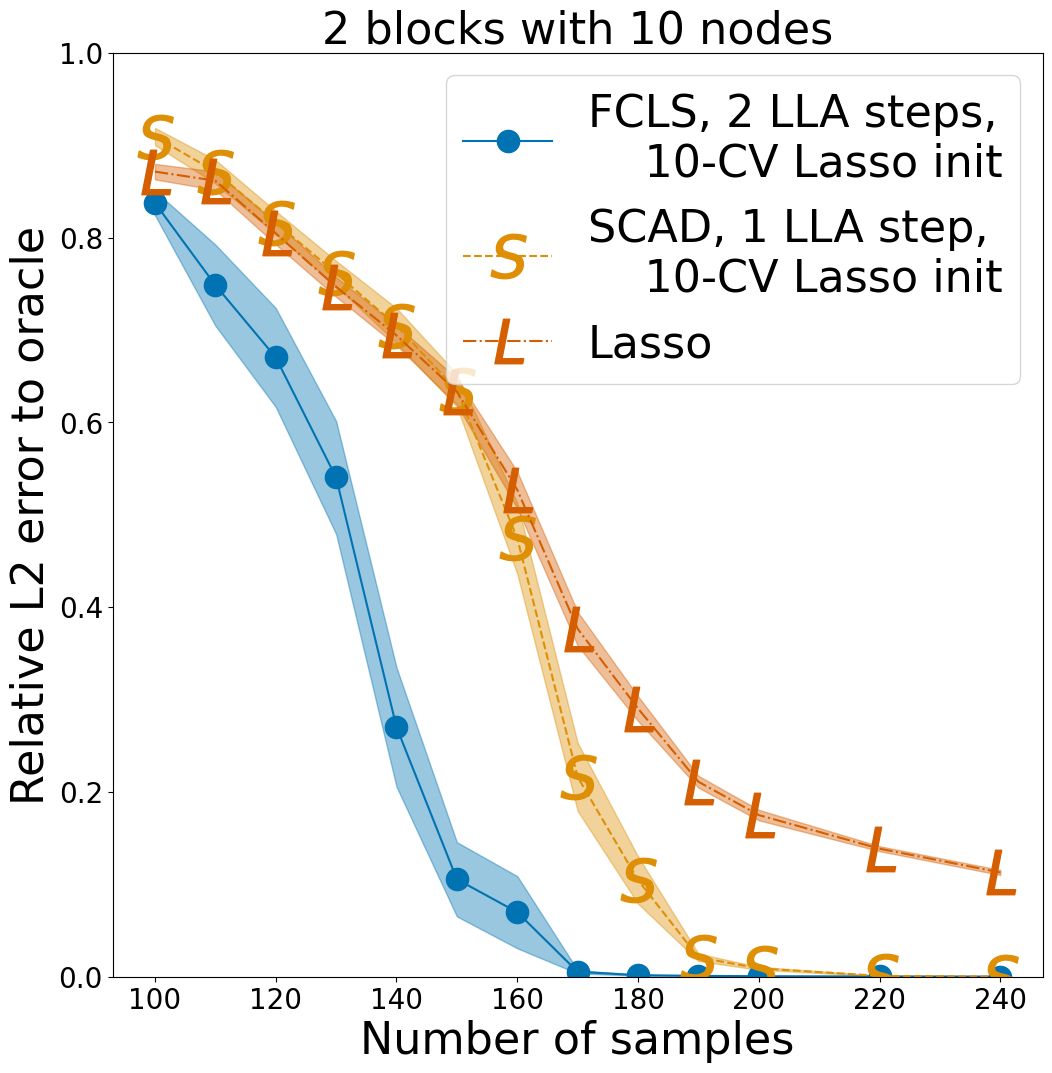}
\caption{
Two $10 \times 10$ blocks
}
\label{fig:lin_reg__bsize=10_2__vs__oracle__L2_rel}
\end{subfigure}
\hfill
\begin{subfigure}[t]{0.3\textwidth}
\centering
\includegraphics[width=\linewidth]{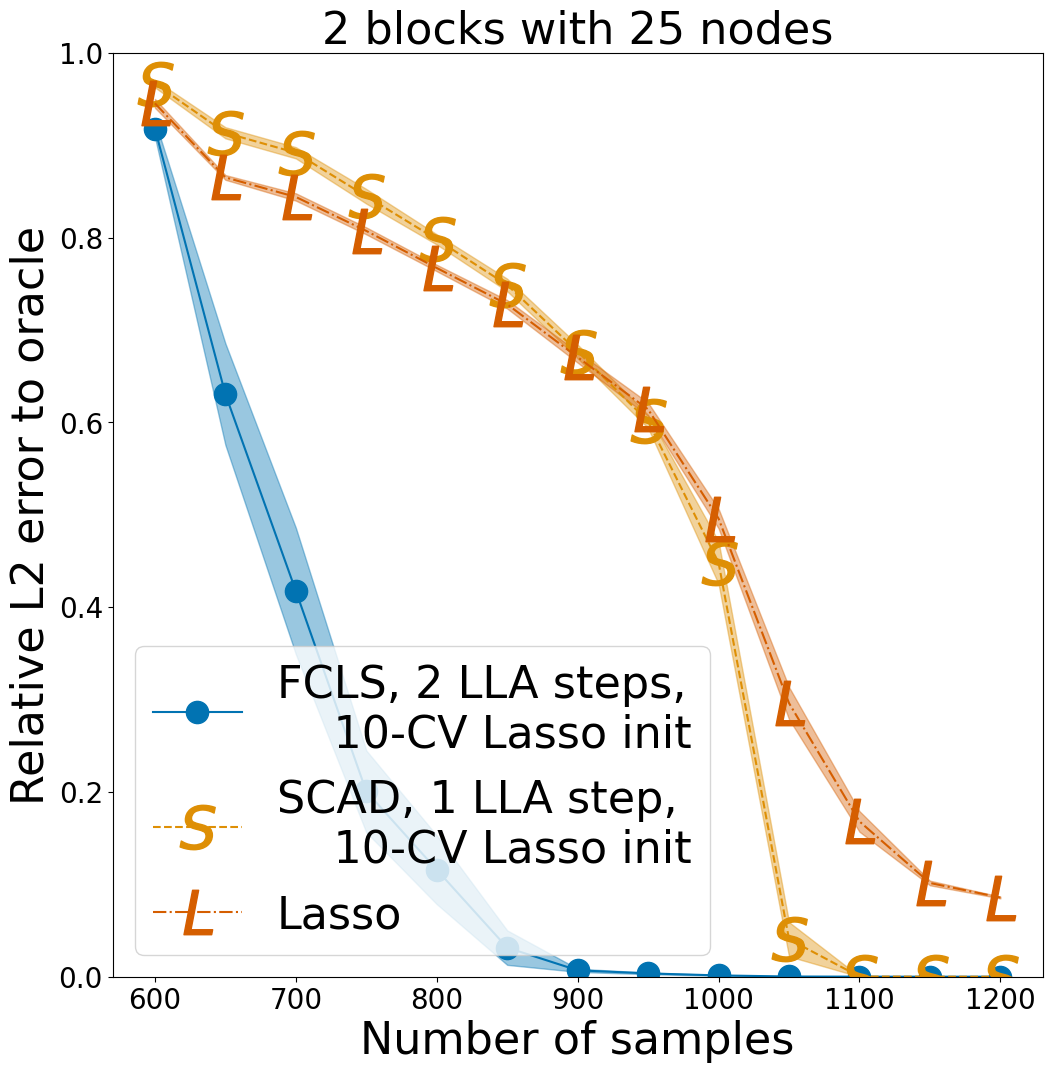}
\caption{
Two $25 \times 25$ blocks
}
\label{fig:lin_reg__bsize=25_2__vs__oracle__L2_rel}
\end{subfigure}
\caption{
Results for linear regression.
}
\label{fig:lin_reg__oracle__L2_rel}
\end{figure}

Figure \ref{fig:lin_reg__oracle__L2_rel} examines these models across a range of number of samples for each of the three graph topologies.
Both the FCLS solution and SCAD solution eventually return the block oracle estimate (with high probability) after enough samples.
The FCLS estimator, however, typically finds the oracle with fewer samples than its entrywise competitor. 
As in the previous sections \penabr\ penalized estimate has the best performance and the gap to the competing method appears to increase with larger block sizes.

\subsection{Logistic regression} \label{ss:sim__log_reg} 

\begin{figure}[H]
 \centering
\begin{subfigure}[t]{0.3\textwidth}
\centering
\includegraphics[width=\linewidth]{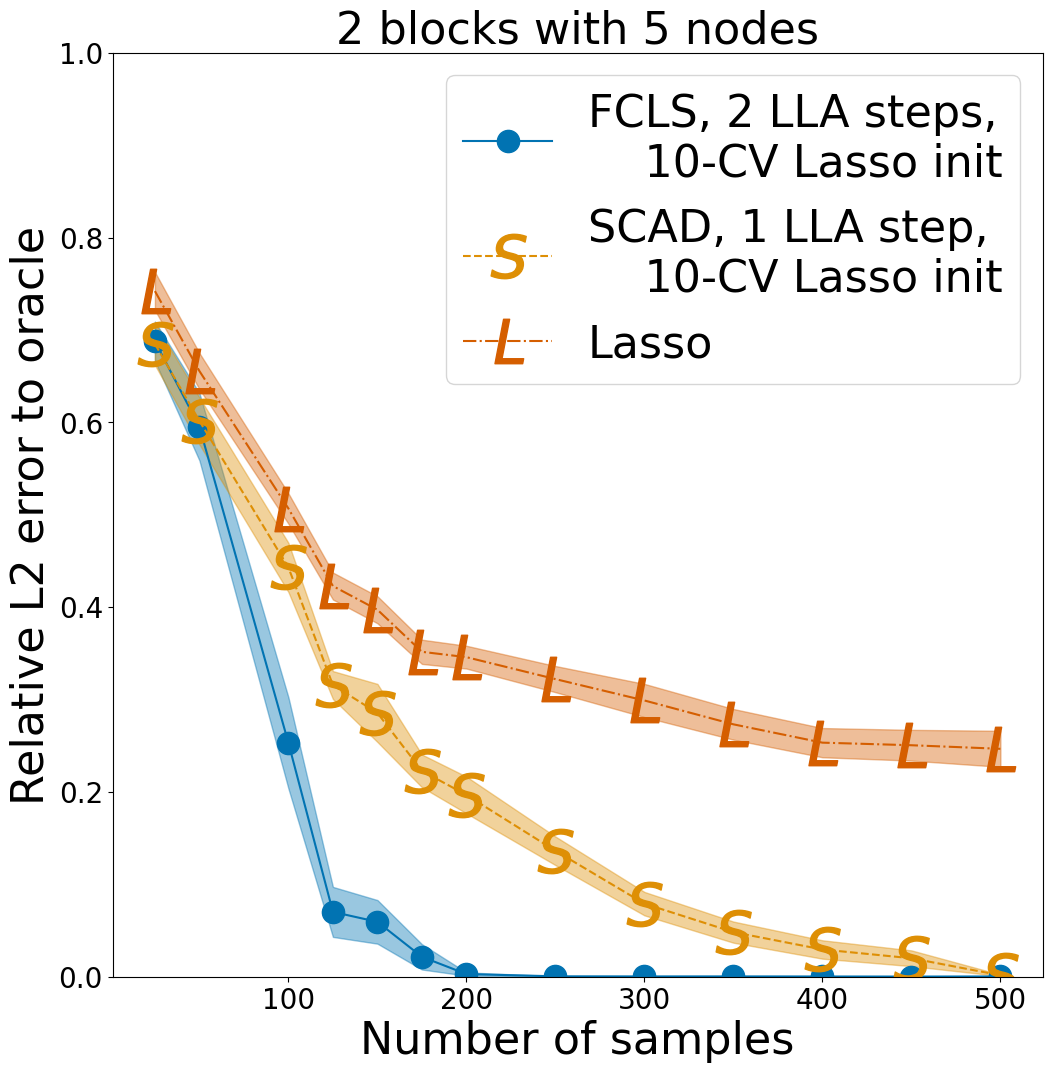}
\caption{
Two $5 \times 5$ blocks
}
\label{fig:log_reg__bsize=5_2__vs__oracle__L2_rel}
\end{subfigure}
\hfill
\begin{subfigure}[t]{0.3\textwidth}
\centering
\includegraphics[width=\linewidth]{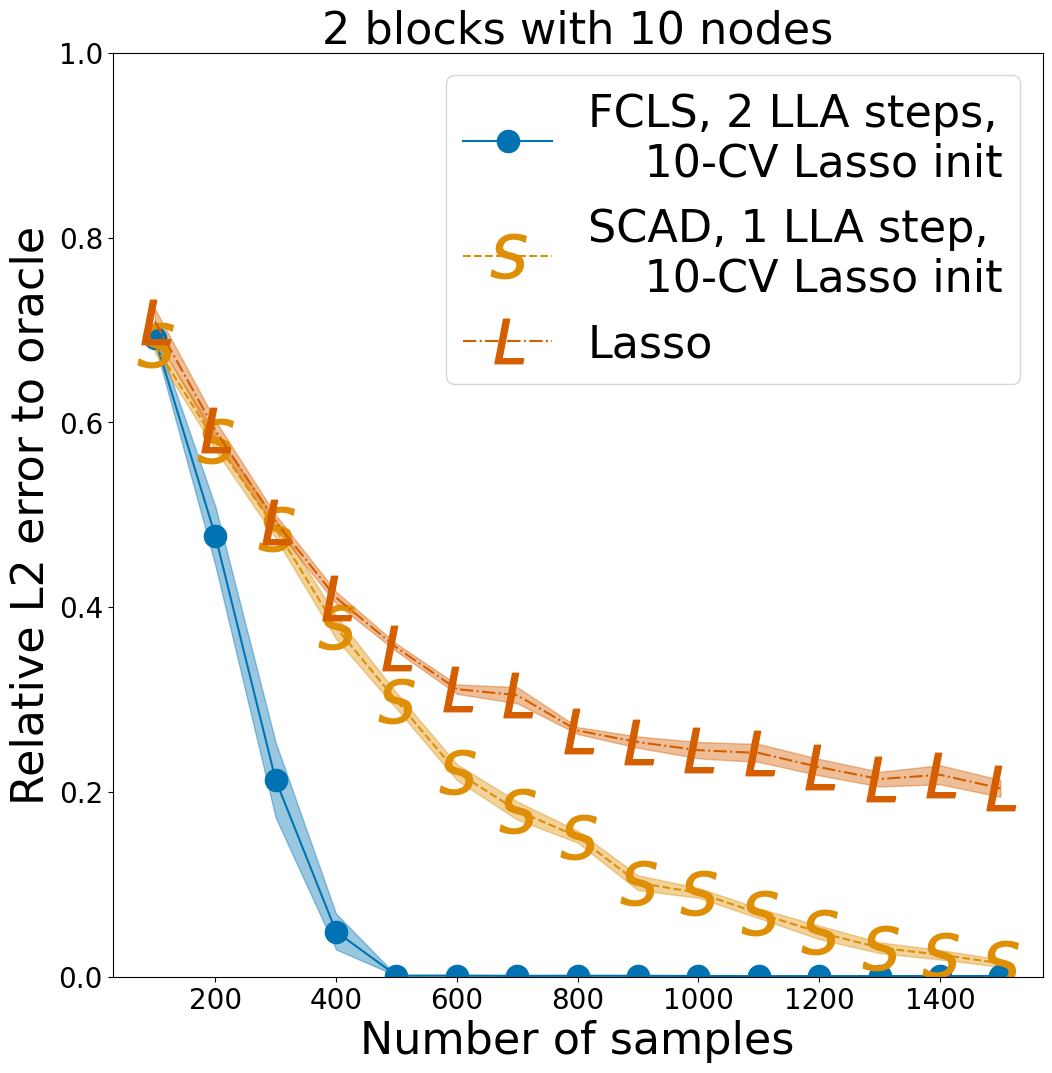}
\caption{
Two $10 \times 10$ blocks
}
\label{fig:log_reg__bsize=10_2__vs__oracle__L2_rel}
\end{subfigure}
\hfill
\begin{subfigure}[t]{0.3\textwidth}
\centering
\includegraphics[width=\linewidth]{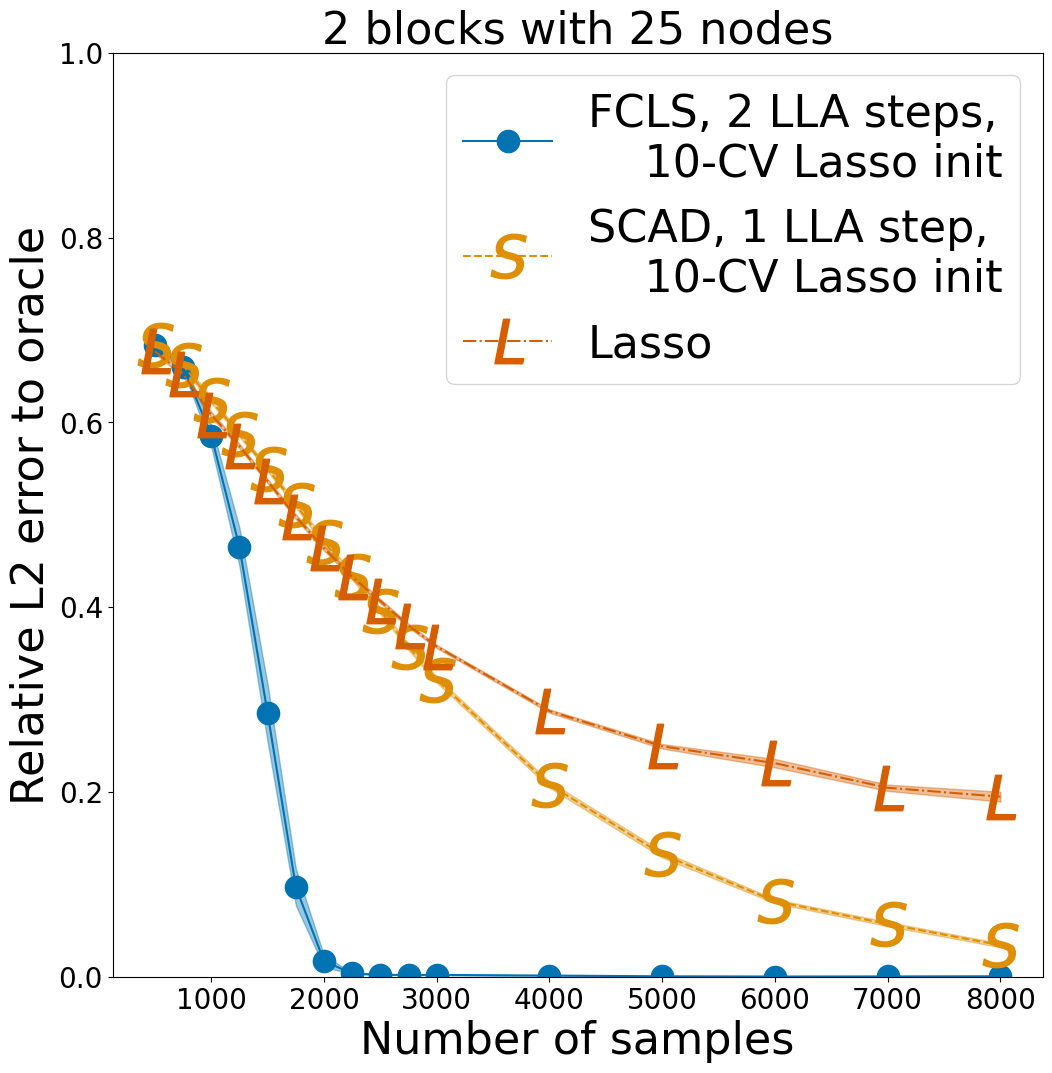}
\caption{
Two $25 \times 25$ blocks
}
\label{fig:log_reg__bsize=25_2__vs__oracle__L2_rel}
\end{subfigure}
\caption{
Results for logistic regression.
}
\label{fig:log_reg__oracle__L2_rel}
\end{figure}

This section examines the block sparse logistic regression example from Section \ref{ss:theo_ex_logisitc}.
Here we follow the setup from the previous section except of course the conditional distribution $y_i|x_i$ follows \eqref{eq:logistic_setup}.
We added a small ridge penalty of $0.01$ to the logistic regression loss function to improve stability.\footnote{Unpenalized logistic regression becomes unstable when the classes are linearly separable. This issue comes up when the FCLS penalized model solution finds the block oracle estimator and the block oracle estimator can perfectly separate the two classes.} 
In other words the loss function is given by
$$
\ell(\est) = \frac{1}{n} \sum_{i=1}  \left(-y_i X(i, :)^T \est + \LogisitcLink{X(i, :)^T \est}  \right)+ 0.01 \cdot \frac{1}{2} ||\est||_2^2
$$
instead of just the sum.
Note the block oracle estimator inherits the  ridge penalty.
Figure \ref{fig:log_reg__oracle__L2_rel} examines the performance of the competing methods and the takeaways are the same as in the previous sections.

\section{Discussion} \label{s:discussion}

This paper makes four contributions:
we propose the \penabr\ penalty, 
provide a majorization-minimization algorithm based on a weighted Lasso surrogate function,
develop a theoretical framework showing this penalty obeys the strong oracle property after 2 steps of the LLA algorithm,
and show this theory applies to several standard statistical models.
The two step convergence guarantees in Sections \ref{s:bd_shrink} and \ref{s:regression_examples} apply in ultra-high dimensional settings with an appropriate initializer (e.g. the Lasso estimate). 
The theory should apply to many statistical models beyond those studied in these two sections.

The simulations in Section \ref{s:sim} show that the \penabr\ penalized models outperform competing methods such as hard-thresholding or the entrywise SCAD penalty.
The rates provided in  Sections \ref{s:bd_shrink} and \ref{s:regression_examples}, however, are often worse than known rates for the entrywise competitors \citep{fan2014strong}.
These two facts suggest there is room for improvement in the theory of \penabr\  penalized models.

Finally we point out a few additional lines of future inquiry.
 Appendix \ref{a:bd_rect_multi_array} sketches extensions of the \penabr\ penalty to rectangular and multi-array settings, which may be interesting for a variety of applications such as multi-response regression. 
While we focus on the penalized estimator, Problem \eqref{prob:con_comp_pen}, the results in Section \ref{s:theory} can be adapted to a constrained formulation of this problem similar to ones considered in \citep{nie2016constrained, nie2017learning, kumar2019unified}.
In reality we may not expect $\esttarg$  to be perfectly block diagonal, but rather approximately block diagonal e.g. when $||\lambda\left( \mathcal{L}(|\esttarg|)\right)||_q$ is small for some $q \in [0, 1]$.
Extending our theory to this approximately block diagonal setting may be useful.
Decomposing $\est = \est_{\text{block diag}} + \est_{\text{sparse}}$ into the sum of a block diagonal plus sparse vector -- in the spirit of \citep{candes2011robust, tan2014learning} -- may also be of interest for applications.

\section*{Acknowledgments}
The author was supported by the National Science Foundation under Award No. 1902440.
We thank Mathurin Massias and Quentin Bertrand for adding a weighted Lasso solver for linear regression to their \textit{andersoncd} package.\footnote{\url{https://github.com/mathurinm/andersoncd/}} 
We also thank Deborah Carmichael, Calum Carmichael, and Ema Perkovic for helpful editorial feedback.


\appendix

\section{Block-diagonal rectangular matrices and multi-arrays} \label{a:bd_rect_multi_array}

In the body of the paper we assumed $\est \in \mathbb{R}^{{d \choose 2}}$ could be viewed as the entries of a hollow-symmetric matrix, $\mathcal{A}(\est) \in \mathbb{R}^{d \times d}$.
Then we used the  \penabr\  penalty to impose a block diagonal structure on $\mathcal{A}(\est)$.
This section extends this idea to block diagonal rectangular matrices and multi-arrays (see Figure \ref{fig:block_diag_examples}).
We sketch the main ideas, but do not pursue a detailed investigation here. 

\subsection{Rectangular matrices and the bipartite \penabr\ penalty}

Suppose $\est \in \mathbb{R}^{RC}$ can be viewed as the entries of a rectangular matrix $\textsc{mat}(\est) \in \mathbb{R}^{R \times C}$.
Let 
\begin{equation}
\mathcal{A}_{\text{bp}}(\est) = \begin{bmatrix} 0 & \textsc{mat}(\est) \\ \textsc{mat}(\est)^T & 0 \end{bmatrix} \in \mathbb{R}^{(R + C) \times (R + C)}
\end{equation}
be the adjacency matrix of the bipartite graph whose real valued edges are given by $\textsc{mat}(x)$.

Forcing $\mathcal{A}_{\text{bp}}(x)$ to have multiple connected components forces $  \textsc{mat}(\est)$ to have a block diagonal structure e.g. see Figure \ref{fig:bd_rect_ex}.
Let
\begin{equation} \label{eq:fcls_bp}
s_{\text{bp}, \tuneparam}(\est) :=  \frac{1}{2} g_{\tuneparam} \circ \lambda (\mathcal{L}(\mathcal{A}_{\text{bp}}(|\est|))),
\end{equation}
be the \textit{bipartite} \penabr\ penalty.
We can then use the LLA algorithm presented in Section \ref{ss:lla} to solve the following problem
\begin{equation} \label{prob:con_comp_pen_bp}
\underset{\est \in \mathbb{R}^{RC}}{\textup{minimize}}  \;\; \ell(\est) + s_{\text{bp}, \tuneparam}(\est).
\end{equation}


\subsection{Multi-arrays and the hypergraph \penabr\ penalty}

In the previous section we saw how to view a rectangular matrix as a bipartite graph.
There we viewed the rows and columns as vertices and the edges as connections between two vertices.
In this section we start with a multi-array $A \in \mathbb{R}^{d\spa{1} \times \dots \times d\spa{V}}$ and show how to view it as a graph whose connected components define a multi-array block diagonal structure.
The key step is to view the entries of $A$ as the hyperedges of a certain hypergraph.
Then we can follow the outline of the previous section and develop a multi-array version of the \penabr\ penalty.
The discussion in this section closely follows \citep{zhou2006learning, carmichael2020learning}.

A \textit{hypergraph} is a generalization of a graph where the hyperedges can connect more than two vertices \citep{berge1984hypergraphs, zhou2006learning}.
Formally a binary hypergraph $\HyperGraph = (\VetSet, \EdgeSet)$ is a collection of vertices $\VetSet$ and hyperedges $\EdgeSet$ where each hyperedge is a subset of $V$ with at least two elements (we do not allow self-loops here!)
For example, a binary graph is a hypergraph where each subset in $\EdgeSet$ has exactly two elements.
A weighted hypergraph  $\HyperGraph = (\VetSet, \EdgeSet, \EdgeWeights)$ has an additional map $\EdgeWeights: \EdgeSet \to \mathbb{R}_+$ that assigns a positive number to each hyperedge (we only consider positive weights for simplicity).

Next we interpret a multi-array\footnote{We will assume the entries of $A$ are positive; we can always take absolute values to make this happen!} $A \in \mathbb{R}^{d\spa{1} \times \dots \times d\spa{V}}_+$ as a hypergraph, $\HyperGraph(A)$ where the vertices are
$$
\VetSet(A) := \{ (v, i) : v \in [V], i \in [d\spa{v}]\}
$$
i.e. the generalization of rows and columns to multi-arrays.
Note $ |\VetSet(A) | = \DimSumMA := \sum_{v=1}^V d\spa{v}$.
The hyperedges correspond to the non-zero entries of $A$.
In other words the hyperedge $\{(1, i\spa{1}), \dots, (V, i\spa{V} ) \} \subseteq \VetSet(A)$ is present if and only if $A_{i\spa{1}, \dots, i\spa{V}} \neq 0$.
Note each present hyperedge has exactly $V$ elements; each edge contains exactly one element from each axis of $A$.
The hyperedge weights are then given by the entries of $A$.

We can identify the blocks of the multi-array $A$ with the connected components of the hypergraph  $\HyperGraph(A)$.
For example, the multi-array in Figure \ref{fig:bd_multi_array_ex} has 4 connected components.
The first connected component is $\{(1, 1), (2, 1), (3, 1)\}$ and corresponds to the bottom left red block i.e. the entry $A_{1 1 1}$.
The second connected component is $\{ (v, j) : j = 2, 3,  v = 1, 2, 3\}$ and corresponds to the green block whose entries are $\{ A_{ijk} : i, j, k = 2, 3\}$.
The isolated vertex $(1, 6)$ is also a connected component corresponding to an empty block.

Next we map the weighted hypergraph of $A$ into a weighted $V$-partite graph that we call the \textit{hypergraph adjacency matrix}.\footnote{This adjacency matrix implicitly defines a weighted, undirected graph.}
Let $\HyperAdjMat(A) \in \mathbb{R}^{\DimSumMA \times \DimSumMA}$ be a hollow symmetric matrix whose rows/columns are indexed by the vertices, $\VetSet(A)$.
The entries of this matrix correspond to a (weighted) graph where edges can only occur between two vertices on different axes of $A$ e.g. $(a, k\spa{a})$ and $(b, k\spa{b})$ where $a \neq b$.
The weight of such an edge is equal to
\begin{equation} \label{eq:hypergraph_adj_mat_edge_sum}
\HyperAdjMat(A)_{(a, k\spa{a}), (b, k\spa{b})} = \sum_{v \in [V] / \{a, b\}} \sum_{j\spa{v}=1}^{d\spa{v}} \HyperAdjMat(A)_{j\spa{1}, \dots, k\spa{a}, \dots, k\spa{b}, \dots j\spa{V}}.
\end{equation}
In other words, we take a two-dimensional slice of $A$ by fixing axes $a, b$ and squash the other dimensions down into this slice via summation.
 See Figure 10.c of \cite{carmichael2020learning}.
Note $\HyperAdjMat(A)$ is equivalent to the hypergraph adjacency matrix defined at the end of Section 2 in \citep{zhou2006learning}.
This graph captures the connected component information of $\HyperGraph(A)$.

\begin{proposition}  \label{prop:hypergraph_to_adj_mat_conn_compts}
There is a one-to-one correspondence between the connected components of the hypergraph $\HyperGraph(A)$ and the connected components of the graph $\HyperAdjMat(A)$.
\end{proposition}
The advantage of working with the graph $\HyperAdjMat(A)$ is we can use the previously developed Laplacian approach to imposing a connected component structure.
In other words, the spectrum of $\mathcal{L}(\HyperAdjMat(A))$ captures the connected components of $\HyperAdjMat(A)$ and therefore also  $\HyperGraph(A)$.

Finally suppose we have a model parameter $\est \in \mathbb{R}^{\DimProdMA}$ where $\DimProdMA := \prod_{v=1}^V d\spa{v}$ that can be viewed as the entries of a multi-array $\textsc{arr}(x) \in \mathbb{R}^{d\spa{1} \times \dots \times d\spa{V}}$.
We then define the \textit{hypergraph folded Laplacian spectral penalty} as
\begin{equation} \label{eq:fcls_multi_array}
s_{\text{hyper}, \tuneparam}(\est) := 
 \frac{1}{2} g_{\tuneparam} \circ \lambda \circ \mathcal{L} \circ \HyperAdjMat \circ \textsc{arr} ( |\est|).
\end{equation}
We can use the LLA algorithm presented in Section \ref{ss:lla} with minor modifications to solve the following problem
\begin{equation} \label{prob:con_comp_pen_multi_array}
\underset{\est \in \mathbb{R}^{\DimProdMA}}{\textup{minimize}}  \;\; \ell(\est) + s_{\text{hyper}, \tuneparam}(\est)
\end{equation}
that will impose a block diagonal structure on $\textsc{arr}( |\est|)$.


\section{Specifying the largest value for the \penabr\ tuning parameter} \label{a:tune_param_ubd}

In practice we need to compute the tuning path of solutions to the \penabr\ penalized Problem  \eqref{prob:con_comp_pen} over a range of $\tuneparam$ values in $[0, \tuneparam_{\text{max}}]$.
For Lasso penalized problems, we can often find an intelligent choice for the maximal tuning parameter \citep{friedman2010regularization}.

\begin{definition}
For a given loss function $\ell(\cdot)$ we say $\LassoKillerBound \in \mathbb{R}_+$ is a \textit{killer Lasso lower bound} if\footnote{This value does not have to exist!} for all weight vectors $c \in \mathbb{R}^D_+$ satisfying $\min (c) \ge \LassoKillerBound$ the solution to the following weighted Lasso problem is 0
\begin{equation} \label{prob:lasso_problem_pen}
\underset{\est \in \mathbb{R}^{D}}{\textup{minimize}}  \;\; \ell(\est) +  \sum_{j=1}^D c_j |\est_j|. 
\end{equation}
\end{definition}
We are not required to find the possible smallest value of $\LassoKillerBound$.
For convex loss functions
\begin{equation}
\LassoKillerBound = || \nabla \ell(0) ||_{\text{max}},
\end{equation}
is a killer Lasso lower bound.\footnote{To see this note $0$ is a stationary point if there is a solution to the first order necessary condition $\nabla_{\est} \ell(0) +  c \odot \nabla  ||0||_1 = 0$. Recalling the form of the sub-gradient of the $L_1$ norm we obtain the stated bound. Of course we need something like strict convexity to guarantee that 0 is the \textit{unique} stationary point.}
For the loss functions considered in Sections \ref{s:bd_shrink} and \ref{s:regression_examples} we have the following explicit expressions:
\begin{itemize}

\item For the Frobenius norm loss function in \eqref{prob:prox_fclsp}, $\LassoKillerBound = ||\estok||_{\text{max}} $.

\item For the least squares loss function without an intercept in \eqref{prob:lin_reg_fclsp}, $\LassoKillerBound = \frac{1}{n} ||X^T y||_{\text{max}}$. 

\item For the logistic loss function without an intercept in \eqref{prob:logistic_reg_fclsp}, $\LassoKillerBound = \frac{1}{n} ||X^T (\frac{1}{2} \mathbf{1}_n - y)||_{\text{max}}$.
\end{itemize}

We can use a killer Lasso lower bound to obtain the largest reasonable tuning parameter value for the LLA algorithm from Section \ref{ss:lla}.
\begin{proposition} \label{prop:fcls_tune_param_max_val}
Let $\estinit \in \mathbb{R}^D$ be the initializer to the LLA algorithm for Problem \eqref{prob:con_comp_pen} where $g_{\tuneparam}(\cdot)$ satisfies Definition \ref{def:scad_like_pen_func}.
Suppose $\LassoKillerBound$ is a killer lasso lower bound and let
\begin{equation} \label{eq:fcls_tune_param_max_val}
\FCLSKillerBound := \max\left(
\frac{ \lambda_{\text{max}} \left(  \mathcal{L}(|\estinit|) \right) }{b_1},
\frac{\ThreshTuneParam_{\text{lasso-max}} }{a_1}
\right).
\end{equation}
Then for all $\tuneparam \ge \FCLSKillerBound$ the LLA algorithm initialized from $\estinit$ converges to 0 in one step.
\end{proposition}


\section{Laplacian spectral bounds} \label{a:lap_spect_bounds}

This section studies quantities and bound related to the Laplacian that play an important role in the theory.
Section \ref{as:graph_spectra} reviews the spectra of some important example graphs.
Section \ref{as:lap_norm_vs_adj_norm} compares norms of $\mathcal{L}(r)$ to norms of the edge vector $r$.
In particular, Remark \ref{rem:lap_comparision_strategy} summarizes facts that we use frequently in the proofs. 
Finally, Section  \ref{ss:lap_coeff_comparision} studies the Laplacian coefficient from Section \ref{ss:lla}.

\subsection{Spectra of some important binary graphs} \label{as:graph_spectra}

The largest and second smallest Laplacian eigenvalues plan an important role in our theory.
Table \ref{tab:lap_examples} shows the largest and second smallest Laplacian eigenvalues for some important graphs binary, undirected graphs.

\begin{table}[H]
\center
\begin{tabular}{|l|l|l|}
\hline
                                      & $\lapSmEval{2}{A}$                   &    $\lambda_{\text{max}} \left(\mathcal{L}(A)\right)$   \\ \hline
Fully connected graph  &    d                                              &    $d$       						                      \\ \hline
Star graph                     &    1                                             &    $1$           						                      \\ \hline
Path graph                    &   $O\left(\frac{1}{d^2} \right) $    &    $O(1)$         							     \\ \hline
Any graph                    &   $ \ge \frac{d}{1 + \graphDiam \cdot \numMissingEdges }$                    &         							                       \\ \hline
\end{tabular}
\caption{
Largest and second smallest Laplacian eigenvalues for connected graphs on $d$ nodes.
Here $A \in \{0, 1\}^{d \times d}$ is the adjacency matrix of graph listed in the first column.
Here we see, unsurprisingly $\lapSmEval{2}{A}$, is generally larger for graphs with more edges though the path graph is a counter example of this trend.
}
\label{tab:lap_examples}
\end{table}
Here $\graphDiam \le d-1$ is the diameter of the graph and let $\numMissingEdges \le {d \choose 2}$ be the number of edges missing from $A$.
The claims about the three simple graphs are standard results in spectral graph theory e.g. see \citep{spielman2004spectral}.
The lower bound on the algebraic connectivity for any graph is from Theorem 1 of \citep{lu2007lower}; a better albeit less transparent lower bound can be found in \citep{rad2011lower}.
Using this lower bound we can generalize \eqref{eq:targ_gap_from_fully_conn} to
\begin{equation} \label{eq:targ_gap_with_some_missing_lu07}
\targgap \ge  \frac{\minccTarg}{1 + \maxMissingTarg \maxDiamTarg}  ||\esttarg_{\nnzSetTarg}||_{\text{min}}.
\end{equation}
where $\nnzSetTarg = \{ (ij) | \esttarg_{(ij)} \neq 0\}$ is the support set of $\esttarg$.

\subsection{Laplacian vs. adjacency matrix norms} \label{as:lap_norm_vs_adj_norm}

The goal of this section is to bound norms of the Laplacian $||\mathcal{L}(\cdot)||$ by norms of $||r||$.
Throughout this section we assume $r \in \mathbb{R}^D$ where $D = {d \choose 2}$.
The bounds\footnote{Similar expressions could be derived by starting with \eqref{eq:lap_bound_1__1} and using the fact $||r||_1 \le \sqrt{D} ||r||_{2} \le D ||r||_{\text{max}}$ however these would be looser than the stated bounds.} in Proposition \ref{prop:lap_bounds} are also shown to be tight (at least up to constants) without further assumptions.

\begin{proposition} \label{prop:lap_bounds}
Let $r \in \mathbb{R}^{D}$, then
\begin{align}
||\mathcal{L}(r)||_{1} & \le 4 ||r||_1  \label{eq:lap_bound_1__1} \\
||\mathcal{L}(r)||_F & \le \sqrt{2d} ||r||_2  \label{eq:lap_bound_F__2} \\
||\mathcal{L}(r)||_{\text{F}}  &  \le \sqrt{d^{3} + d^2 - d} ||r||_{\text{max}} \label{eq:lap_bound_F__max} \\
||\mathcal{L}(r)||_{F} & \le \sqrt{2d} ||\mathcal{A}(r)||_{\text{op}} \label{eq:lap_bound_frob__op} \\
\opn{\mathcal{L}(r)} & \le 2 \opone{\mathcal{A}(r)} \le 2 d ||r||_{\text{max}} \label{eq:lap_bound_op__op1}  \\
\opn{\mathcal{L}(r)} & \le ||\mathcal{A}(r) \mathbf{1}_d||_{\text{max}} + \opn{\mathcal{A}(r)}. \label{eq:lap_bound_op__op_plus_deg} 
\end{align}
If the largest degree of a node in the binary graph of $\matsupp{\mathcal{A}(r)}$ is equal to $\maxDegTarg$ then $d$ can be replaced with $\maxDegTarg$ for claims \eqref{eq:lap_bound_1__1}, \eqref{eq:lap_bound_F__2} and \eqref{eq:lap_bound_op__op1}.

Furthermore we have the lower bound $||\mathcal{A}(r) \mathbf{1}_d ||_{\text{max}}   \le \opn{\mathcal{L}(r)}$.
\end{proposition}

For a symmetric matrix $A \in \mathbb{R}^{d \times d}$ we have the following inequalities:  $\opn{A}  \le ||A||_F \le ||A||_1$, $||A||_F \le \sqrt{d} \opn{A}$, and $\opn{A} \le \opone{A} \le \sqrt{d} \opn{A} \le d ||A||_{\text{max}}$. 
Proof of the second two inequalities can be found in Section 2.3 of \citep{golub2013matrix}.
Using these expressions we can obtain additional upper bounds on $||\mathcal{L}(r)||_{\text{op}}$ from \eqref{eq:lap_bound_1__1},  \eqref{eq:lap_bound_F__2}, and \eqref{eq:lap_bound_op__op1} as well as on on $||\mathcal{L}(r)||_{F}$  from \eqref{eq:lap_bound_op__op1}.
All pairs of upper bounds are summarized in Table \ref{tab:lap_from_adj_upper_bounds}.

\begin{table}[H]
\centering
\begin{tabular}{|l|l|l|l|l|l|l|}
\hline
                                                      & $||r||_{1}$                 & $||r||_{2}$              & $||r||_{\text{max}}$      & $\opone{\mathcal{A}(r)}$       & $\opn{\mathcal{A}(r)}$     &  $||\mathcal{A}(r) \mathbf{1}_d||_{\text{max}} + \opn{\mathcal{A}(r)}$    \\ \hline
$||\mathcal{L}(r)||_{F}$                  &  $O(1)$        & $O(d^{1/2}) *$         &  $O(d^{3/2})$               & $O(d^{1/2}) *$                          & $O(d^{1/2})$                    & $O(d^{1/2})$                   \\ \hline
$\opn{\mathcal{L}(r)}$                   &  $O(1)$                    & $O(d^{1/2}) *$          &   $O(d) *$                       &  $O(1)$                      &  $O(d^{1/2})$                  &  $O(1)$                     \\ \hline
\end{tabular}
\caption{
Upper bound on the ratio $\frac{ || \mathcal{L}(r) || }{ ||\mathcal{A}(r)||}$ for different pairs of norms from Proposition \ref{prop:lap_bounds}.
A $*$ means the dimension $d$ can be replaced with the maximal degree $\maxDegTarg$.
}
\label{tab:lap_from_adj_upper_bounds}
\end{table}

All of the upper bounds in Proposition \ref{prop:lap_bounds} are as good as we can ask for (up to constants) without additional assumptions.
Table \ref{tab:example_graph_lap_and_adj_norms} shows the ratios from Table \ref{tab:lap_from_adj_upper_bounds} for three example graphs.
For example, the star graph shows the $\sqrt{d}$ factor is necessary for the ratio $|| \mathcal{L}(r) ||_{\text{op}} / ||r||_{2}$.
These upper bounds, however, may be loose for graphs observed in practice.
For example, upper bound on the Laplacian operator norm vs. $||r||_1$ or $||r||_2$ is loose for the complete graph.
\begin{table}[H]
\centering
\begin{tabular}{|l|l|l|l|}
\hline
                                                        								& Complete graph       & Star graph                      & Path graph \\ \hline
$||r||_{\text{max}}$                          								 & $O(1)$ (op, F)        &  $O(1)$ (op)                   & $O(1)$  \\ \hline
$||r||_{2}$                                     							        & $O(d)$  (F)              & $O(\sqrt{d})$ (op, F)        & $O(\sqrt{d})$ \\ \hline
$||r||_{1}$                                       							       &   $O(d^2)$                & $O(d)$ (op, F)                 & $O(d)$ \\ \hline
$\opn{\mathcal{A}(r)}$       								       & $O(d)$ (F)               & $O(\sqrt{d})$ (op, F)         & $O(1)$ (F)  \\ \hline
 $||\mathcal{A}(r) \mathbf{1}_d||_{\text{max}} + \opn{\mathcal{A}(r)}$       & $O(d)$ (op, F)           & $O(d)$  (op)                     & $O(1)$ (op, F)  \\ \hline
$\opone{A(r)}$                                							     &  $O(d)$ (op, F)         &  $O(d)$  (op)                     & $O(1)$  (op, F) \\ \hline \hline 
$\opn{\mathcal{L}(r)}$      									     & $O(d)$                     & $O(d)$                              &  $O(1)$ \\ \hline
$||\mathcal{L}(r) ||_{F}$                  							   &  $O(d^{3/2})$           & $O(d)$                               &  $O(d^{1/2})$  \\ \hline
\end{tabular} 
\caption{
Norms of three example graphs on $d$ nodes (up to multiplicative constants).
An $(op)$ means the upper bound on the ratio for $ \opn{\mathcal{L}(r)}$ given in Table \ref{tab:lap_from_adj_upper_bounds} is met.
An $(F)$ means the same thing for $||\mathcal{L}(r) ||_{F} $.
For each upper bound there is an example where this upper bound is met, however, for several of the examples the upper bound is loose.
}
\label{tab:example_graph_lap_and_adj_norms}
\end{table}

Due to the importance of controlling residuals of the Laplacian we explicitly state the bounds we use frequently.
\begin{remark} \label{rem:lap_comparision_strategy}
Let $r, x, y \in \mathbb{R}^D$ then
\begin{equation}
||\mathcal{L}(r)||_F \le 4 ||r||_1  \wedge \sqrt{2d} ||r||_2  
\end{equation}
\begin{equation}
\opn{\mathcal{L}(r)}
\le 2d ||r||_{\text{max}}
\wedge  \opone{\mathcal{A}(r)} 
\wedge \left( ||\mathcal{A}(r) \mathbf{1}_d||_{\text{max}} + \opn{\mathcal{A}(r)}\right).
\end{equation}
We also have the obvious bound that $\opn{\mathcal{L}(r)} \le ||\mathcal{L}(r)||_F$.
If $\mathcal{A}(r)$ has a maximal degree of $\maxDegTarg$  then $d$ can be replaced with $\maxDegTarg$ in these two equations.
If $\mathcal{A}(r)$ has $\nccTarg$ connected components indexed by $\ccTarg{1}, \dots, \ccTarg{\nccTarg}$ then $\opn{\mathcal{L}(r)} = \sup_{k \in [\nccTarg]} \opn{\mathcal{L}(r_{\ccTarg{k}})}$.

We need to be careful about comparing  $\mathcal{L}(|x| - |y|)$ to $\mathcal{L}(x - y)$.
In general we do \textbf{not} have the entrywise relationship\footnote{E.g. consider  $r = [-1, -1, 1]^T$ so $\sum_{i=1}^3 |r_i| = 3 > 2 = |\sum_{i=1} r_i |$.} $|\mathcal{A}(|x| - |y|) \mathbf{1}_d| \le |\mathcal{A}(x - y) \mathbf{1}_d|$.
Therefore we are not guaranteed to be able to immediately compare $||\mathcal{L}(|x| - |y|)||$ to $||\mathcal{L}(x - y)||$.
We do have $|\mathcal{A}(|x| - |y|) \mathbf{1}_d| \le \mathcal{A}(|x - y|) \mathbf{1}_d$ so we can immediately compare $\mathcal{L}(|x| - |y|)$ to $\mathcal{L}(|x - y|)$ using entrywise norms.
\end{remark}

\subsection{Laplacian coefficient comparison} \label{ss:lap_coeff_comparision}


The following proposition shows how much $\lapcoeff{\cdot}{w}$ can change as a function of its first argument.
\begin{proposition}\label{prop:lap_coef_bound_by_V}
Let $w \in \mathbb{R}^{K}_+$ and $V_A, V_B \in \mathbb{R}^{d \times K}$ then
\begin{equation} \label{prop:lap_coef_bound_single_edge}
|\lapcoeff{V_A}{w}_{(ij)}^{1/2}  - \lapcoeff{V_B}{w}_{(ij)}^{1/2} | \le 2 \sqrt{\max(w)} \opnorm{ V_A - V_B}{2}{\infty}
\end{equation}

Let $\mathcal{G}$ be a binary, undirected graph on $d$ nodes and let $\maxccTarg$ be size of the largest connected component, then
\begin{equation}\label{eq:lap_coef_bound_graph_sum_L1}
\sum_{(ij) \in \mathcal{G}} |\lapcoeff{V_A}{w}_{(ij)}^{1/2}  - \lapcoeff{V_B}{w}_{(ij)}^{1/2} |^2 \le 2 \max(w) \maxccTarg  ||V_A - V_B||_{F}^2 
\end{equation}
and
\begin{equation}\label{eq:lap_coef_bound_graph_sum_L2}
\sum_{(ij) \in \mathcal{G}} |\lapcoeff{V_A}{w}_{(ij)}^{1/2}  - \lapcoeff{V_B}{w}_{(ij)}^{1/2} |^4 \le 8 \max(w)^2 \maxccTarg  ||V_A - V_B||_{F}^4
\end{equation}
where each edge is only counted once in the sum.

\end{proposition}

The following proposition shows that the surrogate function $Q(\cdot | \estCurrent)$ only penalizes between block edges when $\estCurrent$ has multiple connected components.
\begin{proposition} \label{prop:lap_coef_with_con_comps}
Let $A \in \mathbb{R}^{d \times d}_+$ be the adjacency matrix of a graph with $K$ connected components.
Suppose $V \in \mathbb{R}^{d \times K}$ is any orthonormal matrix of eigenvectors corresponding to the $K$ 0 eigenvalues of $\mathcal{L}(A)$ (i.e. a basis for the kernel).
For $i, j \in [d]$ we write $i \sim_{CC} j$ if nodes $i$ and $j$ are in the same connected component of $A$.
Then for any $w \in \mathbb{R}^{K}_+$,
\begin{equation*}
\begin{aligned}
\lapcoeff{V}{w}_{ij} &  = 0  & & \text{ if } i \sim_{CC} j  \\
\lapcoeff{V}{w}_{ij} & \ge \min(w) \left( \frac{1}{|C(i)|} + \frac{1}{|C(j)|} \right)  \ge \frac{2 \min(w)}{\maxccTarg} & & \text{ if } i \not\sim_{CC, y} j. 
\end{aligned}
\end{equation*}
where $|C(i)|$ denotes the number of vertices in the connected component that the vertex $i$ belongs to and $\maxccTarg$ is the size of the largest connected component.
\end{proposition}

The following lemma shows that when $\estCurrent$ is close to a graph with multiple connected components, the surrogate function $Q(\cdot | \estCurrent)$ puts large weights on between block edges and small weights on within block edges.

\begin{lemma} \label{lem:lasso_maj_weight_bounds}
Let $x, y \in \mathbb{R}^{D}_+$.
Suppose $\mathcal{A}(y)$ is the adjacency matrix of a graph with $K$ connected components and $\lapSmEval{K + 1}{y} \ge \triangle$.
Let $V_x \in \mathbb{R}^{n \times K}$ be any orthonormal matrix of the smallest $K$ eigenvectors of $\mathcal{L}(x)$ 
and let $w \in \mathbb{R}^{K}_+$.
Let $\suppSet = \blocksupp{\mathcal{A}(y)}$ be the block support of $\mathcal{A}(y)$,
$\maxccTarg$ be the largest number of nodes in a connected component of $\mathcal{A}(y)$,
$$
\frobBound = \min \left( ||\mathcal{L}(x - y)||_F, \nccTargPow{1/2} \opn{\mathcal{L}(x - y)} \right), 
\text{ and }
\opBound =  \opn{\mathcal{L}(x - y)}.
$$

Then
\begin{equation} \label{eq:lasso_maj_weight_bounds_from_frob}
\begin{aligned}
||\lapcoeff{V_x}{w}_{\suppSet^C}||_{\text{min}}  & \ge \left( \sqrt{\frac{2 \min(w)}{\maxccTarg}} -  2^{5/2} \sqrt{\max(w)} \frac{\frobBound}{\triangle}  \right)^2 \\
||\lapcoeff{V_x}{w}_{\suppSet}||_{\text{max}} &  \le  2^{5} \max(w)  \left( \frac{\frobBound}{\triangle}  \right)^2 \\
||\lapcoeff{V_x}{w}_{\suppSet}||_1 & \le 2^4 \maxccTarg \max(w) \left( \frac{\frobBound}{\triangle} \right)^2 \\
||\lapcoeff{V_x}{w}_{\suppSet}||_2 & \le  2^{9/2} \maxccTargPow{1/2} \max(w) \left( \frac{\frobBound}{\triangle} \right)^2
\end{aligned}
\end{equation}

If $\opBound \le \frac{1}{2} \triangle$ then
\begin{equation} \label{eq:lasso_maj_weight_bounds_from_opn}
\begin{aligned}
||\lapcoeff{V_x}{w}_{\suppSet^C}||_{\text{min}}  & \ge \left( \sqrt{\frac{2 \min(w)}{\maxccTarg}}  -  2^{5/2} \sqrt{\max(w)} \frac{\opBound}{\triangle}  \right)^2\\
||\lapcoeff{V_x}{w}_{\suppSet}||_{\text{max}} &  \le  2^{5} \max(w)  \left( \frac{\opBound}{\triangle}  \right)^2 
\end{aligned}
\end{equation}

\end{lemma}

The operator norm part of the previous lemma depends on the following variant of the Davis-Kahan theorem.
\begin{lemma} \label{lem:dk_kernel}
Let $A, B \in \mathbb{R}^{d \times d}$ be symmetric and positive semi-definite.
Suppose $A$ has a $K$ dimensional kernel and $U_{A} \in \mathbb{R}^{d \times K}$ is any orthonormal basis for this kernel.
Let $\lambda_{B, \perp} \in \mathbb{R}^{d - K}_+$ be the leading $d-K$ eigenvalues of $B$ and let $U_{B, \perp} \in \mathbb{R}^{d \times (d - K)}$ be a corresponding orthonormal matrix of eigenvectors.
Finally let $\triangle_A = \lambda_{(K+1)}(A)$ and $\triangle_B = \lambda_{(K+1)}(B)$.

Then there exists an orthonormal matrix $Q \in \mathbb{R}^{K \times K}$ such that
\begin{equation}
\opn{U_B - U_A Q}
\le \sqrt{2} \opn{U_{B, \perp}^T U_A}
%
\le \frac{ \sqrt{2} \opn{B U_A}}{\triangle_B}
\le \frac{  \sqrt{2} \opn{B - A}}{\triangle_B}.
\end{equation}

Furthermore, if $\opn{B - A} \le \frac{1}{2} \triangle_A$ then
\begin{equation}
\opn{U_B - U_A Q}
\le 2 \sqrt{2} \opn{U_{B, \perp}^T U_A}
%
\le \frac{ 2\sqrt{2} \opn{B U_A}}{\triangle_A}
\le \frac{ 2 \sqrt{2} \opn{B - A}}{\triangle_A}.
\end{equation}
\end{lemma}
From the proof of this Lemma we can see that a similar result holds for any unitarily invariant matrix norm. 

\section{Proofs for Sections \ref{s:fcls_pen_opt} and \ref{s:theory}} \label{a:proofs_lla}

\begin{fact} \label{fact:evals_increase_for_pos_mats}
Let $A, B \in \mathbb{R}^{d \times d}_+$ be two hollow symmetric non-negative matrices then
$$
A \ge B \implies \lambda_j(\mathcal{L}(A)) \ge \lambda_j(\mathcal{L}(B)), \text{ for each } j=1, \dots, d,
$$
where the inequality applies entrywise.
\end{fact}
This follows from the fact $\mathcal{L}(A - B)$ is semi-definite (since $A - B \ge 0$ entrywise) and that summing two positive semi-definite matrices increases the eigenvalues.

\subsection{Proofs for Section \ref{ss:lla}} \label{as:proof_lla}

Recall $\lapcoeff{V}{w}$ defined in Section \ref{ss:lla}.
A useful property of this quantity is for any $\est \in \mathbb{R}^D$,
\begin{equation} \label{eq:lap_coef_implicit_def}
\est^T \lapcoeff{V}{w}  = \frac{1}{2} \sum_{i, j = 1, i \neq j}^d \est_{(ij)}  \sum_{k=1}^K w_k ||V(i, k) - V(j, k)||_2^2  = \sum_{k=1}^K w_k V_k^T \mathcal{L}(\est)  V_k =  \text{Tr}(V^T \mathcal{L}(\est) V \text{diag}(w)).
\end{equation}

\begin{proof} of Proposition \ref{prop:spect_pen_maj}
Following the discussion in Section \ref{ss:lla} and the fact that $\lapSpecFunc{\tuneparam}$ is concave we need only check the super-gradient formula \eqref{eq:lap_spect_grad} to verify  $Q$ is a surrogate function. 

Let $G: \mathbb{R}^{d \times d}$ be given by $G(L) = g_{\tuneparam} \circ \lambda(L)$ i.e. $\lapSpec{\tuneparam}{x} = f(\mathcal{L}(x))$.
By Theorem 6 of \citep{lewis1999nonsmooth}, we have the super-differential formula
$$
\nabla G(L) = V \text{diag}(w) V^T
$$
where $V \in \mathbb{R}^{d \times d}$ is any orthonormal matrix whose columns are eigenvectors of $L$ and $w = g_{\tuneparam}' \circ \lambda(L)$ is any super-gradient.

Recall $\mathcal{L}: \mathbb{R}^{d} \to \mathbb{R}^{d \times d}$ is a linear function and $\frac{d \mathcal{L}(x)}{d  x_j} = \mathcal{L}(e_j)$ where $j$ is the $j$th standard basis vector.
Let $V$ be any matrix of eigenvectors of $L(x)$ and $w = g_{\tuneparam}' \circ \lambda(L(x))$ any super-gradient.
Therefore by the chain rule e.g. see \citep{petersen2012matrix}
\begin{equation*}
\frac{d}{d x_j} \lapSpec{\tuneparam}{x}   
= \left\langle \frac{d \mathcal{L}(x)}{d  x_j}, \nabla G( \mathcal{L}(x))  \right \rangle 
  = \langle L(e_j), V \text{diag}(w) V^T \rangle 
 = \text{Tr} \left( V^T \mathcal{L}(e_j) V \text{diag}(w)  \right) 
 = e_j^T \lapcoeff{V}{w},
\end{equation*}
where the final equality is readily verified using \eqref{eq:lap_coef_implicit_def}.
Thus $\frac{d}{d x} \lapSpec{\tuneparam}{x}  =  \lapcoeff{V}{w}$.
We therefore conclude $Q$ is a surrogate function. 
Note that $Q(\cdot | \estCurrent)$ is convex since each entry of $\lapcoeff{V}{w}$ is non-negative by \eqref{eq:lap_coeff_formula}.

Next we show the choice of eigenvectors does not matter.
Let $V,\widetilde{V} \in \mathbb{R}^D$ be two orthonormal matrices whose columns are eigenvectors of $\mathcal{L}(|y|)$ and let $w = g'_{\tuneparam}(\lambda^{y})$.
 Recall from \eqref{eq:lap_coeff_formula} that $\lapcoeff{V}{w}_{(ij)} =   ||V(i, :) - V(j, :)||_{2, w}^2.$
First note that flipping the sign of a column of $V$ leaves $\lapcoeff{V}{w}$ unchanged.
Thus if all of the eigenvalues are unique the claim follows.

WLOG assume that the first $K$ eigenvalues are equal to each other and all other eigenvalues are unique (otherwise the following argument easily generalizes).
Then $V_{1:K} = \widetilde{V}_{1:K} T$ for some orthonormal matrix $T \in \mathbb{R}^{K \times K}$.
Then we can check
\begin{equation}
\begin{aligned}
||V(i, :) - V(j, :)||_{2, w}^2 
& = ||V_{1:K} (i, :) - V_{1:K} (j, :) ||_{2, w_{1:K}}^2  + ||V_{(K+1):D}(i, :) - V_{(K+1):D}(j, :)||_{2, w_{(K +1):D}}^2  \\
& = w_1 || V_{1:K} (i, :) - V_{1:K} (j, :)  ||_{2}^2  + \sameAsPrev
  \\
& = w_1  (V_{1:K} (i, :) - V_{1:K} (j, :) )^T (V_{1:K} (i, :) - V_{1:K} (j, :) )  + \sameAsPrev\\
& = w_1   (V_{1:K} (i, :) - V_{1:K} (j, :) )^T T^T T (V_{1:K} (i, :) - V_{1:K} (j, :) )   + \sameAsPrev \\
& = w_1 ||V_{1:K}T - V_{1:K} T ||_{2}^2   + \sameAsPrev  \\
& = ||\widetilde{V}(i, :) - \widetilde{V}(j, :)||_{2, w}^2,
\end{aligned}
\end{equation}
Thus $\lapcoeff{V}{w} = \lapcoeff{\widetilde{V}}{w}$.

Finally we check that majorizing at $0$ leaves a Lasso penalty.
Note $\mathcal{L}(0) = 0_{D \times D}$ so $\lambda_y= 0$ and $V_y$ is any  orthonormal basis matrix of $\mathbb{R}^D$ thus
$||V_i - V_j||_2^2 = ||V_i||_2^2 + ||V_j||_2^2  + 2 V_i^T V_j = ||V_i||_2^2 + ||V_j||_2^2  = 2$.
Therefore $\lapcoeff{V}{g'_{\tuneparam}(0) \mathbf{1}_D} = g'_{\tuneparam}(0) \lapcoeff{V}{\mathbf{1}_D} =  g'_{\tuneparam}(0) \mathbf{1}_D$.

Finally, suppose $\est$ is a fixed point of the LLA algorithm.
This means $\est$ is a minimizer of 
$$
 \underset{x}{\textup{argmin}} \;\; \ell(x)  + Q(x | \est).
$$
The first order necessary conditions for this problem require $ \nabla \ell (\est) + \nabla  Q(\est | \est)$ where the gradient is applied to the first argument of the second term.
We can check that $ \nabla  Q(\est | \est) = \nabla  \fcls{\tuneparam}{\est}$ by construction.
Thus $\est$ satisfies $ \nabla \ell (\est) +  \nabla  \fcls{\tuneparam}{\est}$, which is the condition for being a stationary point of Problem \eqref{prob:con_comp_pen}.

\end{proof}

\subsection{Proofs for Section \ref{s:theory}}

The proof of Theorems \ref{thm:orc_is_stat_point} and \ref{thm:init_from_close_give_lasso_orc} are based on the proofs of Theorem 1 and 2 of \cite{fan2014strong}.

\begin{proof} of Theorem \ref{thm:orc_is_stat_point}

The assumptions on $\estzero$ imply that $\mathcal{A}(\estzero)$ has exactly $\nccTarg$ connected components and these connected components are exactly the connected components of $\esttarg$.
Let  $\lambda \in \mathbb{R}^{d}_+$ be the eigenvalues of $\mathcal{L}(|\estzero|)$ and let $V \in \mathbb{R}^{d \times d}$ be any orthogonal matrix of eigenvectors.
Also let $w = g_{\tuneparam}'(\lambda) \in \mathbb{R}^{d}_+$.
Assumption \eqref{eq:orc_is_stat_point__big_eval_compwise} further implies that  $\lambda_{(\nccTarg + 1)} \ge b_2 \tuneparam$. 
Thus by Definition \ref{def:scad_like_pen_func}
$$
g_{\lambda}'(\lambda_k) =  a_0\tuneparam , \text{ for }  1 \le k \le \nccTarg, \qquad g_{\lambda}'(\lambda_k) =0, \text{ for }  k \ge \nccTarg + 1,
$$
and we can check
$$
\mathcal{M}(V, w)  = \mathcal{M}(V_{(1:\nccTarg)}, w_{(1:\nccTarg)}) = \mathcal{M}(V_{(1:\nccTarg)}, a_0 \mathbf{1}_{\nccTarg})   =:  M \in \mathbb{R}^D.
$$

Now let $\estone$ be the output of taking one LLA step from $\estzero$.
Then $\estone$ is a solution to 
\begin{equation} \label{eq:est_orc_stat_lla_one_lasso}
\underset{\est}{\textup{minimize}} \;\; \ell(\est)  + \frac{1}{2}  M^T |\est| 
\end{equation}
and let $L(\cdot)$ be the loss function for this problem.

Recalling that $\estzero$ has the same $\nccTarg$ connected components as $\esttarg$ and applying Proposition \ref{prop:lap_coef_with_con_comps} we see
\begin{equation} \label{eq:orc_is_stat_pt_lap_coef}
M_{\suppSetTarg}  = 0, \text{ and }
M_{\suppSetTargPow{C}} = \frac{2 a_0 \tuneparam}{\maxccTarg}.
\end{equation}

Recall $\estorc$ is the solution to the constrained Problem \eqref{prob:block_oracle}. 
Therefore for any $\est \in \mathbb{R}^D$,
\begin{equation} \label{eq:loss_func_ineq_for_orc}
\begin{aligned}
\ell(\est) - \ell(\estorc) 
& \ge \nabla \ell(\estorc)^T (\est - \estorc) \\
& = \nabla_{\suppSetTarg} \ell(\estorc)^T (\est_{\suppSetTarg} - \estorc_{\suppSetTarg}) +  \nabla_{ \suppSetTargPow{C}} \ell(\estorc)^T (\est_{ \suppSetTargPow{C}} - \estorc_{ \suppSetTargPow{C}}) \\
& = \nabla_{ \suppSetTargPow{C}} \ell(\estorc)^T \est_{ \suppSetTargPow{C}}
\end{aligned}
\end{equation}
where the inequality comes from the convexity of $\ell$ and the second equality comes from the facts $\estorc_{\suppSetTargPow{C}} = 0$ and $\nabla_{\suppSetTarg} \ell(\estorc) = 0$ comes from the first order necessary conditions for Problem \eqref{prob:block_oracle}. 

Then for any $\est \in \mathbb{R}^D$,
\begin{equation*}
\begin{aligned}
L(\est) - L(\estorc) 
& = \ell(\est) - \ell(\estorc) +  \frac{1}{2}  M^T \left( |\est| - |\estorc|  \right)\\
& = \ell(\est) - \ell(\estorc) +  \frac{1}{2}  M_{\suppSetTargPow{C}}^T  |\est_{\suppSetTargPow{C}}|  \\
& \ge  \nabla_{\suppSetTargPow{C}} \ell(\estorc) \est_{\suppSetTargPow{C}} +  \frac{1}{2} M_{\suppSetTargPow{C}}^T  |\est_{\suppSetTargPow{C}}|  \\
& =  \left(\nabla_{\suppSetTargPow{C}} \ell(\estorc) \text{sign}(\est_{\suppSetTargPow{C}})  +  \frac{1}{2} M_{\suppSetTargPow{C}}\right)^T |\est_{\suppSetTargPow{C}}|  \\
& =  \left(  \frac{1}{2} M_{\suppSetTargPow{C}} - |\nabla_{\suppSetTargPow{C}} \ell(\estorc)| \right)^T |\est_{\suppSetTargPow{C}}|  \\
& =  \left( \frac{ a_0 \tuneparam}{\maxccTarg} - |\nabla_{\suppSetTargPow{C}} \ell(\estorc)| \right)^T |\est_{\suppSetTargPow{C}}|  \\
& \ge 0
\end{aligned}
\end{equation*}
The second equality uses the facts $M_{\suppSetTarg} = 0$ by \eqref{eq:orc_is_stat_pt_lap_coef} and $\estorc_{\suppSetTargPow{C}} = 0$ by construction.
The first inequality uses \eqref{eq:loss_func_ineq_for_orc}.
The final equality comes from  \eqref{eq:orc_is_stat_pt_lap_coef}.
The final inequality comes from Assumption \eqref{eq:orc_has_nice_grad}; this inequality is strict unless $\est_{\suppSetTargPow{C}} = 0$.
Putting this together with the uniqueness of Problem \eqref{prob:block_oracle} we conclude that $\estone = \estorc$.

The claim that fixed points are stationary points follows from Proposition \ref{prop:spect_pen_maj}.

\end{proof}

\begin{proof} of Theorem \ref{thm:init_from_close_give_lasso_orc}

Let  $\lambda \in \mathbb{R}^{d}_+$ be the eigenvalues of $\mathcal{L}(|\estzero|)$ and let $V \in \mathbb{R}^{d \times d}$ be any orthogonal matrix of eigenvectors.
Also let $w = g_{\tuneparam}'(\lambda) \in \mathbb{R}^{d}_+$.
Assumption \eqref{eq:init_from_close_give_lasso_orc__big_gap} says that  $\lambda_{(\nccTarg)} \le b_1 \tuneparam$ and $\lambda_{(\nccTarg + 1)} \ge b_2 \tuneparam$.
Thus by Definition \ref{def:scad_like_pen_func}
\begin{equation} \label{eq:any_init_grad_bounds}
g_{\lambda}'(\lambda_k) \in  [a_1, a_0]\tuneparam , \text{ for }  1 \le k \le \nccTarg, \qquad g_{\lambda}'(\lambda_k) =0, \text{ for }  k \ge \nccTarg + 1
\end{equation}
and we can check
$$
\mathcal{M}(V, w)  = \mathcal{M}(V_{(1:\nccTarg)}, w_{(1:\nccTarg)}) = \mathcal{M}(V_{(1:\nccTarg)}, a_0 \mathbf{1}_{\nccTarg})   =:  M \in \mathbb{R}^D_+
$$
where by construction  $\estone$ is a solution to 
\begin{equation} \label{eq:any_init_lla_one_lasso}
\underset{\est}{\textup{minimize}} \;\; \ell(\est)  + \frac{1}{2} M^T |\est|.
\end{equation}

Note Assumption \eqref{eq:init_from_close_give_lasso_orc__small_ratio__op} guarantees $\opn{\mathcal{L}(|\estzero| - |\esttarg|)}  \le \frac{1}{2}\targgap$ thus
\begin{equation} \label{eq:lasso_weight_lbd_for_between_block_edges}
\begin{aligned}
||M_{\suppSetTargPow{C}}||_{\text{min}}
& \ge \left( \sqrt{\frac{2 \min w^{(0)}_{1:\nccTarg}}{\maxccTarg}} -  2^{5/2} \sqrt{\max w^{(0)}_{1:\nccTarg}} \frac{\opn{\mathcal{L}(|\estzero| - |\esttarg|)}} {\targgap}  \right)^2 \\
& \ge \left( \sqrt{\frac{2 a_1 \tuneparam}{\maxccTarg} } -  2^{5/2} \sqrt{a_0 \tuneparam} \frac{\opn{\mathcal{L}(|\estzero| - |\esttarg|)}}{\targgap}  \right)^2 \\
& =  \frac{2 a_1 \tuneparam}{\maxccTarg}  \left( 1 -  2^{5} \sqrt{\frac{a_0 \maxccTarg}{a_1 }} \frac{\opn{\mathcal{L}(|\estzero| - |\esttarg|)}}{\targgap}  \right)^2\\
& \ge \frac{a_1 \tuneparam}{2\maxccTarg}
\end{aligned}
\end{equation}
where the first inequality comes from \eqref{eq:lasso_maj_weight_bounds_from_opn} of Lemma \ref{lem:lasso_maj_weight_bounds},
the second inequality comes from \eqref{eq:any_init_grad_bounds},
and the final inequality comes from Assumption  \eqref{eq:init_from_close_give_lasso_orc__small_ratio__op}.
Similarly,
\begin{equation*}
\begin{aligned}
||M_{\suppSetTarg}||_{\text{max}}
 \le  2^{5} \max(w)  \left( \frac{\opn{\mathcal{L}(|\estzero| - |\esttarg|)}}{\triangle}  \right)^2 
\le 2^{5} a_0 \tuneparam \residDiffGapRatio^2
\end{aligned}
\end{equation*}
where we have again used Assumption  \eqref{eq:init_from_close_give_lasso_orc__small_ratio__op} and \eqref{eq:lasso_maj_weight_bounds_from_opn}.

Let $\estOrcLasso$  be the solution to the Lasso oracle problem,
\begin{equation} \label{eq:any_init_lasso_orc_prob}
\begin{aligned}
& \underset{\est \in \mathbb{R}^{D} }{\textup{minimize}}  & & \ell(\est) + \frac{1}{2} M_{\suppSetTarg}^T |\est_{\suppSetTarg}| \\ 
& \text{subject to } & & \est_{\suppSetTargPow{C}} = 0.
\end{aligned}
\end{equation}
The above discussion shows $\estOrcLasso \in \lassoOrcSetOp{\residDiffGapRatio, \tuneparam}$ under Assumption \eqref{eq:init_from_close_give_lasso_orc__small_ratio__op}.
We next show $\estone = \estOrcLasso$.

Let $L(\est)$ be the loss function for Problem \eqref{eq:any_init_lla_one_lasso} and let $L^{\suppSetTarg}(\est)$ be the loss function for \eqref{eq:any_init_lasso_orc_prob}.
By convexity for any $\est \in \mathbb{R}^D$,
\begin{equation*}
\begin{aligned}
L^{\suppSetTarg}(\est) - L^{\suppSetTarg}(\estOrcLasso)
& \ge \nabla L^{\suppSetTarg}(\estOrcLasso)^T \left( \est - \estOrcLasso \right) \\
& = \nabla_{\suppSetTarg} L^{\suppSetTarg}(\estOrcLasso)^T \left( \est_{\suppSetTarg} - \estOrcLasso_{\suppSetTarg} \right)   \\
& \;\;\;\; +  \nabla_{\suppSetTargPow{C}} L^{\suppSetTarg}(\estOrcLasso)^T \left( \est_{\suppSetTargPow{C}} - \estOrcLasso_{\suppSetTargPow{C}} \right)\\
& =  \nabla_{\suppSetTargPow{C}} L^{\suppSetTarg}(\estOrcLasso)^T \left( \est_{\suppSetTargPow{C}} - \estOrcLasso_{\suppSetTargPow{C}} \right)\\
& =  \nabla_{\suppSetTargPow{C}} L^{\suppSetTarg}(\estOrcLasso)^T \est_{\suppSetTargPow{C}},
\end{aligned}
\end{equation*}
where we have use the fact that $\nabla_{\suppSetTarg} L^{\suppSetTarg}(\estOrcLasso) = 0$ by the first order necessary conditions of \eqref{eq:any_init_lasso_orc_prob} and $\estOrcLasso_{\suppSetTargPow{C}} = 0$ by construction.
Therefore
\begin{equation*}
\begin{aligned}
L(\est) - L(\estOrcLasso) 
& = L^{\suppSetTarg}(\est) - L^{\suppSetTarg}(\estOrcLasso) +  \frac{1}{2}M_{ \suppSetTargPow{C}}^T \left(|\est_{ \suppSetTargPow{C}}| - |\estOrcLasso_{ \suppSetTargPow{C}}| \right)\\
& = L^{\suppSetTarg}(\est) - L^{\suppSetTarg}(\estOrcLasso) + \frac{1}{2} M_{ \suppSetTargPow{C}}^T |\est_{ \suppSetTargPow{C}}| \\
& \ge  \nabla_{\suppSetTargPow{C}} L^{\suppSetTarg}(\estOrcLasso)^T \est_{\suppSetTargPow{C}} +\frac{1}{2}  M_{ \suppSetTargPow{C}}^T |\est_{ \suppSetTargPow{C}}| \\
&  =  \left( \nabla_{\suppSetTargPow{C}} L^{\suppSetTarg}(\estOrcLasso) \text{sign} (\est_{\suppSetTargPow{C}} )+ \frac{1}{2} M_{ \suppSetTargPow{C}}\right)^T |\est_{ \suppSetTargPow{C}}| \\
&  \ge  \left( \frac{1}{2} M_{ \suppSetTargPow{C}} - | \nabla_{\suppSetTargPow{C}}L^{\suppSetTarg}(\estOrcLasso)|\right)^T |\est_{ \suppSetTargPow{C}}| \\
&  \ge  \left( \frac{a_1 \tuneparam}{4 \maxccTarg} - | \nabla_{\suppSetTargPow{C}}L^{\suppSetTarg}(\estOrcLasso)|\right)^T |\est_{ \suppSetTargPow{C}}| \\
&  =  \left( \frac{a_1 \tuneparam}{4 \maxccTarg} - | \nabla_{\suppSetTargPow{C}}\ell(\estOrcLasso)|\right)^T |\est_{ \suppSetTargPow{C}}| \\
& \ge 0
\end{aligned}
\end{equation*}
where the penultimate inequality comes from \eqref{eq:lasso_weight_lbd_for_between_block_edges}.
The final inequality follows from Assumption \eqref{eq:init_from_close_give_lasso_orc__nice_grad} and the fact that  $\estOrcLasso \in \lassoOrcSetOp{\residDiffGapRatio, \tuneparam}$; this inequality is strict unless $\est_{ \suppSetTargPow{C} } = 0$.
Putting this together with the uniqueness of \eqref{eq:any_init_lasso_orc_prob} we conclude $\estone = \estOrcLasso$.

Next we show taking one more step results in the block oracle estimator.
By Assumption \eqref{eq:init_from_close_give_lasso_orc__nice_grad}
$$
\estorc \in \lassoOrcSetOp{0, \tuneparam} \subseteq \lassoOrcSetOp{\residDiffGapRatio, \tuneparam} \subseteq \setNiceGrad{\frac{a_1 \tuneparam}{4 \maxccTarg}}  \subseteq \setNiceGrad{a_0 \tuneparam},
$$
so the condition \eqref{eq:orc_has_nice_grad} is satisfied. 
By Assumption \eqref{eq:init_from_close_give_lasso_orc__lasso_orc_uniformly_big_eval}
$$
\estorc \in \lassoOrcSetOp{0, \tuneparam} \subseteq  \lassoOrcSetOp{\residDiffGapRatio, \tuneparam}  \subseteq \setBigEvalComp{b_2 \tuneparam}
$$
so the condition \eqref{eq:orc_is_stat_point__orc_in_set_for_stat_point} is satisfied.
Thus by the second claim of Theorem \ref{thm:orc_is_stat_point} we conclude that $\estorc$ is a fixed point of the LLA algorithm.

By construction $ \blocksupp{\mathcal{A}(\est)} \subseteq  \blocksupp{\mathcal{A}(\esttarg)}$ for any $\est \in  \lassoOrcSetOp{\residDiffGapRatio, \tuneparam}$.
Therefore the fact that $\estone \in \lassoOrcSet{\residDiffGapRatio, \tuneparam}$ and Assumption \eqref{eq:init_from_close_give_lasso_orc__lasso_orc_uniformly_big_eval}  together imply that the conditions of the first claim of Theorem \ref{thm:orc_is_stat_point} are satisfied.
Therefore taking a LLA step from $\estone$ results in $\estorc$, which is a fixed point of the LLA algorithm. 
Thus all the claims of the theorem follow except the final Frobenius norm bound.

If Assumption \eqref{eq:init_from_close_give_lasso_orc__small_ratio__op} is replaced with \eqref{eq:init_from_close_give_lasso_orc__small_ratio__op} then using \eqref{eq:lasso_maj_weight_bounds_from_frob} of  Lemma \ref{lem:lasso_maj_weight_bounds} we obtain the same bounds on $||M_{\suppSetTargPow{C}}||_{\text{min}}$ and $||M_{\suppSetTarg}||_{\text{max}}$ as above.
Let $\frobBound$ be the left hand side of \eqref{eq:init_from_close_give_lasso_orc__small_ratio__op}.
Similarly, applying the $L_1$ and $L_2$ claims in \eqref{eq:lasso_maj_weight_bounds_from_frob} we obtain
\begin{equation*}
\begin{aligned}
||M_{\suppSetTarg}||_{1}
& \le  2^4 d_{\text{max}} \max(w) \left( \frac{\frobBound}{\triangle} \right)^2 
\le  2^4 d_{\text{max}} a_0 \residDiffGapRatio^2 \\
||M_{\suppSetTarg}||_{2}
& \le  2^{9/2} \maxccTargPow{1/2} \max(w) \left( \frac{\frobBound}{\triangle} \right)^2 
\le  2^{9/2}  \maxccTargPow{1/2}a_0 \residDiffGapRatio^2.
\end{aligned}
\end{equation*}
The remainder of the proof is the same thus the claim about the Frobenius norm bound follows.

\end{proof}

\begin{proof} of Corollary \ref{cor:suff_cond_whp_for_init_from_close_give_lasso_orc}
We verify the conditions of Theorem \ref{thm:init_from_close_give_lasso_orc}  are satisfied.
Under the (complement of) the event in $\probGoodInitLLA$ we have
$$
\opn{\mathcal{L}(|\estinit| - |\esttarg|)}
 \le \min \left(\residDiffGapRatio,  \frac{1}{2^6} \sqrt{\frac{a_1}{a_0 \maxccTarg}}\right)\targgap \wedge b_1 \tuneparam
$$
where we have used the assumption $\targgap \ge (b_1 + b_2) \tuneparam$.
The first term on the right hand side shows condition \eqref{eq:init_from_close_give_lasso_orc__small_ratio__op} holds.
Using the second term on the right hand side and Weyl's inequality we see
$$
|\lapSmEval{j}{|\estinit|} - \lapSmEval{j}{|\estorc|} | \le b_1 \tuneparam
$$
for each $j$.
Setting $j = \nccTarg$ we get
$$
\lapSmEval{\nccTarg}{|\estinit|} \le b_1 \tuneparam
$$
and setting $j =\nccTarg + 1$ we get
$$
\lapSmEval{\nccTarg + 1}{|\estinit|} \ge \targgap -  b_1 \tuneparam \ge (b_1 + b_2) \tuneparam - b_1 \tuneparam  = b_2 \tuneparam.
$$
Thus condition \eqref{eq:init_from_close_give_lasso_orc__big_gap} holds.
The remaining conditions of the Theorem are controlled by $\probNiceGradLLA$ and $\probSmallResidLLA$.

The Frobenius norm claims follow analogously.
Note condition  \eqref{eq:init_from_close_give_lasso_orc__big_gap} follows as above since the Frobenius norm upper bounds the operator norm.
\end{proof}

\section{Proofs for Section \ref{s:bd_shrink}}  \label{a:proofs__bd_shrink}

\subsection{Preliminary results}\label{as:proofs__bd_shrink__prelim}

We first give concentration bounds on $\opn{\mathcal{L}(|X|)}$ and $\opn{\mathcal{L}(X)}$ where $X$ is a random vector with sub-Gaussian entries.
For the next two lemmas we assume the entries of $X \in \mathbb{R}^D$ are $\sigma$ sub-Gaussian.
The first lemma is based on bounding  $\prob{ \opone{\mathcal{A}(|X|)} \ge t}$.
\begin{lemma} \label{lem:op_bound_via_opone}
Let $\opoEBoundAbs := \opone{\mathcal{A}(\expect{|X|})} \le d || \expect{|X|}||_{\text{max}} \le c d \sigma$.
Then
\begin{equation*}
\prob{ \opn{\mathcal{L}(|X|)}  \ge t} 
\le  2 d \exp \left( - \frac{[t - \opoEBoundAbs ]_+^2}{C d^2 \sigma^2} \right)
%
\end{equation*}
and $\expect{\mathcal{L}(|X|)}  \le c  \sigma d \sqrt{\log d}$. 
If the entries of $X$ are independent then
\begin{equation*}
\prob{ \opn{\mathcal{L}(|X|)}  \ge t} 
\le  2 d \exp \left( - \frac{[t - \opoEBoundAbs ]_+^2}{C d \sigma^2} \right)
%
\end{equation*}
and if   $|| \expect{|X|}||_{\text{min}} \ge c \sigma$ then  $\expect{\mathcal{L}(|X|)}  = C d  \sigma$.
\end{lemma}

The next lemma is based on bounding $\prob{ \opn{\mathcal{A}(X)}  + ||\mathcal{A}(X) \mathbf{1}||_{\text{max}} \ge t}$.
\begin{lemma} \label{lem:op_bound_via_deg_plus_op}
Assume the entries of $X \in \mathbb{R}^D$ are independent.  
Let $\opoEBound := \opone{\mathcal{A}(\expect{X})} \le d || \expect{X}||_{\text{max}} \le c d \sigma$.
Then
\begin{equation*}
\prob{ \opn{\mathcal{L}(X)}  \ge t} 
 \le d \exp\left(-\frac{ [\frac{t}{2} - \opoEBound  ]^2}{d c_1 \sigma^2} \right)  
+ 4 \exp \left(- \frac{1}{c_4^2 \sigma^2 } \left[ \frac{t}{2} - \opoEBound - c_4 \sigma \sqrt{d} \right]_+^2 \right)
\end{equation*}
and $\expect{ \opn{\mathcal{L}( X) }} \le 2 \opoEBound + c \sigma \sqrt{d\log(d)} $.

If $\mathcal{A}(X)$ has $\nccNzTarg$ non-trivial connected components, $\NumNonIso$ non-isolated vertices and the largest connected component has $\maxccTarg$ vertices then  $\opoEBound \le \maxccTarg || \expect{X}||_{\text{max}} \le c \maxccTarg \sigma$,
\begin{equation*}
\prob{ \opn{\mathcal{L}(X)}  \ge t} 
 \le \NumNonIso \exp\left(-\frac{ [\frac{t}{2} - \opoEBound  ]^2}{\maxccTarg c_1 \sigma^2} \right)  
+ 4  \nccNzTarg \exp \left(- \frac{1}{c_4^2 \sigma^2 } \left[ \frac{t}{2} - \opoEBound - c_4 \sigma \sqrt{\maxccTarg} \right]_+^2 \right) 
\end{equation*}
and $\expect{ \opn{\mathcal{L}( X})} \le 2 \opoEBound + c \sigma \sqrt{\maxccTarg \log(\NumNonIso) }$.
\end{lemma}

The next claim improves upon the naive bound $\expect{|x| - |\mu|} \le c \sigma$ that comes from standard sub-Gaussian properties.
\begin{proposition} \label{prop:abs_resid_sub_g}
Let  $\mu := \expect{x}$ and suppose $x - \mu$ is a $\sigma$ sub-Gaussian random variable.
Then $|x| - |\mu|$ is $\sigma$ sub-Gaussian and 
\begin{equation} \label{eq:abs_resid_sub_g_expect_bound}
\expect{|x| - |\mu|} \le C( \sigma + |\mu|) \exp \left( -\frac{\mu^2}{4 \sigma^2}\right).
\end{equation}
\end{proposition} 




The following lemma characterizes the error after applying a generalized thresholded operator to $\estok$.
Recall from Section \ref{s:bd_shrink} $\maxDegTarg$ is the maximal degree of a node in the binary target $\binGraphTarg$.
Also let $\nnzSizeTarg$  be the number of non-zero elements of $\esttarg$ i.e. the number of edges in $\binGraphTarg$.

\begin{lemma} \label{lem:gen_thresh_tail_cond_init_resid}
Suppose $\estok \in \mathbb{R}^D$ satisfies a $\tailclass{f}{\tailconst}$ residual tail condition.
Then for any $\ThreshTuneParam \in (0,  \frac{1}{\tailconst})$ and $\alpha \in (0, 1)$,
\begin{equation}\label{eq:gen_thresh_tail_cond_init_resid_L1_bound}
\prob{||\threshop{\estok}{\ThreshTuneParam} - \esttarg||_1 > 6 \nnzSizeTarg \ThreshTuneParam }  \le 
 \nnzSizeTarg f\left(n, \ThreshTuneParam \right) + D f\left(n, (1- \alpha) \ThreshTuneParam \right),
\end{equation}
and
\begin{equation}\label{eq:gen_thresh_tail_condition_init_resid_opone_bound}
\prob{ \opone{\mathcal{A} \left(\threshop{\estok}{\ThreshTuneParam}  - \esttarg\right)}>  6 \maxDegTarg \ThreshTuneParam } \le 
 \nnzSizeTarg f\left(n, \ThreshTuneParam \right) + D f\left(n, (1- \alpha) \ThreshTuneParam \right).
 \end{equation}
\end{lemma}

The following are sufficient conditions for hard thresholding to return the oracle estimate.
Let $\estorc_{\text{entrywise}}$ be the entrywise oracle for the target parameter $\esttarg$ i.e.  the solution to \eqref{prob:block_oracle} where $\suppSetTarg$ is replaced with $\nnzSetTarg = \{j \in [D] \text{ s.t. } \esttarg_j \neq 0\}$.
\begin{proposition} \label{prop_hard_thresh}
Suppose $\estok$ in Section obeys a $\sigma$ sub-Gaussian tail residual condition.
Then $\estHT{\ThreshTuneParam} = \estOrcEntry$ with probability at least $1 - \probHTBigTrue -\probHTSmallNoise$ where
$$
\probHTBigTrue 
= \prob{||\estok_{\nnzSetTarg} - \esttarg_{\nnzSetTarg} ||_{\text{max}} 
\ge \minNnzMagTarg- \ThreshTuneParam }
\le 2 \nnzSizeTarg \exp \left(\frac{- n \left[\minNnzMagTarg - \ThreshTuneParam \right]_+^2 }{2 \sigma^2} \right)
$$
and
$$
\probHTSmallNoise 
= \prob{||\estok_{\nnzSetTargPow{C}}||_{\text{max}} \le \ThreshTuneParam}
\le 2 (D - \nnzSizeTarg) \exp\left( \frac{-n \ThreshTuneParam^2}{2 \sigma^2} \right).
$$
\end{proposition}

\subsection{Proofs for main results}\label{as:proofs__bd_shrink__main}

\begin{proof} of Theorem \ref{thm:bd_shrink__stat_pt__nice} 

We will appeal to Corollary \ref{cor:orc_is_stat_pt_whp}.
The objective function of Problem \eqref{prob:prox_fclsp},  $\ell(\est) = \frac{1}{2} || \estok - \est||_2^2$, is strongly convex so the solution is always unique.
Thus Assumption \ref{assu:basic} is satisfied.

Note $\nabla \ell(\est) = \estok - \est$ and $\estorc_{\suppSetTargPow{C}} = 0$ thus $\nabla_{\suppSetTargPow{C}} \ell(\estorc) = \estok_{\suppSetTargPow{C}}$.
Therefore by the union bound
\begin{equation}
\prob{ ||\nabla_{\suppSetTargPow{C}} \ell(\estorc) ||_{\text{max}} > t}  
=  \prob{||\estok_{\suppSetTargPow{C}} ||_{\text{max}} > t}  
\le 2 (D - \suppSizeTarg) \exp \left( - \frac{n t^2}{2 \sigma^2} \right),
\end{equation}
thus the claim about $\probNiceGradStatPt$ follows by plugging in $t = \frac{a_0 \tuneparam}{\maxccTarg}$.

Note $\estorc_{ \ccTarg{k}} = \estok_{ \ccTarg{k}}$.
By assumption $\expect{\estorc_{ \ccTarg{k}}} = \esttarg_{ \ccTarg{k}}$ for any  $k \in [\nccNzTarg]$ and the entries of $\estorc_{\suppSetTarg}$ are independent and have sub-Gaussian variance proxy $\frac{\sigma^2}{n}$.
Therefore
\begin{equation} \label{eq:fuzzy_sleepy_cat}
\begin{aligned}
\opone{\mathcal{A}(\expect{|\estorc|} - |\esttarg |)} 
&  \le  \maxccTarg ||\expect{|\estorc|} - |\esttarg |||_{\text{max}}  \\
& \le C\maxccTarg ( \frac{\sigma}{\sqrt{n}} + ||\esttarg||_{\text{max}}) \exp \left( -\frac{c n  \minSuppMagTarg^2 }{ \sigma^2}\right)
=: \subGResidExpectBound
\end{aligned}
\end{equation}
where the first inequality comes from the fact that the row sums have at most $\maxccTarg$ non-zero terms 
and the second inequality comes from  Proposition \ref{prop:abs_resid_sub_g}.
Thus by Lemma \ref{lem:op_bound_via_deg_plus_op} we have
\begin{equation*}
\begin{aligned}
\prob{\opn{\mathcal{L}(|\estorc| - |\esttarg |)} \ge t}
 \le \NumNonIso \exp\left(-\frac{ n [\frac{t}{2} - \subGResidExpectBound  ]^2}{\maxccTarg c_1 \sigma^2} \right)  
+ 4  \nccNzTarg \exp \left(- \frac{n}{c_4^2 \sigma^2 } \left[ \frac{t}{2} - \subGResidExpectBound - c_4 \sigma \sqrt{\frac{\maxccTarg}{n}} \right]_+^2 \right) 
\end{aligned}
\end{equation*}
From the assumptions of the theorem $\targgap - b_2 \tuneparam  \ge b_2 \tuneparam$.
Therefore plugging $t = b_2 \tuneparam$ and using Assumption \eqref{eq:bd_shrink__indep_unbiased__stat_pt__targ_gap__assumpt} the claim about $\probSmallResidStatPt$ follows.


\end{proof}

\begin{proof} of Theorem \ref{thm:bd_shrink__two_step_general}
To verify the two step convergence claims we will appeal to Corollary \ref{cor:suff_cond_whp_for_init_from_close_give_lasso_orc}.

For any $\est \in \lassoOrcSetOp{\residDiffGapRatio, \tuneparam}$, $\est_{\suppSetTargPow{C}} = 0$ so  $\nabla_{\suppSetTargPow{C}} \ell(\est) = \estok_{\suppSetTargPow{C}}$.
Therefore by the union bound
\begin{equation}
\prob{ \sup_{\est \in \lassoOrcSetOp{\residDiffGapRatio, \tuneparam}} ||\nabla_{\suppSetTargPow{C}} \ell(\est) ||_{\text{max}} > t}  
=  \prob{||\estok_{\suppSetTargPow{C}} ||_{\text{max}} > t}  
\le  (D - \suppSizeTarg) f(n, t)
\end{equation}
 for any   $t \in (0, 1/\tailconst)$.
Thus the claim about $\probNiceGradStatPt$ follows after plugging in $t = \frac{a_1 \tuneparam}{4\maxccTarg}$.

Note the entries of any $\est \in \lassoOrcSetOp{\residDiffGapRatio, \tuneparam}$ are obtained by soft-thresholding the entries of $\estok$ so $||\est_{\suppSetTarg} - \estok_{\suppSetTarg}||_{\text{max}}  \le  2^5 a_0 \tuneparam \residDiffGapRatio^2$ by the definition of $ \lassoOrcSetOp{\residDiffGapRatio, \tuneparam}$.
Thus
$$
\sup_{\est \in \lassoOrcSet{\residDiffGapRatio, \tuneparam}}  ||\est_{\suppSetTarg} - \estok_{\suppSetTarg}||_{\text{max}}  
\le 2^5 a_0 \tuneparam \residDiffGapRatio^2
$$
and we see
\begin{equation*}
\begin{aligned}
\prob{\sup_{\est \in \lassoOrcSet{\residDiffGapRatio, \tuneparam}}  \opn{\mathcal{L}(|\est| - |\esttarg| ) } \ge t} 
& \le \prob{  \opn{\mathcal{L}(|\estorc| - |\esttarg| ) } \ge t - \sup_{\est \in \lassoOrcSet{\residDiffGapRatio, \tuneparam}}  \opn{\mathcal{L}(|\est| - |\estorc| ) } }  \\
& \le \prob{  \opn{\mathcal{L}(|\estorc| - |\esttarg| ) } \ge t - \sup_{\est \in \lassoOrcSet{\residDiffGapRatio, \tuneparam}} 2  \maxccTarg  ||\est - \estorc ||_{\text{max}}  }  \\
& \le \prob{  \opn{\mathcal{L}(|\estorc| - |\esttarg| ) } \ge t -  \maxccTarg 2^6 a_0 \tuneparam \residDiffGapRatio^2 }  \\
& \le \prob{  ||\estorc - \esttarg ||_{\text{max}} \ge \frac{t -  \maxccTarg 2^6 a_0 \tuneparam \residDiffGapRatio^2}{2 \maxccTarg} }  \\
& \le D f(n, \frac{t -  \maxccTarg 2^6 a_0 \tuneparam \residDiffGapRatio^2}{2 \maxccTarg})
\end{aligned}
\end{equation*}
where the first inequality uses the triangle inequality,
the second and fourth inequality bound the operator norm using the max norm (Proposition \ref{prop:lap_bounds}),
and the final inequality uses the union bound.
Setting $t = b_1 \tuneparam$ and applying Assumption \eqref{eq:thoe_ex__bd_shrink__reg_sat_whp__init_tol}  the claim about $\probSmallResidLLA$ follows.

For the claim in Remark \ref{rem:two_step_reates_under_indep_unbaised} we can argue as in the proof of Theorem \ref{thm:bd_shrink__stat_pt__nice} to bound $ \opn{\mathcal{L}(|\estorc| - |\esttarg| ) } $.

\end{proof}


\begin{proof} of Corollary \ref{cor:bd_shrink__initializer__general} 

We need only verify the initialization condition of Corollary \ref{cor:suff_cond_whp_for_init_from_close_give_lasso_orc}.
Note
\begin{equation*}
\begin{aligned}
\prob{\opn{ \mathcal{L} \left(|\estok| - |\esttarg| \right)} \ge t} 
& \le \prob{||\estok -  \esttarg ||_{\text{max}} \ge \frac{t}{2 d}}
& \le D f(n, \frac{t}{2 d}).
\end{aligned}
\end{equation*}
where the first inequality comes from Proposition \ref{prop:lap_bounds}.
Thus the claim \eqref{eq:bd_shrink_est_ok_init} follows by setting $t =\initRatioConsts \tuneparam $.
For the claim in Remark \ref{rem:two_step_reates_under_indep_unbaised} we can improve this rate by using Lemma \ref{lem:op_bound_via_deg_plus_op} to control the within block terms and Lemma \ref{lem:op_bound_via_opone} to control the between block terms.

If $\ThreshTuneParam \le  \frac{1}{2} \cdot \frac{1}{6 \nnzSizeTarg} \cdot  \initRatioConsts \tuneparam$ then 
\begin{equation*}
\begin{aligned}
\prob{\opn{ \mathcal{L} \left(|\threshop{\estok}{\ThreshTuneParam}| - |\esttarg| \right)} \ge  \initRatioConsts \tuneparam} 
& \le \prob{\opone{ \mathcal{A} \left(\threshop{\estok}{\ThreshTuneParam} - \esttarg \right)} \ge  \frac{1}{2} \initRatioConsts \tuneparam}  \\
& \le \prob{\opone{ \mathcal{A} \left(\threshop{\estok}{\ThreshTuneParam} - \esttarg \right)}  \ge  6 \nnzSizeTarg \ThreshTuneParam }\\
& \le \text{right hand side of \eqref{eq:bd_shrink_gen_thresh_init}}
\end{aligned}
\end{equation*}
where the first inequality comes Proposition \ref{prop:lap_bounds} and the final inequality comes from the operator one norm bound of Lemma \eqref{lem:gen_thresh_tail_cond_init_resid}.
Thus \eqref{eq:bd_shrink_gen_thresh_init} follows.
Note we can bound $\prob{ ||\mathcal{L} \left(|\estok| - |\esttarg| \right)||_F \ge t} $ similarly by making use of the $L_1$ bound in Lemma \eqref{lem:gen_thresh_tail_cond_init_resid}.

\end{proof}

\begin{remark}
The claims in Remark \ref{rem:two_step_reates_under_indep_unbaised} follow by enlisting Lemmas \ref{lem:op_bound_via_opone} and \ref{lem:op_bound_via_deg_plus_op} in the proofs of Theorem \ref{thm:bd_shrink__two_step_general} and Corollary \ref{cor:bd_shrink__initializer__general}.
\end{remark}


\subsection{Proofs for preliminary results}\label{as:proofs__bd_shrink__prelim}

\begin{proof} of Lemma \ref{lem:op_bound_via_opone}
Recall $\opn{\mathcal{L}(|X|)} \le 2 \opone{\mathcal{A}(X)}$ from Proposition \ref{prop:lap_bounds}.

Note
\begin{equation}
\begin{aligned}
\sum_{j=1, j \neq \ell} |X_{(\ell j)}|
\le \sum_{j=1, j \neq \ell} \left( |X_{(\ell j)}| - \expect{ |X_{(\ell j)}| } \right) + \opone{\mathcal{A}(\expect{|X|)}}.
\end{aligned}
\end{equation}
thus setting $R := \mathcal{A}(|X| - \expect{|X|})$ we see that
$$
\opone{\mathcal{A}(|X|)} \le || R \mathbf{1}_d ||_{\text{max}} + \opone{\mathcal{A}(\expect{|X|)}}.
$$
Note the entries of $R$ are mean zero and $\sigma$ sub-Gaussian.
Thus we have
\begin{equation}
\begin{aligned}
\prob{\opone{\mathcal{A}(|X|)}  \ge t}
& \le \prob{ || R \mathbf{1}_d ||_{\text{max}}  \ge t - \opone{\mathcal{A}(\expect{|X|)}}} \\
& \le d \exp\left( - \frac{[ t - \opone{\mathcal{A}(\expect{|X|)}}]_+^2}{c d^2 \sigma^2}\right)
\end{aligned}
\end{equation}
where the second inequality uses the union bound and the fact the row sums of $R$ have sub-Gaussian parameter $c \sigma d$. 
Thus the first claim follows.
The second claim follows analogously after noting that the rows sums of $R$ have sub-Gaussian parameter $c \sqrt{d} \sigma$ by Hoeffding's inequality.

Recall by standard sub-Gaussian properties, $|| \expect{X}||_{\text{max}} \le \widetilde{c} \sigma$.
Integrating the two concentration inequalities gives the stated upper bounds on $ \expect{\mathcal{L}(|X|)}$.
The lower bound in the second claim follows by Jensen's inequality which gives $\expect{\opn{\mathcal{L}(|X|})}\ge\opn{\mathcal{L}(\expect{|X|})} = C \sigma d$.

\end{proof}

\begin{proof} of Lemma \ref{lem:op_bound_via_deg_plus_op}

Recall $\opn{\mathcal{L}(X)} \le || \mathcal{A}(X) \mathbf{1}_d||_{\text{max}} + \opn{\mathcal{A}(X)}$ from Proposition \ref{prop:lap_bounds}.
We have the bound
$$\max(\opn{\expect{\mathbf{A}(X)}}, ||\expect{\mathbf{A}(X) \mathbf{1}_d }||_{\text{max}}) \le \opone{\expect{\mathbf{A}(X)}} :=  \opoEBound.$$

We can control the maximal degree as
\begin{equation*}
\begin{aligned}
\prob{||\mathbf{A}(X) \mathbf{1}_d ||_{\text{max}} \ge t }
& \le \prob{||\mathbf{A}(X - \expect{X}) \mathbf{1}_d ||_{\text{max}} \ge t - \opoEBound  } 
& \le  d \exp\left(\frac{[t -\opoEBound ]_+^2}{c d \sigma^2} \right).
 \end{aligned}
\end{equation*}
The first inequality comes from linearity, the triangle inequality, and the above expectation bound.
The second inequality comes from Hoeffding's bound and the union bound.
Similarly we can control the operator by
\begin{equation*} 
\begin{aligned}
\prob{\opn{\mathcal{A}(X)} \ge t }
& \le \prob{\opn{\mathcal{A}(X - \expect{X})} \ge t -  \opoEBound} \\
& \le 4 \exp \left(-\left[  \frac{t - \opoEBound}{c_5 \sigma} - \sqrt{d} \right]_+^2 \right) 
= 4 \exp \left(- \frac{1}{C \sigma^2} \left[t - \opoEBound- c \sqrt{d} \right]_+^2 \right) 
 \end{aligned}
\end{equation*}
where the second inequality comes from Corollary 4.4.8. of \citep{vershynin2018high}.
Putting these together with the union bound the first claim follows.

The second concentration bound follows similarly by using the fact that the operator norm of a block diagonal matrix is equal to the largest operator norm of the blocks.
Both expectation claims follow by integrating the two concentration bounds.
\end{proof}

\begin{proof} of Proposition \ref{prop:abs_resid_sub_g}

The first claim follows from the reverse triangle inequality e.g. $| |x| - |\mu|| \le |x  - \mu|$.
Now let $y := x - \mu$ so $|x| - |\mu| = |y + \mu| - |\mu|$.
Note $y$ is mean 0 and  $c\sigma$ sub-Gaussian e.g. by Lemma 2.6.8 of \citep{vershynin2018high}.
Let $B(y, \mu) := \{(y, m) \text{ s.t. } |y| \ge |\mu|, \text{sign}(y) \neq \text{sign}(\mu) \}$.
We can check  for any pair of real numbers,
\begin{equation} \label{eq:abs_ym_minus_abs_m_decomp}
\begin{aligned}
|y + \mu| - |\mu| 
& = y \cdot \text{sign}(\mu)  \ind{B(y, \mu)^C } -  \left( y \cdot \text{sign}(\mu) + 2 |\mu| \right) \ind{B(y, \mu) }\\
& = y \cdot \text{sign}(\mu) -  \left( 2 y \cdot \text{sign}(\mu) + 2 |\mu| \right) \ind{B(y, \mu) }
\end{aligned}
\end{equation}
where the first equality comes from writing out cases and the second equality comes from adding 0.
Since $y$ is mean zero the claim follows by bounded the expectation of the second term in \eqref{eq:abs_ym_minus_abs_m_decomp}.
\begin{equation}
\begin{aligned}
\expect{ \left( 2 y \cdot \text{sign}(\mu) + 2 |\mu| \right) \ind{B(y, \mu) }} 
 & \le  2 \expect{ y^2 }^{1/2}  \prob{B(y, \mu)}^{1/2} + 2 |\mu| \prob{B(y, \mu)} \\
& \le  c (\sigma + |\mu| )\prob{B(y, \mu)}^{1/2} \\
& \le  c (\sigma + |\mu| )\prob{|y| \ge |\mu|}^{1/2} \\
& \le c (\sigma + |\mu| ) \exp\left(-\frac{\mu^2}{4 \sigma^2} \right),
\end{aligned}
\end{equation}
where the first inequality uses Cauchy-Schwartz,
the second and fourth inequality use standard sub-Gaussian properties,
and the third inequality follows by expanding the event $B(y, \mu)$.

\end{proof}

\begin{proof} of Proposition \ref{prop_hard_thresh}
The event that all the within block entries are larger than $\ThreshTuneParam$ occurs with probability at least $1 - \probHTBigTrue$ by the union bound and sub-Gaussian assumption.
The event that all the between block entries are smaller than  $\ThreshTuneParam$ occur with probability at least $1 - \probHTSmallNoise$ by the union bound and sub-Gaussian assumption.
\end{proof} 

\begin{proof} of Lemma \ref{lem:gen_thresh_tail_cond_init_resid}

We first prove the $L_1$ bound following an argument of \citep{donoho1994ideal, bickel2008covariance}. 

\begin{equation} \label{lem:gen_thresh_tail_cond_init_resid_1_and_2}
||\threshop{\estok}{\ThreshTuneParam} - \esttarg||_1 \le ||\threshop{\esttarg}{\ThreshTuneParam} - \esttarg||_1 + ||\threshop{\estok}{\ThreshTuneParam} - \threshop{\esttarg}{\ThreshTuneParam}||_1 
\end{equation}

We can control the first term of \eqref{lem:gen_thresh_tail_cond_init_resid_1_and_2} via
\begin{equation} \label{lem:gen_thresh_tail_cond_init_resid_1}
||\threshop{\esttarg}{\ThreshTuneParam} - \esttarg||_1 \le \sum_{\ell=1}^D \ThreshTuneParam \ind{|\esttarg_{\ell}| > \ThreshTuneParam} \le  \nnzSizeTarg \ThreshTuneParam
\end{equation}
which follows by properties of \ref{def:gen_thresh} and the definition of $\nnzSizeTarg$.

Next we break the second term of \eqref{lem:gen_thresh_tail_cond_init_resid_1_and_2} into four terms
\begin{equation} \label{lem:gen_thresh_tail_cond_init_resid_2}
\begin{aligned}
||\threshop{\estok}{\ThreshTuneParam} - \threshop{\esttarg}{\ThreshTuneParam}||_1   \le  \sum_{\ell=1}^D |\threshop{\estok_{\ell}}{\ThreshTuneParam}  - \threshop{\esttarg_{\ell}}{\ThreshTuneParam}| & 
\Big[ 
\ind{|\estok_{\ell}| \ge \ThreshTuneParam, |\esttarg_{\ell}| < \ThreshTuneParam} + 
\ind{|\estok_{\ell}| < \ThreshTuneParam, |\esttarg_{\ell}| \ge \ThreshTuneParam} + \\ 
&  \ind{|\estok_{\ell}| \ge \ThreshTuneParam, |\esttarg_{\ell}| \ge \ThreshTuneParam} + 
\ind{|\estok_{\ell}| < \ThreshTuneParam, |\esttarg_{\ell}| < \ThreshTuneParam} 
 \Big].
\end{aligned}
\end{equation}
The fourth term is 0 by Definition \ref{def:gen_thresh}.

Let $R_{\text{max}} := \sup_{\ell \in [D] \text{s.t. } |\esttarg_{\ell}| > 0 } |\estok_{\ell} - \esttarg_{\ell}|$ be the maximal residual over the non-zero entries of $\esttarg$
The third term of \eqref{lem:gen_thresh_tail_cond_init_resid_2} is bounded as
\begin{equation} \label{lem:gen_thresh_tail_cond_init_resid_2__3}
 \sum_{\ell=1}^D |\threshop{\estok_{\ell}}{\ThreshTuneParam}  - \threshop{\esttarg_{\ell}}{\ThreshTuneParam}| \ind{|\estok_{\ell}| \ge \ThreshTuneParam, |\esttarg_{\ell}| \ge \ThreshTuneParam}   \le
\sum_{\ell=1}^D |\estok_{\ell}  -\esttarg_{\ell}| \ind{|\estok_{\ell}| \ge \ThreshTuneParam, |\esttarg_{\ell}| \ge \ThreshTuneParam}  
 \le \nnzSizeTarg R_{\text{max}},
\end{equation}
by Definition \ref{def:gen_thresh} and construction.

The second term of \eqref{lem:gen_thresh_tail_cond_init_resid_2} is bounded as
\begin{equation} \label{lem:gen_thresh_tail_cond_init_resid_2__2}
\begin{aligned}
 \sum_{\ell=1}^D |\threshop{\estok_{\ell}}{\ThreshTuneParam}  - \threshop{\esttarg_{\ell}}{\ThreshTuneParam}| \ind{|\estok_{\ell}| < \ThreshTuneParam, |\esttarg_{\ell}| \ge \ThreshTuneParam} 
& =  \sum_{\ell=1}^D |\threshop{\esttarg_{\ell}}{\ThreshTuneParam}| \ind{|\estok_{\ell}| < \ThreshTuneParam, |\esttarg_{\ell}| \ge \ThreshTuneParam} \\
& \le  \sum_{\ell=1}^D |\esttarg_{\ell}| \ind{|\estok_{\ell}| < \ThreshTuneParam, |\esttarg_{\ell}| \ge \ThreshTuneParam} \\
& \le  \sum_{\ell=1}^D\left( |\esttarg_{\ell} - \estok_{\ell}| +   |\estok_{\ell}|  \right) \ind{|\estok_{\ell}| < \ThreshTuneParam, |\esttarg_{\ell}| \ge \ThreshTuneParam} \\
%
 & \le (R_{\text{max}} + \ThreshTuneParam) \nnzSizeTarg 
 \end{aligned}
\end{equation}
where the first two claims follow from Definition \ref{def:gen_thresh}.

We bound the first term of \eqref{lem:gen_thresh_tail_cond_init_resid_2} via
\begin{equation} \label{lem:gen_thresh_tail_cond_init_resid_2__1}
\begin{aligned}
\sum_{\ell=1}^D |\threshop{\estok_{\ell}}{\ThreshTuneParam}  - \threshop{\esttarg_{\ell}}{\ThreshTuneParam}| \ind{|\estok_{\ell}| \ge \ThreshTuneParam, |\esttarg_{\ell}| < \ThreshTuneParam} 
&  =  \sum_{\ell=1}^D |\threshop{\estok_{\ell}}{\ThreshTuneParam}| \ind{|\estok_{\ell}| \ge \ThreshTuneParam, |\esttarg_{\ell}| < \ThreshTuneParam} \\
&  \le  \sum_{\ell=1}^D |\estok_{\ell}| \ind{|\estok_{\ell}| \ge \ThreshTuneParam, |\esttarg_{\ell}| < \ThreshTuneParam} \\
&  \le  \sum_{\ell=1}^D \left(|\esttarg_{\ell}| +  |\estok_{\ell} - \esttarg_{\ell}|  \right) \ind{|\estok_{\ell}| \ge \ThreshTuneParam, |\esttarg_{\ell}| < \ThreshTuneParam} \\
&  \le \ThreshTuneParam \nnzSizeTarg  +  \sum_{\ell=1}^D  |\estok_{\ell} - \esttarg_{\ell}|  \ind{|\estok_{\ell}| \ge \ThreshTuneParam, |\esttarg_{\ell}| < \ThreshTuneParam} 
\end{aligned}
\end{equation}
where the first two claims follow from Definition \ref{def:gen_thresh}.

We control the final term in \eqref{lem:gen_thresh_tail_cond_init_resid_2__1} as follows.
Take $\alpha \in (0, 1)$.
\begin{equation} \label{lem:gen_thresh_tail_cond_init_resid_2__1_final}
\begin{aligned}
\sum_{\ell=1}^D  |\estok_{\ell} - \esttarg_{\ell}|  \ind{|\estok_{\ell}| \ge \ThreshTuneParam, |\esttarg_{\ell}| < \ThreshTuneParam} 
& = \sum_{\ell=1}^D  |\estok_{\ell} - \esttarg_{\ell}| \left[ \ind{|\estok_{\ell}| \ge \ThreshTuneParam, \alpha \ThreshTuneParam < |\esttarg_{\ell}| <  \ThreshTuneParam} +  \ind{|\estok_{\ell}| \ge \ThreshTuneParam, |\esttarg_{\ell}| < \alpha \ThreshTuneParam}  \right]\\
 & \le  R_{\text{max}} \nnzSizeTarg +  \sum_{\ell=1}^D  |\estok_{\ell} - \esttarg_{\ell}|\ind{|\estok_{\ell}| \ge \ThreshTuneParam, |\esttarg_{\ell}| < \alpha \ThreshTuneParam} \\
 & \le  R_{\text{max}}  \left( \nnzSizeTarg + N(1 - \alpha)  \right)
\end{aligned}
\end{equation}
where $N(a) := \sum_{\ell=1}^D \ind{ |\estok_{\ell} - \esttarg_{\ell}|> a \ThreshTuneParam}$.
Note
\begin{equation} \label{eq:gen_thresh_tail_cond_init_resid_N_bound}
\mathcal{P}\left( N(1 - \alpha) > 0 \right) \le  \mathcal{P}\left(  \max_{\ell \in [D]} |\estok_{\ell} - \esttarg_{\ell}|> (1- \alpha) \ThreshTuneParam \right) \le  D f\left(n, (1 - \alpha) \ThreshTuneParam \right)
\end{equation}
by Definition \ref{def:tail_error}, the union bound and the assumption that $\ThreshTuneParam \le   \frac{1}{\tailconst}$.

Combining \eqref{lem:gen_thresh_tail_cond_init_resid_1_and_2}, \eqref{lem:gen_thresh_tail_cond_init_resid_1}, \eqref{lem:gen_thresh_tail_cond_init_resid_2}, \eqref{lem:gen_thresh_tail_cond_init_resid_2__3}, \eqref{lem:gen_thresh_tail_cond_init_resid_2__2}, \eqref{lem:gen_thresh_tail_cond_init_resid_2__1},  \eqref{lem:gen_thresh_tail_cond_init_resid_2__1_final} and \eqref{eq:gen_thresh_tail_cond_init_resid_N_bound} we have
\begin{equation}
\begin{aligned}
\prob{ ||\threshop{\estok}{\ThreshTuneParam} - \esttarg||_1 > t }
& \le \prob{3 \nnzSizeTarg (\ThreshTuneParam  + R_{\text{max}}) + N(1 -\alpha) R_{\text{max}} > t } \\
& \le \prob{R_{\text{max}} > \frac{t}{3 \nnzSizeTarg}  - \ThreshTuneParam \right)  + \mathcal{P} \left(N(1 -\alpha)  > 0} \\
& \le \nnzSizeTarg f\left(n, \frac{t}{3 \nnzSizeTarg}  - \ThreshTuneParam \right) + D f\left(n, (1- \alpha) \ThreshTuneParam \right) 
\end{aligned}
\end{equation}
where the last inequality comes from the union bound.
Plugging in $t = 6 \nnzSizeTarg \ThreshTuneParam$ leaves us with \eqref{eq:gen_thresh_tail_cond_init_resid_L1_bound}.
Note these bounds are valid by the assumption that $\ThreshTuneParam \le   \frac{1}{\tailconst}$.

Repeating this argument by summing over a single row of $\mathcal{A}\left(\threshop{\estok}{\ThreshTuneParam} - \esttarg \right)$ we obtain 
\begin{equation}
\prob{ \opone{ \threshop{\estok}{\ThreshTuneParam} - \esttarg }> t }  \le \prob{3 \maxDegTarg (\ThreshTuneParam  + R_{\text{max}}) + N(1 -\alpha) R_{\text{max}} > t },
\end{equation}
and \eqref{eq:gen_thresh_tail_condition_init_resid_opone_bound} follows.

\end{proof}

\section{Proofs for Section \ref{s:regression_examples}} \label{a:proofs__regression_examples}

\subsection{Preliminary facts for Section \ref{ss:theo_ex_lin_reg}} \label{as:lin_reg__prelim}

We will make use of the following result.
\begin{proposition} \label{prop:lip_perturb}
Let $f, g: \mathbb{R}^d \to \mathbb{R}$.
Suppose $f$ is differentiable and $\mu$ strongly convex.
Also suppose $g$ is $L$ Lipchitz continuous . 
Let $x =  \underset{z \in C}{\text{argmin}} \; f(z)$ and let $y  \in \underset{z \in C}{\text{argmin}} \;   f(z)  + g(z)$ for some set $C$.
Then
$$
||x - y||_2 \le \frac{2L}{\mu} 
$$
\end{proposition}
%
%
%

The following corollary may be useful for settings beyond linear regression.
Often standard concentration bounds can be used to control quantities related to $\estorc$.
If we know the loss function is convex, then we can guarantee $\lassoOrcSet{\residDiffGapRatio, \tuneparam}$ is uniformly close to $\estorc$ as follows.
\begin{corollary} \label{cor:convex_loss_bounds_lasso_orc_dist_to_orc}
Suppose the loss function $\ell(\cdot)$ is $\mu$ strongly convex when restricted to the variables in $\suppSetTarg$.
Then
$$
\sup_{\est \in \lassoOrcSetFrob{\residDiffGapRatio, \tuneparam}} || \est - \estorc ||_2
\le \frac{2^{9/2}a_0  \sqrt{\maxccTarg} \tuneparam \residDiffGapRatio^2 }{\mu} 
$$
and
$$
\sup_{\est \in \lassoOrcSetFrob{\residDiffGapRatio, \tuneparam}} \opn{\mathcal{L}(|\est| - |\esttarg|)} \le
\opn{\mathcal{L}(|\esttarg| - |\estorc|)} + \frac{32 a_0 \maxccTarg \tuneparam \residDiffGapRatio^2 }{\mu}.
$$
Furthermore
$$
\sup_{\est \in \lassoOrcSetOp{\residDiffGapRatio, \tuneparam}} || \est - \estorc ||_2
\le \frac{2^{5}a_0  \sqrt{\suppSizeTarg} \tuneparam \residDiffGapRatio^2 }{\mu} 
$$
and a similar claim about $\opn{\mathcal{L}(|\esttarg| - |\estorc|)}$ holds.

\end{corollary}

The next proposition controls the error for the linear regression oracle.
\begin{proposition} \label{prop:lin_reg_orc_lap_resid_bound_via_Linfty}
For the linear regression setting of Section  \ref{ss:theo_ex_lin_reg}
$$
\prob{ ||\estorc_{\suppSetTarg}  - \esttarg_{\suppSetTarg}||_{\text{max}}  \ge t}
\le 2 \suppSizeTarg \exp \left(  \frac{-n\minevalxsshort  t^2}{2 \sigma^2} \right)
$$
thus
$$
\prob{\opn{\mathcal{L}(|\estorc| - |\esttarg|)} \ge t}
\le 2 \suppSizeTarg \exp \left(  \frac{-n\minevalxsshort t^2}{8 \maxccTargPow{2} \sigma^2} \right).
$$
\end{proposition}

\subsection{Proofs for Section \ref{ss:theo_ex_lin_reg}} \label{as:lin_reg__prelim}

\begin{proof} of Theorem \ref{thm:thoe_ex__lin_reg__reg_sat_whp}
To verify the two step convergence claims we will appeal to Corollary \ref{cor:suff_cond_whp_for_init_from_close_give_lasso_orc}.
When restricted to the variables in  $\suppSetTarg$ the least squares loss $\ell_{\suppSetTarg}(\est_{\suppSetTarg})  = \frac{1}{2} ||X_{\suppSetTarg} \est_{\suppSetTarg} - y||_2^2$ is $\minevalxsshort$ strongly convex.
Therefore below we will first establish concentration inequalities for $\estorc$ then extend them to $\lassoOrcSet{\residDiffGapRatio, \tuneparam}$ using Corollary \ref{cor:convex_loss_bounds_lasso_orc_dist_to_orc}.

\ProofSection{Tail bounds for gradient}

We closely follow the proof of Theorem 3 of \cite{fan2014strong} to obtain tail bounds on the gradient at the oracle.
Note  $\nabla \ell(\est) = \frac{1}{n}X^T (X \est - y)$ and 
\begin{equation} \label{eq:lin_reg_org_equals_targ_plus_noise}
\estorc = \esttarg + (X_{\suppSetTarg}^T X_{\suppSetTarg})^{-1} X_{\suppSetTarg}^T \epsilon
\end{equation} 
 thus 
$$\nabla_j \ell(\estorc)  = \frac{1}{n}X_j^T (H_{\supp} - I_n) \epsilon$$
 where $H_{\suppSetTarg} := X_{\suppSetTarg} (X_{\suppSetTarg}^T X_{\suppSetTarg})^{-1} X_{\suppSetTarg}$ is a projection matrix.
Using this fact we see $|| \frac{1}{n}X_j^T (H_{\suppSetTarg} - I_n)||_2^2 \le \frac{1}{n} ||X_j||_2^2$, thus Hoeffding's bound for sub-Gaussian random variables we have 
\begin{equation} \label{eq:lin_reg_orc_grad_concentration} 
\prob{ |\nabla_j \ell(\estorc) | \ge t}
\le 2 \exp \left( -\frac{n^2 t^2}{2 \sigma^2 ||X_j||^2}\right)
\le 2 \exp \left( -\frac{n t^2}{2 \sigma^2 \maxColNorm}\right).
\end{equation} 

Since  $\nabla \ell(\est) = \frac{1}{n}X^T (X \est - y)$, we see
\begin{equation*}
\begin{aligned}
|\nabla_j \ell(\est) - \nabla_j \ell(\widetilde{\est})|  & = \frac{1}{n} | X_j^T X (\est - \widetilde{\est})|\\
& \le  \frac{1}{n} ||X^T X_j||_2  ||\est - \widetilde{\est}||_2 \\
& \le \frac{||X_j||_2}{\sqrt{n}} \sqrt{\maxevalxexpl}||\est - \widetilde{\est}||_2\\
& \le  \sqrt{\maxColNorm \maxevalxshort } ||\est - \widetilde{\est}||_2
\end{aligned}
\end{equation*}
where the first inequality is from Cauchy-Schwartz and the second inequality is from the Courrant-Fischer theorem.
In other words $\nabla_j \ell(\est)$ is $\sqrt{\maxColNorm \maxevalxshort }$ Lipschitz.
Therefore for any $j$,
\begin{equation*}
\begin{aligned}
\sup_{\est \in \lassoOrcSet{\residDiffGapRatio, \tuneparam}} |\nabla_j \ell(\est) |
 \le |\nabla_j \ell(\estorc) | +\sqrt{\maxColNorm \maxevalxshort }  \sup_{\est \in \lassoOrcSet{\residDiffGapRatio, \tuneparam}} || \est - \estorc||_2,
\end{aligned}
\end{equation*}
thus by the first claim of Corollary \ref{cor:convex_loss_bounds_lasso_orc_dist_to_orc},
\begin{equation*}
\begin{aligned}
\sup_{\est \in \lassoOrcSet{\residDiffGapRatio, \tuneparam}} ||\nabla_{\suppSetTargPow{C} } \ell(\est)  ||_{\text{max}} 
 \le ||\nabla_{\suppSetTargPow{C} }  \ell(\estorc) ||_{\text{max}} +\sqrt{\maxColNorm \maxevalxshort } \frac{2^{9/2}a_0  \sqrt{\maxccTarg} \tuneparam \residDiffGapRatio^2 }{\minevalxsshort}.
\end{aligned}
\end{equation*}
Putting this together with \eqref{eq:lin_reg_orc_grad_concentration} and the union bound,
\begin{equation*}
\begin{aligned}
\prob{\sup_{\est \in \lassoOrcSet{\residDiffGapRatio, \tuneparam}} ||\nabla_{\suppSetTargPow{C} } \ell(\est)  ||_{\text{max}}  \ge t}
& \le \prob{ ||\nabla_{\suppSetTargPow{C} } \ell(\estorc)  ||_{\text{max}} 
\ge t - \sqrt{\maxColNorm \maxevalxshort } \frac{2^{9/2}a_0  \sqrt{\maxccTarg} \tuneparam \residDiffGapRatio^2 }{\minevalxsshort}
} \\
& \le 2 (D - \suppSizeTarg) \exp\left(\frac{-n}{2 \sigma^2 \maxColNorm} \left[
 t - \sqrt{\maxColNorm \maxevalxshort } \frac{2^{9/2}a_0  \sqrt{\maxccTarg} \tuneparam \residDiffGapRatio^2 }{\minevalxsshort}
\right]_+^2 \right)
\end{aligned}
\end{equation*}
Setting in $t = \frac{a_1\tuneparam}{2}$ and using Assumption \eqref{eq:theo_ex_lin_reg_init_tol} we obtain \eqref{eq:thoe_ex__lin_reg__reg_sat_whp__prob_nice_grad}. 


\ProofSection{Tail bounds for Laplacian operator norm}

By the second claim of Corollary \ref{cor:convex_loss_bounds_lasso_orc_dist_to_orc},
\begin{equation*}
\begin{aligned}
\prob{\sup_{\est \in \lassoOrcSet{\residDiffGapRatio, \tuneparam}} \opn{\mathcal{L}(|\est| - |\esttarg|)}  \ge t}
& \le \prob {\opn{\mathcal{L}(|\esttarg| - |\estorc|)}  \ge t -  \frac{32 a_0 \maxccTarg \tuneparam \residDiffGapRatio^2 }{\minevalxsshort} } \\
& \le  2 \suppSizeTarg \exp \left(  \frac{-n\minevalxsshort }{8 \maxccTargPow{2} \sigma^2}
\left[
 t -  \frac{32 a_0 \maxccTarg \tuneparam \residDiffGapRatio^2 }{\minevalxsshort}
\right]_+^2 \right).
\end{aligned}
\end{equation*}
where the second inequality comes from Proposition \ref{prop:lin_reg_orc_lap_resid_bound_via_Linfty}.
Setting $t = b_1 \tuneparam$ and using Assumption \eqref{eq:theo_ex_lin_reg_init_tol} we obtain \eqref{eq:thoe_ex__lin_reg__reg_sat_whp__prob_small_resid}. 


\end{proof}

\begin{proof} of Proposition \ref{prop:lin_reg_orc_lap_resid_bound_via_Linfty}
Note
$$
r := \estorc_{\suppSetTarg}  - \esttarg_{\suppSetTarg} = (X_{\suppSetTarg}^T X_{\suppSetTarg})^{-1}X_{\suppSetTarg}^T \epsilon
$$
and
$$
\opn{ (X_{\suppSetTarg}^T X_{\suppSetTarg})^{-1}X_{\suppSetTarg}^T} 
= \sqrt{ \opn{(X_{\suppSetTarg}^T X_{\suppSetTarg})^{-1}X_{\suppSetTarg}^T X_{\suppSetTarg} (X_{\suppSetTarg}^T X_{\suppSetTarg})^{-1} } }
=  \sqrt{ \opn{ (X_{\suppSetTarg}^T X_{\suppSetTarg})^{-1} }}
=  \sqrt{\frac{1}{n \minevalxsshort}}.
$$

Note $r_j = (X_{\suppSetTarg} (X_{\suppSetTarg}^T X_{\suppSetTarg})^{-1} e_j)^T \epsilon$ so by Hoeffding's inequality $r_j$ is sub-Gaussian with parameter $||X_{\suppSetTarg} (X_{\suppSetTarg}^T X_{\suppSetTarg})^{-1} e_j||_2 \sigma \le \opn{X_{\suppSetTarg} (X_{\suppSetTarg}^T X_{\suppSetTarg})^{-1}} \sigma =  \sqrt{\frac{\sigma^2}{n \minevalxsshort}}$.
By the union bound we obtain
$$
\prob{||r||_{\text{max}} \ge t} \le 2 \suppSizeTarg \exp \left(  \frac{-n\minevalxsshort t^2}{2 \sigma^2} \right)
$$
and the first claim follows.
The second claim follows from Proposition \ref{prop:lap_bounds}.
\end{proof}

\begin{proof} of Proposition \ref{prop:lip_perturb}
By strong convexity $f(y)  \ge f(x) + \nabla f(x)^T (y - x) + \frac{\mu}{2}||y - x||_2^2 = f(x) + \frac{\mu}{2}||y - x||_2^2 $.
By  Lipchitz continuity $g(y) \ge g(x) - L||y - x||_2$.
Therefore
\begin{equation}
f(y) + g(y) \ge f(x) + g(x) + \frac{\mu}{2}||x - y||^2 - L||x - y||_2.
\end{equation}
Since $y$ is a minimizer of $f + g$, $f(y) + g(y) - (f(x) + g(x)) \le 0$ so we get $\frac{\mu}{2}||x - y||^2 - L||x - y||_2 \le 0$, thus
$$
||x - y||_2 \le \frac{2L}{\mu}
$$
and the result follows.

\end{proof}

\begin{proof} of Corollary \ref{cor:convex_loss_bounds_lasso_orc_dist_to_orc}
Let $\est \in \lassoOrcSet{\residDiffGapRatio, \tuneparam}$ be a solution to the weighted Lasso oracle problem \eqref{prob:weighted_lasso_oracle}.
The weighted Lasso penalty function is given by $g(\est) := \frac{1}{2} w^T |\est|$ where $w$ satisfies the conditions of Definition \eqref{eq:def_lassoOrcSetFrob}.
Applying Cauchy-Schwartz
\begin{equation*}
\begin{aligned}
|g(x) - g(y)| 
& = \frac{1}{2}  w^T | |x| - |y| | \\
& \le \frac{1}{2}  w^T  |x - y| \\
& \le \frac{1}{2}  ||w||_2 ||x - y||_2 \\
& \le \frac{1}{2}  2^{9/2}a_0  \sqrt{\maxccTarg} \tuneparam \residDiffGapRatio^2 ||x - y||_2 ,
\end{aligned}
\end{equation*}
we conclude that $g(\cdot)$ has a Lipschitz constant given by $2^{7/2}a_0  \sqrt{\maxccTarg} \tuneparam \residDiffGapRatio^2 $.
The first claim then follows by Proposition \ref{prop:lip_perturb}. 

Recall the supports of $\est, \estorc$ are subsets of the block support of $\esttarg$.
Therefore the largest number of nodes in a connected component of $\mathcal{A}(|\est| - |\estorc| )$  is given by $\maxccTarg$ (which bounds the largest degree).
Applying the first claim,
 the fact that $|| |\est| - |\estorc| ||_2 \le || \est - \estorc ||_2$
 and Proposition \ref{prop:lap_bounds} we see
$$
\sup_{\est \in \lassoOrcSet{\residDiffGapRatio, \tuneparam} } \opn{\mathcal{L}(|\est| - |\estorc| )} \le \sqrt{2 \maxccTarg} \cdot \frac{2^{9/2}a_0  \sqrt{\maxccTarg} \tuneparam \residDiffGapRatio^2 }{\mu}
$$
and the second claim follows.

The operator norm claims follow similarly after noting for any $w$ satisfying \eqref{eq:def_lassoOrcSetOp} we have
\begin{equation}
\begin{aligned}
|g(x) - g(y)| 
%
%
& \le \frac{1}{2}  w^T  |x - y| \\
& \le \frac{1}{2}  ||w||_{\text{max}} ||x - y||_1 \\
& \le \frac{1}{2}  2^{5}a_0 \tuneparam \residDiffGapRatio^2 \sqrt{\suppSizeTarg} ||x - y||_2 ,
\end{aligned}
\end{equation}
where the second inequality comes from Holder's and the final inequality uses the fact that $x - y$ is supported on $\suppSizeTarg$ entries.

\end{proof}

\subsection{Preliminary facts  for Section \ref{ss:theo_ex_logisitc}} \label{as:logistic_reg__prelim}

To obtain Corollary \ref{cor:theo_ex_logisitc_lasso_init} we need a bound on the $L_1$ error for logistic regression with a Lasso penalty,
\begin{equation} \label{eq:logisic_regression_lasso_pen}
\begin{aligned}
& \underset{\est \in \mathbb{R}^{{D}}}{\textup{minimize}}  & & \frac{1}{n} \sum_{i=1}^n  \left(-y_i X(i, :)^T \est + \LogisitcLink{X(i, :)^T \est}  \right)+\ThreshTuneParam ||\est||_1.
\end{aligned}
\end{equation}
The following theorem is a consequence of the proof of Theorem 5 of \citep{fan2014strong}, which gives the analogous result for the $L_2$ norm.
We work under the following restricted eigenvalue condition.
\begin{equation} \label{eq:logistic_restr_eval}
\begin{aligned}
\LogisticRestrEval := \; & \underset{u \neq 0}{\textup{minimize}}  & & \frac{u^T \nabla^2 \ell(\esttarg) u }{||u||_2} && \in (0, \infty)\\ 
& \text{subject to } & & ||u_{\nnzSetTargPow{C}}||_1 \le 3 || u_{\nnzSetTarg}||_1.  &&
\end{aligned}
\end{equation}

\begin{theorem} \label{thm:logisitc_regession_lasso_pen_L1_err}
Let $\estLasso$ be a solution to Problem \eqref{eq:logisic_regression_lasso_pen}.
Under assumptions \eqref{eq:logistic_setup} and \eqref{eq:logistic_restr_eval} if 
$$\ThreshTuneParam  \le \frac{\LogisticRestrEval}{20 \LogisticMaxX \nnzSizeTarg }$$
then with probability at least $1 - 2D  \exp\left(-\frac{n}{2 \LogisticMaxVarNorm}\ThreshTuneParam^2 \right)$ we have
$$
||\estLasso - \esttarg||_1 \le \frac{20 \nnzSizeTarg}{\LogisticRestrEval}\ThreshTuneParam
$$
\end{theorem}

Next make some definitions for the Lasso penalized logistic regression problem restricted to the block support set, $\suppSetTarg$ (i.e. under the constraint $\est_{\suppSetTargPow{C}} = 0$).
We will assume with out loss of generality that $\suppSizeTarg = [\suppSetTarg]$ i.e. $\esttarg_j = 0$ if $j > \suppSizeTarg$.
Let $w \in \mathbb{R}^{\supptarg}_+$ be the weight vector satisfying the constraints of Definition \eqref{eq:def_lassoOrcSetOp} for the Lasso oracle problem,
\begin{equation}  \label{prob:lasso_orc_prob}
\begin{aligned}
& \underset{\est \in \mathbb{R}^{\suppSizeTarg}}{\textup{minimize}}  & & \ell_{\suppSetTarg}(\est)  + \sum_{j=1}^{\suppSizeTarg} w_j | \est_j| 
\end{aligned}
\end{equation}
where $\ell_{\suppSetTarg}: \mathbb{R}^{\suppSizeTarg} \to \mathbb{R}$, is given by $\ell_{\suppSetTarg}(\est) = \frac{1}{n} \sum_{i=1}^n  -y_i X_{\suppSetTarg}(i, :)^T \est + \LogisitcLink{X_{\suppSetTarg}(i, :)^T \est}$.
Note 
$$\nabla\ell_{\suppSetTarg}(\est) =  \frac{1}{n}X^T_{\suppSetTarg} \left( \LogisitcMuS{\est}  - y  \right) ,$$
where $\LogisitcMuS{\cdot}: \mathbb{R}^{\suppSizeTarg} \to \mathbb{R}^n$ is given by
$$
\LogisitcMuS{\est}_i =  \LogisitcLinkGrad{X_{\suppSetTarg}(i, :)^T \est}, i = 1, \dots, n.
$$
The Jacobian of $\LogisitcMuS{\cdot}$ is 
$$
D \LogisitcMuS{\est} = \text{diag}(\LogisitcHS{\est}) X_{\suppSetTarg}
$$
where $\LogisitcHS{\cdot}: \mathbb{R}^{\suppSizeTarg} \to \mathbb{R}^n$ is given by
$$
\LogisitcHS{\est}_i = \LogisitcLinkGG{X_{\suppSetTarg}(i, :)^T \est}, i = 1, \dots, n,
$$
and we see
$$
\nabla^2 \ell_{\suppSetTarg}(\est) = \frac{1}{n} X_{\suppSetTarg}^T \text{diag}(\LogisitcHS{\est}) X_{\suppSetTarg}.
$$
The Jacobian of $\LogisitcHS{\cdot}$ is 
$$
D \LogisitcHS{\est} = \text{diag}(\LogisitcTS{\est}) X_{\suppSetTarg}
$$
where $\LogisitcTS{\cdot}: \mathbb{R}^{\suppSizeTarg} \to \mathbb{R}^n$ is given by
$$
\LogisitcTS{\est}_i = \LogisitcLinkGGG{X_{\suppSetTarg}(i, :)^T \est}, i = 1, \dots, n.
$$.

We also need the analogs, $\LogisitcMu{\cdot}, \LogisitcH{\cdot}, \LogisitcT{\cdot}$, for the full problem e.g. where $\LogisitcH{\cdot}: \mathbb{R}^D \to \mathbb{R}^n$ is given by $\LogisitcH{\est}_i = \LogisitcLinkGG{X(i, :)^T \est}$ for each $i \in [n]$.

\subsection{Proofs for Section \ref{ss:theo_ex_logisitc}} \label{as:logistic_reg__proofs}

\begin{proof} of Theorem \ref{thm:thoe_ex__logistic__reg_sat_whp}

Our proof is similar to the proof of Theorem 4 in \citep{fan2014strong}.

\ProofSection{1. Reduce the theorem to statements about the residual}

Throughout the proof we let $\est \in \PrevlassoOrcSet{\ParmMax}$ be any element.
We reduce the theorem's claims to statements about the residual $\est_{\suppSetTarg} - \esttarg_{\suppSetTarg}$.
We first show sufficient conditions for Claims \eqref{eq:theo_ex_logistic__reg_sat__prob_small_resid} and  \eqref{eq:theo_ex_logistic__reg_sat__prob_nice_grad}.

By Proposition \ref{prop:lap_bounds}, we have $\opn{\mathcal{L}(|\est| - \esttarg|)} \le 2 \maxccTarg ||\est - \esttarg||_{\text{max}}$.
Therefore Claim \eqref{eq:theo_ex_logistic__reg_sat__prob_small_resid} will hold if we can ensure  
\begin{equation}\label{eq:logistic_resid_bound_nts_gap}
||\est- \esttarg||_{\text{max}} = ||\est_{\suppSetTarg} - \esttarg_{\suppSetTarg}||_{\text{max}} \le  \frac{b_1 \tuneparam}{2\maxccTarg}
\end{equation}
 with the stated probability. 

Fix $j \in \suppSetTargPow{C}$ and let $g_j: \mathbb{R}^D \to \mathbb{R}$ be given by $g_j(b) := \nabla_{j} \ell(b) =  \frac{1}{n} X_j^T (\LogisitcMu{b} - y)$.
By a Taylor expansion for any $b \in \mathbb{R}^D$ with $b_{\suppSetTargPow{C}} = 0$,
\begin{equation}  \label{eq:soft_chair}
\begin{aligned}
|g_j(b) - g_j(\esttarg)| & \le |\nabla g_j(\esttarg)^T (b - \esttarg)| +  \frac{1}{2} |(b - \esttarg)^T \nabla^2 g_j(\widetilde{b}) (b - \esttarg)| \\
&  = |\nabla_{\suppSetTarg} g_j(\esttarg)^T (b_{\suppSetTarg} - \esttarg_{\suppSetTarg})| +  \frac{1}{2} |(b_{\suppSetTarg} - \esttarg_{\suppSetTarg})^T \nabla_{\suppSetTarg}^2 g_j(\widetilde{b}) (b_{\suppSetTarg} - \esttarg_{\suppSetTarg})|
\end{aligned} 
\end{equation}
for some $\widetilde{b}$ on the line segment between $b$ and $\esttarg$ where $\nabla_{\suppSetTarg}^2$ means we restrict the Hessian to the sub-matrix whose rows/columns are index by $\suppSetTarg$.
Note 
$$ \nabla_{\suppSetTarg} g_j(\esttarg)=  \frac{1}{n} X_{j}^T \text{diag}(\LogisitcH{\esttarg}) X_{\suppSetTarg},
\qquad
\nabla^2_{\suppSetTarg} g_j(\widetilde{b})= \frac{1}{n} X_{\suppSetTarg} \text{diag}(\LogisitcT{\widetilde{b}} \odot  X_j ) X_{\suppSetTarg}.
$$
For the first term in the Taylor expansion we use Holder's inequality and the definition of $\LogisticXsHXscOne$ to obtain
\begin{equation} \label{eq:logistic_taylor_grad_control}
|\nabla_{\suppSetTarg} g_j(\esttarg)^T (\est_{\suppSetTarg} - \esttarg_{\suppSetTarg})|
\le ||\nabla_{\suppSetTarg} g_j(\esttarg)||_1 ||\est_{\suppSetTarg} - \esttarg_{\suppSetTarg}||_{\text{max}}  
\le  \suppSizeTarg \LogisticXsHXscOne ||\est_{\suppSetTarg} - \esttarg_{\suppSetTarg}||_{\text{max}}.
\end{equation}

Next we bound the second term.
We can check $|\LogisitcLinkGGG{t}| \le \frac{1}{4}$ for all $t$ as in the proof of Theorem 4 of \citep{fan2014strong}.
Using this fact and the form of $\nabla_{\suppSetTarg}^2 g_j(\widetilde{b}) $ we can verify that for any $y \in \mathbb{R}^{\supptarg}$
\begin{equation} \label{eq:snowy_tree}
| y^T \nabla^2_{\suppSetTarg} g_j(\widetilde{b})  y|  \le \frac{1}{4} y^T \tempMatFrank_j y \le  \frac{\lambda_{\text{max}}( \tempMatFrank_j) }{4} ||y||_2^2 \le \frac{\LogisticEvalMax}{4} ||y||_2^2,
\end{equation}
where $\tempMatFrank_j := \frac{1}{n} X_{\suppSetTarg}^T \text{diag}(|X_j|) X_{\suppSetTarg} $.
Note the second inequality uses the fact that $\tempMatFrank_j$ is positive semi-definite.

Putting this together with  \eqref{eq:soft_chair} and \eqref{eq:logistic_taylor_grad_control} we obtain
\begin{equation} \label{eq:logistic_nice_grad_det_cond}
\begin{aligned}
|| \nabla_{\suppSetTargPow{C} } \ell(\est) ||_{\text{max}} 
&  \le \max_{j \in \suppSetTargPow{C} } \left|\frac{1}{n} X_j^T (\LogisitcMu{\esttarg} - y) \right| +  \suppSizeTarg \LogisticXsHXscOne ||\est_{\suppSetTarg} - \esttarg_{\suppSetTarg}||_{\text{max}} + \frac{\LogisticEvalMax}{8} ||\est_{\suppSetTarg} - \esttarg_{\suppSetTarg}||_2^2  \\
& \le \max_{j \in \suppSetTargPow{C} } \left|\frac{1}{n} X_j^T (\LogisitcMu{\esttarg} - y) \right| +  \left(\suppSizeTarg \LogisticXsHXscOne + \frac{\sqrt{\suppSizeTarg} \LogisticEvalMax}{8} \right)  ||\est_{\suppSetTarg} - \esttarg_{\suppSetTarg}||_{\text{max}}.
\end{aligned}
\end{equation}

Claim \eqref{eq:theo_ex_logistic__reg_sat__prob_nice_grad}  therefore holds if we can ensure both
\begin{equation} \label{eq:logistic_resid_bound_nts_grad}
||\est_{\suppSetTarg} - \esttarg_{\suppSetTarg}||_{\text{max}} \le \frac{a_1 \tuneparam }{8\maxccTarg \left(\suppSizeTarg \LogisticXsHXscOne + \frac{\sqrt{\suppSizeTarg} \LogisticEvalMax}{8} \right)}
\end{equation}
and
\begin{equation}  \label{eq:logistic_grad_at_targ_nts_grad}
\max_{j \in \suppSetTargPow{C} } \left|\frac{1}{n} X_j^T (\LogisitcMu{\esttarg} - y) \right| \le \frac{a_1 \tuneparam }{8 \maxccTarg}
\end{equation}
holds with high probability. 
We can ensure \eqref{eq:logistic_grad_at_targ_nts_grad}  holds with high probability by the union bound and Proposition 4a of \citep{fan2011nonconcave},
\begin{equation} \label{eq:logistic_grad_at_targ_grad_holds_whp}
\begin{aligned}
\prob{\max_{j \in \suppSetTargPow{C} } \left|\frac{1}{n} X_j^T (\LogisitcMu{\esttarg} - y) \right| \ge \frac{a_1 \tuneparam}{8 \maxccTarg} } 
\le
(D - \suppSizeTarg) \exp\left(\frac{-na_1^2 \tuneparam^2}{64 \maxccTargPow{2} \LogisticMaxVarNorm} \right).
\end{aligned}
\end{equation}
Thus the remainder of the proof is devoted to \eqref{eq:logistic_resid_bound_nts_gap} and \eqref{eq:logistic_resid_bound_nts_grad}.
Note since $\est \in \PrevlassoOrcSet{\ParmMax}$ is arbitrary, we will have actually proved bounds on $\sup_{\est \in \PrevlassoOrcSet{\ParmMax}} ||\est- \esttarg||_{\text{max}} $.

\ProofSection{2. Use a fixed point to argument control the residual}

The stationary points of  \eqref{prob:lasso_orc_prob} are characterized by the first order necessary conditions 
 $\LogisticGradCond{\est} = 0$ where 
\begin{equation}
\LogisticGradCond{\est}  = \nabla \ell_{\suppSetTarg}(\est) + w \odot \zeta, 
\end{equation}
and $\zeta \in \mathbb{R}^{\supptarg}$ is any sub-gradient of the weighted lasso penalty evaluated at $\est $ i.e.
\begin{equation} \label{eq:logistic_lasso_orc_subgrad}
\zeta_j = 
\begin{cases}
\text{sign}(\est_j) & \text{if } \est \neq 0 \\
\in [-1, 1] & \text{if } \est_j = 0.
\end{cases}
\end{equation}
Note  \eqref{prob:lasso_orc_prob}  has a unique solution since $h(\cdot) > 0$ and the columns of $X_{\suppSetTarg}$ are linearly independent.

Let $F: \mathbb{R}^{\supptarg} \to \mathbb{R}^{\supptarg}$ be the following map
\begin{equation}
\begin{aligned}
F(\triangle_{\suppSetTarg}) 
& = - \nabla^2 \ell_{\suppSetTarg}(\esttarg_{\suppSetTarg}) ^{-1} \LogisticGradCond{\esttarg_{\suppSetTarg} +\triangle_{\suppSetTarg} } + \triangle_{\suppSetTarg} \\
& = \left(\frac{1}{n} X_{\suppSetTarg}^T \text{diag}(\LogisitcHS{\esttarg_{\suppSetTarg}}) X_{\suppSetTarg} \right)^{-1} \left(\frac{1}{n}X^T_{\suppSetTarg} (y - \LogisitcMuS{\esttarg_{\suppSetTarg} +\triangle_{\suppSetTarg}}) -  w \odot \zeta \right) + \triangle_{\suppSetTarg} .
\end{aligned}
\end{equation}
Note $F$ is continuous since since $h(\cdot) > 0$ and the columns of $X_{\suppSetTarg}$ are linearly independent. 

We can check $F(\triangle_{\suppSetTarg}) = \triangle_{\suppSetTarg}$ if and only if $\LogisticGradCond{\esttarg_{\suppSetTarg} +\triangle_{\suppSetTarg} } = 0$. 
In other words, $\triangle$ is a fixed point of $F$ if and only if $\esttarg_{\suppSetTarg} +\triangle_{\suppSetTarg} $ is the unique solution to \eqref{prob:lasso_orc_prob}.
Let $B(r) := \{b \in  \mathbb{R}^{\supptarg} \text{ s.t. }  ||b||_{\text{max}} \le r \} $. 
Suppose for some $r$ we can show  $F$ maps $B(r)$ onto itself. 
Since $F$ is continuous and $B(r)$ is convex and compact, by Brouwer's  fixed point theorem (e.g. p161 of  \citealt{ortega2000iterative}) $F$ must have a fixed point inside $B(r)$.
Putting this all together we see 
\begin{equation} \label{eq:contraction_to_implication}
F(B(r)) \subseteq B(r) \implies ||\est_{\suppSetTarg} - \esttarg_{\suppSetTarg}||_{\text{max}} \le r,
\end{equation}
for any solution $\est_{\suppSetTarg}$ of \eqref{prob:lasso_orc_prob}.

\ProofSection{3. Find conditions under which the desired contraction holds}

From the discussion above we need to verify \eqref{eq:contraction_to_implication} with $r = \LogisticTempResidConst \tau$ where $\LogisticTempResidConst$ is given in the right hand side of either \eqref {eq:logistic_resid_bound_nts_gap} or \eqref{eq:logistic_resid_bound_nts_grad}.
Recall the weight vector used in Problem  \eqref{prob:lasso_orc_prob} satisfies the conditions of $\lassoOrcSetFrob{\residDiffGapRatio, \tuneparam}$ which gives us that 
$$|| \zeta \odot w||_{\text{max}} \le  2^5 a_0 \tuneparam \residDiffGapRatio^2$$
where $\zeta$ are the sub-gradients from \eqref{eq:logistic_lasso_orc_subgrad}.
Note $\residDiffGapRatio$ is selected in Assumption  \eqref{eq:theo_ex_logisitc_init_tol} to ensure 
\begin{equation}\label{eq:logistic_lasso_weight_condition}
|| \zeta \odot w||_{\text{max}}  \le \frac{\LogisticTempResidConst}{4 \LogisticXshXsInv} \tuneparam
\end{equation}
for either of these choices of $\LogisticTempResidConst$.
We will further assume the condition
\begin{equation} \label{eq:logistic_resid_condition}
\begin{aligned}
\frac{1}{n} ||  X^T_{\suppSetTarg}(\LogisitcMuS{ \esttarg_{\suppSetTarg} } - y)  ||_{\text{max}}
\le
  \frac{\LogisticTempResidConst}{4 \LogisticXshXsInv} \tuneparam,
\end{aligned}
\end{equation}
and verify this condition holds with high probability later on.

By bringing the $ \triangle_{\suppSetTarg}$ term inside and adding zero we see
\begin{equation}
\begin{aligned}
F(\triangle_{\suppSetTarg}) 
 =   \left(\frac{1}{n}X_{\suppSetTarg}^T \text{diag}(\LogisitcHS{\esttarg_{\suppSetTarg}}) X_{\suppSetTarg} \right)^{-1} 
\Big[ &
\frac{1}{n} X^T_{\suppSetTarg} (y -  \LogisitcMuS{ \esttarg_{\suppSetTarg} })
- w \odot \zeta \\
& + \frac{1}{n} X^T_{\suppSetTarg}  \left\{ \LogisitcMuS{\esttarg_{\suppSetTarg}} +  \text{diag}(\LogisitcHS{\esttarg_{\suppSetTarg}}) X_{\suppSetTarg}   \triangle_{\suppSetTarg}   -\LogisitcMuS{\esttarg_{\suppSetTarg} +\triangle_{\suppSetTarg} }  \right\}
 \Big].
 \end{aligned}
\end{equation}

By the definition of the $(\infty, \infty)$ operator norm and the triangle inequality we obtain
\begin{equation} \label{eq:meow_nyan_woof_wanwan}
\begin{aligned}
|| F(\triangle_{\suppSetTarg})||_{\text{max}}  &   \le
 \opnorm{ \left( \frac{1}{n} X_{\suppSetTarg}^T \text{diag}(\LogisitcHS{\esttarg}) X_{\suppSetTarg} \right)^{-1}}{\infty}{\infty}
  \Big[ 
  ||  \frac{1}{n} X^T_{\suppSetTarg} (\LogisitcMuS{ \esttarg_{\suppSetTarg}} - y)  ||_{\text{max}}  
 +  || w \odot \zeta ||_{\text{max}}  \\
& \qquad +  \frac{1}{n} ||  X^T_{\suppSetTarg}  \left(  \LogisitcMuS{\esttarg_{\suppSetTarg}} +  \text{diag}(\LogisitcHS{\esttarg_{\suppSetTarg}}) X_{\suppSetTarg}   \triangle_{\suppSetTarg}   - \LogisitcMuS{ \esttarg_{\suppSetTarg} +\triangle_{\suppSetTarg} }  \right)  ||_{\text{max}} 
 \Big].
 \end{aligned}
\end{equation}

Let $v :=  \frac{1}{n} \left[ X^T_{\suppSetTarg}  \left( \LogisitcMuS{\esttarg_{\suppSetTarg}} +  \text{diag}(\LogisitcHS{\esttarg_{\suppSetTarg}}) X_{\suppSetTarg}   \triangle_{\suppSetTarg}   - \LogisitcMuS{ \esttarg_{\suppSetTarg} +\triangle_{\suppSetTarg} }  \right) \right] $ be the third term  of \eqref{eq:meow_nyan_woof_wanwan}.
We use a Taylor series argument to bound $||v||_{\text{max}}$.
For each $j \in [\supptarg]$ let $g_j: \mathbb{R}^{\supptarg} \to \mathbb{R}$, $g_j(\est) = \frac{1}{n} X_j^T \LogisitcMuS{\est}$.
By the Taylor remainder theorem we have
\begin{equation} \label{eq:logistic_F_third_term_j}
v_j
=  |g_j( \esttarg_{\suppSetTarg} +\triangle_{\suppSetTarg})  - g_j(\esttarg_{\suppSetTarg} ) - \nabla g_j(\esttarg_{\suppSetTarg} )^T  \triangle_{\suppSetTarg} | \le \frac{1}{2} | \triangle_{\suppSetTarg}^T \nabla^2 g_j(\widetilde{\est})  \triangle_{\suppSetTarg}|
%
\end{equation}
for some $\widetilde{b}$ on the line segment between $ \esttarg$  and $ \esttarg +\triangle$.
Applying the argument used to obtain  \eqref{eq:snowy_tree} we have
$$
| \triangle_{\suppSetTarg}^T \nabla^2 g_j(\widetilde{b})  \triangle_{\suppSetTarg}| \le \frac{1}{4} \LogisticEvalMax  ||\triangle_{\suppSetTarg}||_2^2 \le \frac{\supptarg}{4} \LogisticEvalMax  ||\triangle_{\suppSetTarg}||_{\text{max}}^2.
$$
Putting this together with \eqref{eq:logistic_F_third_term_j} we obtain
\begin{equation}\label{eq:logistic_F_third_term_bound}
||v||_{\text{max}} \le \frac{\supptarg}{8} \LogisticEvalMax  ||\triangle_{\suppSetTarg}||_{\text{max}}^2.
\end{equation}

Finally, for any $\triangle_{\suppSetTarg} \in B(\LogisticTempResidConst \tuneparam)$
\begin{equation}
\begin{aligned}
|| F(\triangle_{\suppSetTarg})||_{\text{max}}  & \le 
\LogisticXshXsInv
\left(
 \frac{1}{n} ||  X^T_{\suppSetTarg}(  \LogisitcMuS{ \esttarg_{\suppSetTarg}} -y) ||_{\text{max}}  
 + 
|| \zeta \odot w||_{\text{max}}
 + 
  \frac{\supptarg}{8} \LogisticEvalMax  ||\triangle_{\suppSetTarg}||_{\text{max}}^2
 \right)\\
 %
& \le 
\LogisticXshXsInv
\left(
\frac{\LogisticTempResidConst}{4\LogisticXshXsInv} \tuneparam
 + 
\frac{\LogisticTempResidConst}{4\LogisticXshXsInv} \tuneparam
 + 
 \frac{\suppSizeTarg}{8} \LogisticEvalMax  \LogisticTempResidConst^2 \tuneparam^2
 \right)\\
 & \le \LogisticTempResidConst \tuneparam.
 \end{aligned}
\end{equation}
The first inequality comes from  \eqref{eq:meow_nyan_woof_wanwan},  \eqref{eq:logistic_F_third_term_bound} and the definition of $\LogisticXshXsInv$.
The second equality comes from \eqref{eq:logistic_lasso_weight_condition} and \eqref{eq:logistic_resid_condition}.
The final inequality comes from Assumption \eqref{eq:theo_ex_logistic_tune_param_ubd} and the form of $\LogisticTempResidConst$ which together imply $ \frac{\suppSizeTarg}{8} \LogisticEvalMax  \LogisticTempResidConst  \tuneparam \le \frac{1}{2\LogisticXshXsInv }$.



\ProofSection{4. The assumed conditions that ensure the contraction hold with high probability}

We complete the proof by showing \eqref{eq:logistic_resid_condition} holds with high probability.
Let  $\LogisticTempResidConst$ be given by \eqref{eq:logistic_resid_bound_nts_gap}.
Applying the union bound and Proposition 4a of \citep{fan2011nonconcave}
\begin{equation}
\begin{aligned}
 \prob{|| \frac{1}{n}X^T_{\suppSetTarg} (\LogisitcMuS{ \esttarg_{\suppSetTarg} }  - y) ||_{\text{max}} \ge  \frac{b_1 \tuneparam}{2 \maxccTarg} }
\le 2 \suppSizeTarg \exp\left(- \frac{n}{\LogisticMaxVarNorm} \frac{b_1^2 \tuneparam^2}{4 \maxccTargPow{2}} \right)
\end{aligned}
\end{equation}
and Claim \eqref{eq:theo_ex_logistic__reg_sat__prob_small_resid} follows.
Claim \eqref{eq:theo_ex_logistic__reg_sat__prob_nice_grad} follows from \eqref{eq:logistic_grad_at_targ_grad_holds_whp} and the same argument when $\LogisticTempResidConst$ is given by \eqref{eq:logistic_resid_bound_nts_grad}.

\end{proof}

\begin{proof} of theorem \ref{thm:logisitc_regession_lasso_pen_L1_err}
From the proof of Theorem 5 in \citep{fan2014strong}, under the event (1) conditioned on in the proof we have $||\estLasso - \esttarg||_1 \le 4 ||\estLassoSubscr{\supp} - \esttarg_{\suppSetTarg}||_1$.
In this case,
\begin{equation}
||\estLasso - \esttarg||_1 
 \le 4 ||\estLassoSubscr{\nnzSetTarg} - \esttarg_{\nnzSetTarg}||_1 
 \le 4 \supptargpow{1/2} ||\estLassoSubscr{\nnzSetTarg} - \esttarg_{\nnzSetTarg}||_2 
 \le 4 \supptargpow{1/2} ||\estLasso  - \esttarg||_2 
 \le \frac{20 \nnzSizeTarg}{\LogisticRestrEval}\ThreshTuneParam,
\end{equation}
where the last inequality comes from the conclusion of Theorem 5 in \citep{fan2014strong}.
\end{proof}

\section{Proofs for Appendices} \label{proof__appendices}

\subsection{Proofs for Appendix \ref{a:lap_spect_bounds}}

\begin{proof} of Proposition \ref{prop:lap_bounds}

\ProofSection{Entrywise norm comparisons}

For any entrywise $q$ norm, $q \in [1, \infty)$ we have
\begin{equation} \label{eq:lap_q_norm_deg_plus_edges}
\begin{aligned}
||\mathcal{L}(r)||_q^q  & = ||\text{diag}(\mathcal{A}(r) \mathbf{1}_d) -  \mathcal{A}(r)||_q^q  = ||\mathcal{A}(r) \mathbf{1}_d||_q^q  + ||\mathcal{A}(r)||_q^q \\
& = ||\mathcal{A}(r) \mathbf{1}_d||_q^q + 2 ||r||_q^q.
\end{aligned}
\end{equation}
Recall  (e.g. by H\"{o}lder's inequality) for any $q \ge 1$ and $x \in \mathbb{R}^d$, $||x||_1^q \le d^{q - 1} ||x||_q^q$. 
Therefore
\begin{equation} \label{eq:deg_Lq_ubd_adj_Lq}
\begin{aligned}
||\mathcal{A}(r) \mathbf{1}_d||_q^q  & = \sum_{i=1}^d \left| \sum_{j=1, j \neq i}^d \mathcal{A}(r)_{ij} \right|^q \\
& \le \sum_{i=1}^d \left| \sum_{j=1, j \neq i}^d |\mathcal{A}(r) _{ij}| \right|^q \\
& \le (d-1)^{q-1} \sum_{i=1}^d \sum_{j=1, j \neq i}^d |\mathcal{A}(r) _{ij}|^q\\
& =  2 (d-1)^{q-1} ||r||_q^q,
\end{aligned}
\end{equation}
thus
\begin{equation}\label{eq:lap_entry_Lq_from_Lq} 
||\mathcal{L}(r)||_q^q \le 2\left( (d-1)^{q-1} + 1 \right) ||r||_q^q,
\end{equation}
so  \eqref{eq:lap_bound_1__1}  and \eqref{eq:lap_bound_F__2} follow.

If each node has at most $\maxDegTarg - 1$ non-zero edges then \eqref{eq:lap_entry_Lq_from_Lq}  becomes
\begin{equation}\label{eq:lap_entry_Lq_from_Lq__max_deg} 
||\mathcal{L}(r)||_q^q \le 2\left( (\maxDegTarg - 1)^{q-1} + 1 \right) ||r||_q^q
\end{equation}
so the analogous claims about the maximal degree follow.


For \eqref{eq:lap_bound_F__max} note $||r||_2^2 \le D ||r||_{\text{max}}^2 $. 
Similarly, $|| \mathcal{A}(r) \mathbf{1}_d ||_{\text{max}} \le d  ||r||_{\text{max}}$. 
Also note $ ||\mathcal{A}(r) \mathbf{1}_d||_2^2 \le d (d  ||r||_{\text{max}})^2 $. 
Putting this together with \eqref{eq:lap_q_norm_deg_plus_edges} we get
$$
||\mathcal{L}(r)||_F^2 \le d^3  ||r||_{\text{max}}^2 +  2 {d \choose 2}  ||r||_{\text{max}}^2,
$$
and  \eqref{eq:lap_bound_F__max} follows. 

\ProofSection{Frobenius norm from Adjacency matrix operator norm}


Note $ ||\mathcal{A}(r) \mathbf{1}_d||_2 =  \sqrt{d}  ||\mathcal{A}(r)\frac{ \mathbf{1}_d}{\sqrt{d}}||_2 \le \sqrt{d}  ||\mathcal{A}(r)||_{\text{op}}$.
Furthermore $||A||_F \le d^{1/2} ||A||_{\text{op}}$ e.g. see Section 2.3.2 of \citep{golub2013matrix}.
Therefore starting with \eqref{eq:lap_q_norm_deg_plus_edges} we have
\begin{equation}
||\mathcal{L}(r)||^2_{F}  \le  ||\mathcal{A}(r) \mathbf{1}_d||_2^2 +  || \mathcal{A}(r)||_2^2 \le 2 d  ||A||_{\text{op}}^2
\end{equation}
and  \eqref{eq:lap_bound_frob__op} follows.

\ProofSection{Operator 2 from operator 1 and max norm}

Note $||\mathcal{L}(r)||_{\text{op}} \le \opone{\mathcal{L}(r)} $  e.g. by Corollary 2.3.2 of \citep{golub2013matrix} and the fact the Laplacian is symmetric. 
Furthermore,
\begin{equation}
\sum_{j=1}^d |\mathcal{L}(r)_{ij}|  
 = |\sum_{j=1, j \neq i} r_{ij}| + \sum_{j=1, j \neq i} |r_{ij}| 
 \le \sum_{j=1, j \neq i} |r_{ij}| + \sum_{j=1, j \neq i} |r_{ij}| 
= 2 \sum_{j=1, j \neq i} |\mathcal{A}(r)_{ij}|.
\end{equation}
Thus $ \opone{\mathcal{L}(r)} \le 2\opone{\mathcal{A}(r)} $ and \eqref{eq:lap_bound_op__op1} follows.

\ProofSection{Maximal degree upper bounds}

Using the triangle inequality we see
$$
\opn{\mathcal{L}(r))}  \le \opn{ \text{diag}(\mathcal{A}(r) \mathbf{1}_d)} + \opn{\mathcal{A}(r)} = ||\mathcal{A}(r) \mathbf{1}_d||_{\text{max}} + \opn{\mathcal{A}(r)} ,
$$
thus \eqref{eq:lap_bound_op__op_plus_deg} follows.

\ProofSection{Lower bound}
The lower bound claim follows from the fact that the largest singular value is at least the largest diagonal entry. 
\end{proof}

\begin{proof} of Proposition \ref{prop:lap_coef_bound_by_V}
WLOG $w$ is strictly positive.
Note $||\cdot||_{2, w}$ is a norm and $||\cdot||_{2} \le \sqrt{\max(w)} ||\cdot||_2$.
Starting with \eqref{eq:lap_coeff_formula} and applying the triangle inequality,
\begin{equation}
\begin{aligned}
\lapcoeff{V_A}{w}_{(ij)}^{1/2} & = ||V_A(i, :) - V_A(j, :)||_{2, w}  \\
& \le ||V_A(i, :) - V_B(i, :)||_{2, w}  + ||V_A(j, :) - V_B(j, :)||_{2, w}  + ||V_B(i, :) - V_B(j, :)||_{2, w} \\
& \le \sqrt{\max(w)} \left( ||V_A(i, :) - V_B(i, :)||_{2}  + ||V_A(j, :) - V_B(j, :)||_{2} \right) + \lapcoeff{V_B}{w}_{(ij)}^{1/2}
\end{aligned}
\end{equation}
Applying the same argument to $\lapcoeff{Y}{w}_{ij}^{1/2} $ we get
\begin{equation} \label{eq:lap_coef_bound_by_V_freddy}
\begin{aligned}
|\lapcoeff{V_A}{w}_{(ij)}^{1/2} - \lapcoeff{V_B}{w}_{(ij)}^{1/2} | 
& \le \sqrt{\max(w)}  \left( ||V_A(i, :) - V_B(i, :)||_{2}  + ||V_A(j, :) - V_B(j, :)||_{2}\right)\\
& \le 2 \sqrt{\max(w)} \opnorm{V_A- V_B}{2}{\infty}\\
%
%
\end{aligned}
\end{equation}
Thus the first claim follows.


Let $\cc{1}, \dots, \cc{k} \subseteq \mathcal{G}$ be the connected components of the graph, then
\begin{equation} 
\begin{aligned}
\sum_{(ij) \in \mathcal{G} } |\lapcoeff{V_A}{w}_{(ij)}^{1/2}  - \lapcoeff{V_B}{w}_{(ij)}^{1/2} |^2
& = \sum_{k=1}^K \sum_{(ij) \in \ccTarg{k}} |\lapcoeff{V_A}{w}_{(ij)}^{1/2}  - \lapcoeff{V_B}{w}_{(ij)}^{1/2} |^2 \\
& \le  \max(w) \sum_{k=1}^K \sum_{(ij) \in \ccTarg{k}}    \left( ||V_A(i, :) - V_B(i, :)||_{2}  + ||V_A(j, :) - V_B(j, :)||_{2}\right )^2  \\
& \le  \max(w) \sum_{k=1}^K \sum_{(ij) \in \ccTarg{k}}    \left( \sqrt{2} \sqrt{ ||V_A(i, :) - V_B(i, :)||_{2}^2  + ||V_A(j, :) - V_B(j, :)||_{2}^2 }\right )^2  \\
& =  2 \max(w) \sum_{k=1}^K \sum_{(ij) \in \ccTarg{k}}    \left(  ||V_A(i, :) - V_B(i, :)||_{2}^2  + ||V_A(j, :) - V_B(j, :)||_{2}^2 \right )  \\
& \le  2 \max(w) (\maxccTarg - 1)  \sum_{k=1}^K \sum_{i \in \ccTarg{k}}    ||V_A(i, :) - V_B(i, :)||_{2}^2   \\
& =  2 \max(w) (\maxccTarg - 1)   ||V_A - V_B||_{F}^2   \\
\end{aligned}
\end{equation}
where the first inequality comes from \eqref{eq:lap_coef_bound_by_V_freddy},
the second inequality comes from $||\cdot||_1 \le \sqrt{d} ||\cdot||_2$,
and
the third inequality comes from the fact that each node has at most $\maxccTarg - 1$ neighbors 
Thus \eqref{eq:lap_coef_bound_graph_sum_L1} follows.

For a matrix $R \in \mathbb{R}^{d \times d}$ note $||R||_F^4 = (||R||_F^2)^2 = (\sum_{i=1}^d ||R_i||_2^2)^2 \ge \sum_{i=1}^d ||R_i||_2^4$.
Applying the previous argument to $ |\lapcoeff{V_A}{w}_{(ij)}^{1/2}  - \lapcoeff{V_B}{w}_{(ij)}^{1/2} |^4$ (and using the inequality $||\cdot||_1 \le d^{3/4} ||\cdot||_2$) we obtain claim \eqref{eq:lap_coef_bound_graph_sum_L2}.


\end{proof}

\begin{proof} of Proposition \ref{prop:lap_coef_with_con_comps}

Recall from Proposition \ref{prop:spect_pen_maj} $\lapcoeff{V}{w}$ does not depend on the particular choice of basis, $V$, so we are free to pick any basis for the kernel.
Let $C_1, \dots, C_K \subseteq [d]$ be the indices of the connected components of $A$.
Let $\mathbf{1}_{C_k} \in \{0, 1\}^d$ be the indicator vector of the $k$th connected component and let $\widetilde{\mathbf{1}}_{C_k}  := \frac{1}{\sqrt{|C_k|}} \mathbf{1}_{C_k} $ be the normalized versions of these vectors.
Suppose $V =  [\widetilde{\mathbf{1}}_{C_1} | \dots | \widetilde{\mathbf{1}}_{C_K}] \in \mathbb{R}^{d \times K}$ e.g. by standard results about the graph Laplacian (Proposition 2 of \cite{von2007tutorial}) this is an orthonormal basis for the kernel.
Note $V(i, :)  =  \frac{1}{\sqrt{|C(i)|}} e_i $ for each $i$ where $C(i)$ is the connected component containing the $i$th vertex.
From this we can check that if $i, j$ are in the same connected component then $V(i, :) = V(j, :)$, thus $\lapcoeff{V}{w}_{ij} = 0$.
If $i, j$ are in different connected components we see
\begin{equation} \label{eq:lap_coef_lbd}
\begin{aligned}
\lapcoeff{V}{w}_{ij} 
& = ||  \frac{1}{\sqrt{|C(i)|}}e_i - \frac{1}{\sqrt{|C(j)|}} e_j ||_{w, 2}^2  
 \ge \min(w)  ||  \frac{1}{\sqrt{|C(i)|}}e_i - \frac{1}{\sqrt{|C(j)|}} e_j ||_{2}^2 \\
& = \min(w) \left(\frac{1}{|C(i)|} +  \frac{1}{|C(i)|}\right) 
 \ge  2 \frac{\min(w)}{\maxccTarg}
\end{aligned}
\end{equation}

\end{proof}

\begin{proof} of Lemma \ref{lem:lasso_maj_weight_bounds}

Let $V_y \in \mathbb{R}^{d \times K}$ be an orthonormal matrix whose columns are eigenvectors corresponding to the smallest $K$ eigenvalues of $\mathcal{L}(y)$. 
By a variant of the Davis-Kahan theorem (Theorem 2 of \cite{yu2015useful}) and the assumptions we see there exists an orthonormal matrix $\Theta \in \mathbb{R}^{K \times K}$ such that
\begin{equation} \label{eq:dk_consequence}
||V_x - V_y\Theta||_F 
%
%
\le 2^{3/2} \frac{\frobBound}{\triangle},
\end{equation}
thus
\begin{equation*} 
|\lapcoeff{V_x}{w}_{(ij)}^{1/2}  - \lapcoeff{V_y\Theta}{w}_{(ij)}^{1/2} |
 \le 2 \sqrt{\max(w)} \opnorm{V_x - V_y\Theta}{2}{\infty}
 \le 2 \sqrt{\max(w)} ||V_x - V_y\Theta||_F 
 \le 2^{5/2} \sqrt{\max(w)} \frac{\frobBound}{\triangle}
\end{equation*}
where the first inequality comes from Proposition \ref{prop:lap_coef_bound_by_V}.
The first two claims in \eqref{eq:lasso_maj_weight_bounds_from_frob} follow by applying Proposition  \ref{prop:lap_coef_with_con_comps} to $\lapcoeff{V_y\Theta}{w}_{(ij)}$ after after noting that $V_y\Theta$ is a basis for the kernel of $\mathcal{L}(y)$. 

Similarly we see
\begin{equation*} 
\begin{aligned}
\sum_{(ij) \in \mathcal{\suppSet}} |\lapcoeff{V_x}{w}_{(ij)}^{1/2}  |^2
& = \sum_{(ij) \in \mathcal{\suppSet}}  |\lapcoeff{V_x}{w}_{(ij)}^{1/2}  - \lapcoeff{V_y \Theta}{w}_{(ij)}^{1/2} |^2\\
& \le 2 \max(w) d_{\text{max}} ||V_x - V_y\Theta||_F^2\\
& \le 2 \max(w) 2^{3} \left( \frac{\frobBound}{\triangle} \right)^2
\end{aligned}
\end{equation*}
where the equality comes from Proposition  \ref{prop:lap_coef_with_con_comps},
the first inequality comes from \eqref{eq:lap_coef_bound_graph_sum_L1} of Proposition \ref{prop:lap_coef_bound_by_V},
and the final inequality comes from  \eqref{eq:dk_consequence}. 
Thus the $L_1$ claim  in \eqref{eq:lasso_maj_weight_bounds_from_frob} follows.

Once again we have
\begin{equation*}
\begin{aligned}
\sum_{(ij) \in \mathcal{\suppSet}} |\lapcoeff{V_x}{w}_{(ij)}^{1/2}  |^4
& = \sum_{(ij) \in \mathcal{\suppSet}}  |\lapcoeff{V_x}{w}_{(ij)}^{1/2}  - \lapcoeff{V_y \Theta}{w}_{(ij)}^{1/2} |^4\\
& \le 8 \max(w)^2 d_{\text{max}} ||V_x - V_y\Theta||_F^2\\
& \le 8 \max(w)^2 2^{6}  \left(\frac{\frobBound}{\triangle} \right)^4
\end{aligned}
\end{equation*}
where the equality comes from Proposition  \ref{prop:lap_coef_with_con_comps},
the first inequality comes from \eqref{eq:lap_coef_bound_graph_sum_L2} of Proposition \ref{prop:lap_coef_bound_by_V},
and the final inequality comes from  \eqref{eq:dk_consequence}. 
Thus $L_2$ claim  in \eqref{eq:lasso_maj_weight_bounds_from_frob} follows.

Finally assume $\opBound \le \frac{1}{2} \triangle$.
By another variant of the Davis-Kahan theorem, Lemma \ref{lem:dk_kernel}, there exists an orthonormal matrix $\Theta \in \mathbb{R}^{K \times K}$ such that
\begin{equation*} 
\opn{V_x - V_y\Theta} \le 2^{3/2} \frac{\opBound}{\triangle},
\end{equation*}
so applying  \ref{prop:lap_coef_bound_by_V} as above,
\begin{equation*} 
|\lapcoeff{V_x}{w}_{(ij)}^{1/2}  - \lapcoeff{V_y\Theta}{w}_{(ij)}^{1/2} |
 \le 2 \sqrt{\max(w)} \opnorm{V_x - V_y\Theta}{2}{\infty}
 \le 2 \sqrt{\max(w)} \opn{V_x - V_y\Theta} 
 \le 2^{5/2} \sqrt{\max(w)} \frac{\opBound}{\triangle}
\end{equation*}
where the first inequality comes from Proposition \ref{prop:lap_coef_bound_by_V} and we have used the fact that $\opnorm{\cdot}{2}{\infty} \le \opn{\cdot}$.
Thus \eqref{eq:lasso_maj_weight_bounds_from_opn} follows.

%


\end{proof}

\begin{proof} of Lemma \ref{lem:dk_kernel}
We will make frequent reference to standard matrix analysis results stated in Section 2 of \citep{chen2020spectral}.
Let $\arbNorm{\cdot}$ be an arbitrary, unitarily invariant matrix norm; see Definition 2.6.1 of \cite{chen2020spectral}.

Note 
$$
U_{B, \perp} (B - A) U_A = \text{diag}(\lambda_{B, \perp})U_{B, \perp}^T U_A
$$
since $U_A$ is a basis for the kernel of $A$ so
$$
\arbNorm{U_{B, \perp} (B - A) U_A}
 = \arbNorm{ \text{diag}(\lambda_{B, \perp})U_{B, \perp}^T U_A}
\ge  \min(\lambda_{B, \perp})  \arbNorm{U_{B, \perp}^T U_A}
= \triangle_B \arbNorm{U_{B, \perp}^T U_A}
$$
by Lemma 2.6.1 of \citep{chen2020spectral}.
Therefore
\begin{equation} \label{eq:fuzzy_cat}
 \arbNorm{U_{B, \perp}^T U_A} \le \frac{\arbNorm{U_{B, \perp} (B - A) U_A}}{\triangle_B} 
 \le \frac{\arbNorm{(B - A) U_A}}{\triangle_B} 
 = \frac{\arbNorm{B  U_A}}{\triangle_B} 
 \le \frac{\arbNorm{B - A}}{\triangle_B} 
\end{equation}
where the second and final inequalities are from the fact $\arbNorm{}$ is unitarily invariant, Lemma 2.6.1 of \citep{chen2020spectral} and the facts $\opn{U_{B, \perp}} = \opn{U_A} = 1$.
The equality comes from the fact $A U_A = 0$.

Under the assumption $\opn{B - A} \le \frac{1}{2} \triangle_A$
we have $\triangle_B \ge \triangle_A - \opn{B - A} \ge \frac{1}{2}  \triangle_A$ by Weyl's inequality.
Thus in this case we can rewrite \eqref{eq:fuzzy_cat} as 
\begin{equation} \label{eq:fuzzy_cat_A_gap}
 \arbNorm{U_{B, \perp}^T U_A} \le \frac{2 \arbNorm{U_{B, \perp} (B - A) U_A}}{\triangle_A} 
 \le \frac{2 \arbNorm{(B - A) U_A}}{\triangle_A} 
 = \frac{2\arbNorm{B  U_A}}{\triangle_A} 
 \le \frac{2\arbNorm{B - A}}{\triangle_A} 
\end{equation}

By Lemma 2.1.3 of \citep{chen2020spectral} and the fact that the set of orthogonal matrices is compact (thus the minimizer in this lemma is attained by some matrix $Q$), we see there exists an orthonormal $Q \in \mathbb{R}^{K \times K}$ such that
$$
\opn{U_{B, \perp} - U_A Q} \le \sqrt{2} \opn{U_{B, \perp}U_{B, \perp}^T - U_A U_A^T} = \sqrt{2} \opn{U_{B, \perp}^T U_A } 
$$
where the equality comes from Lemma 2.1.2 of  \citep{chen2020spectral}.
The two claims now follow.

\end{proof}

\subsection{Proofs for Appendix \ref{a:bd_rect_multi_array} }  \label{a:proofs__bd_rect_multi_array}

\begin{proof} of Proposition \ref{prop:hypergraph_to_adj_mat_conn_compts}
Suppose $\{ (v_i, k\spa{v_i}_i)\}_{i=1}^M$ is a path in $\HyperGraph(A)$ between nodes $(v_1, k\spa{v_1}_1)$ and $(v_M, k\spa{v_M}_M)$.
This means for each pair consecutive pair $(v_i, k\spa{v_i}_i)$ and $(v_{i + 1}, k\spa{v_{i + 1}}_{i + 1})$ there is an entry $A_{j\spa{1}, \dots, j\spa{V}} \neq 0$ where $j\spa{v_i} =  k\spa{v_i}_i$ and $j\spa{v_{i + 1}} =  k\spa{v_{i + 1}}_{i + 1}$.
By inspecting the right hand side of \eqref{eq:hypergraph_adj_mat_edge_sum}, we see that $\HyperAdjMat(A)_{(v_i, k\spa{v_i}_i), (v_{i + 1}, k\spa{v_{i + 1}}_{i + 1})} \neq 0$ (recall each entry of $A$ is nonnegative).
Thus $\{ (v_i, k\spa{v_i}_i)\}_{i=1}^M$  is a path in $\HyperAdjMat(A)$.

This argument can be reversed to show that if $\{ (v_i, k\spa{v_i}_i)\}_{i=1}^M$ is a path in $\HyperAdjMat(A)$ then it is also a path in $\HyperGraph(A)$.
\end{proof}

\subsection{Proofs for Appendix \ref{a:tune_param_ubd} }  \label{a:proofs__tune_param_ubd}

\begin{proof} of Proposition \ref{prop:fcls_tune_param_max_val}
Let $\lambda \in \mathbb{R}^d_+$ be the eigenvalues of $\mathcal{L}(|\estinit|) $ and let $V \in \mathbb{R}^{d \times d}$ be a matrix whose columns are a corresponding set of orthonormal eigenvectors.
Also let $w \in \mathbb{R}^d_+$ be given by $w_j = g_{\tuneparam}(\lambda_j)$. 
Fix $\tuneparam \ge \FCLSKillerBound$.

By construction $\tuneparam b_1 \ge  \lambda_{\text{max}} \left(  \mathcal{L}(|\estinit|) \right)$ thus by Definition \ref{def:scad_like_pen_func} 
$$w_j =  g_{\tuneparam}'\left(\lambda_{j} \left(  \mathcal{L}(|\estinit|) \right) \right) \ge a_1 \tuneparam$$
 for each $j \in [D]$.
Recalling \eqref{eq:lap_coeff_formula} we have
 \begin{equation} 
 \begin{aligned}
\lapcoeff{V}{w}_{(ij)}  &=   ||V(i, :) - V(j, :)||_{2, w}^2 \ge \min(w)  ||V(i, :) - V(j, :)||_{2}^2  \ge  a_1 \tuneparam ||V(i, :) - V(j, :)||_{2}^2  = 2 a_1 \tuneparam 
\end{aligned}
\end{equation}
where the last equality comes from the orthonormality of $V$.
Recalling \eqref{prob:weighted_lasso_maj} we see $\estone$ is a solution to the following weighted Lasso problem \eqref{prob:lasso_problem_pen} where each $c_j \ge a_1$. 
By \eqref{eq:fcls_tune_param_max_val}, $\tuneparam a_1 \ge \LassoKillerBound$ so by the definition of $\LassoKillerBound$ we have $\estone = 0$.

Since $\estone = 0$, $\lambda_{\text{max}} (\mathcal{L}(|\estone|)) = 0 \le \lambda_{\text{max}} (\mathcal{L}(|\estinit|)) $ so the same argument shows $\esttwo = 0$.
 

\end{proof}


\bibliographystyle{apalike}
\bibliography{refs}

\begin{thebibliography}{}

\bibitem[Ando and Zhang, 2007]{ando2007learning}
Ando, R.~K. and Zhang, T. (2007).
\newblock Learning on graph with laplacian regularization.
\newblock {\em Advances in neural information processing systems}, 19:25.

\bibitem[Asteris et~al., 2015]{asteris2015sparse}
Asteris, M., Papailiopoulos, D., Kyrillidis, A., and Dimakis, A.~G. (2015).
\newblock Sparse pca via bipartite matchings.
\newblock In {\em Advances in Neural Information Processing Systems}, pages
  766--774.

\bibitem[Beck, 2017]{beck2017first}
Beck, A. (2017).
\newblock {\em First-order methods in optimization}.
\newblock SIAM.

\bibitem[Berge, 1984]{berge1984hypergraphs}
Berge, C. (1984).
\newblock {\em Hypergraphs: combinatorics of finite sets}, volume~45.
\newblock Elsevier.

\bibitem[Bertrand and Massias, 2021]{bertrand2021anderson}
Bertrand, Q. and Massias, M. (2021).
\newblock Anderson acceleration of coordinate descent.
\newblock In {\em International Conference on Artificial Intelligence and
  Statistics}, pages 1288--1296. PMLR.

\bibitem[Bickel et~al., 2008]{bickel2008covariance}
Bickel, P.~J., Levina, E., et~al. (2008).
\newblock Covariance regularization by thresholding.
\newblock {\em The Annals of Statistics}, 36(6):2577--2604.

\bibitem[Broto et~al., 2019]{broto2019block}
Broto, B., Bachoc, F., Clouvel, L., and Martinez, J.-M. (2019).
\newblock Block-diagonal covariance estimation and application to the shapley
  effects in sensitivity analysis.
\newblock {\em arXiv preprint arXiv:1907.12780}.

\bibitem[Cand{\`e}s et~al., 2011]{candes2011robust}
Cand{\`e}s, E.~J., Li, X., Ma, Y., and Wright, J. (2011).
\newblock Robust principal component analysis?
\newblock {\em Journal of the ACM (JACM)}, 58(3):1--37.

\bibitem[Candes et~al., 2008]{candes2008enhancing}
Candes, E.~J., Wakin, M.~B., and Boyd, S.~P. (2008).
\newblock Enhancing sparsity by reweighted l1 minimization.
\newblock {\em Journal of Fourier analysis and applications}, 14(5-6):877--905.

\bibitem[Carmichael, 2020]{carmichael2020learning}
Carmichael, I. (2020).
\newblock Learning sparsity and block diagonal structure in multi-view mixture
  models.
\newblock {\em arXiv preprint arXiv:2012.15313}.

\bibitem[Carmichael, 2021]{carmichael2021yet}
Carmichael, I. (2021).
\newblock {idc9/ya\_glm} yet another penalized generalized linear model package
  in python.
\newblock {\em \url{https://doi.org/10.5281/zenodo.5076358}}.

\bibitem[Chen et~al., 2020]{chen2020spectral}
Chen, Y., Chi, Y., Fan, J., and Ma, C. (2020).
\newblock Spectral methods for data science: A statistical perspective.
\newblock {\em arXiv preprint arXiv:2012.08496}.

\bibitem[De~Abreu, 2007]{de2007old}
De~Abreu, N. M.~M. (2007).
\newblock Old and new results on algebraic connectivity of graphs.
\newblock {\em Linear algebra and its applications}, 423(1):53--73.

\bibitem[Devijver and Gallopin, 2018]{devijver2018block}
Devijver, E. and Gallopin, M. (2018).
\newblock Block-diagonal covariance selection for high-dimensional gaussian
  graphical models.
\newblock {\em Journal of the American Statistical Association},
  113(521):306--314.

\bibitem[Dewaskar et~al., 2020]{dewaskar2020finding}
Dewaskar, M., Palowitch, J., He, M., Love, M.~I., and Nobel, A. (2020).
\newblock Finding stable groups of cross-correlated features in multi-view
  data.
\newblock {\em arXiv preprint arXiv:2009.05079}.

\bibitem[Donoho and Johnstone, 1994]{donoho1994ideal}
Donoho, D.~L. and Johnstone, J.~M. (1994).
\newblock Ideal spatial adaptation by wavelet shrinkage.
\newblock {\em biometrika}, 81(3):425--455.

\bibitem[Egilmez et~al., 2017]{egilmez2017graph}
Egilmez, H.~E., Pavez, E., and Ortega, A. (2017).
\newblock Graph learning from data under laplacian and structural constraints.
\newblock {\em IEEE Journal of Selected Topics in Signal Processing},
  11(6):825--841.

\bibitem[Fan and Li, 2001]{fan2001variable}
Fan, J. and Li, R. (2001).
\newblock Variable selection via nonconcave penalized likelihood and its oracle
  properties.
\newblock {\em Journal of the American statistical Association},
  96(456):1348--1360.

\bibitem[Fan et~al., 2020]{fan2020statistical}
Fan, J., Li, R., Zhang, C.-H., and Zou, H. (2020).
\newblock {\em Statistical foundations of data science}.
\newblock CRC press.

\bibitem[Fan and Lv, 2011]{fan2011nonconcave}
Fan, J. and Lv, J. (2011).
\newblock Nonconcave penalized likelihood with np-dimensionality.
\newblock {\em IEEE Transactions on Information Theory}, 57(8):5467--5484.

\bibitem[Fan et~al., 2014]{fan2014strong}
Fan, J., Xue, L., and Zou, H. (2014).
\newblock Strong oracle optimality of folded concave penalized estimation.
\newblock {\em Annals of statistics}, 42(3):819.

\bibitem[Feng et~al., 2014]{feng2014robust}
Feng, J., Lin, Z., Xu, H., and Yan, S. (2014).
\newblock Robust subspace segmentation with block-diagonal prior.
\newblock In {\em Proceedings of the IEEE conference on computer vision and
  pattern recognition}, pages 3818--3825.

\bibitem[Fiedler, 1973]{fiedler1973algebraic}
Fiedler, M. (1973).
\newblock Algebraic connectivity of graphs.
\newblock {\em Czechoslovak mathematical journal}, 23(2):298--305.

\bibitem[Friedman et~al., 2010]{friedman2010regularization}
Friedman, J., Hastie, T., and Tibshirani, R. (2010).
\newblock Regularization paths for generalized linear models via coordinate
  descent.
\newblock {\em Journal of statistical software}, 33(1):1.

\bibitem[Gasso et~al., 2009]{gasso2009recovering}
Gasso, G., Rakotomamonjy, A., and Canu, S. (2009).
\newblock Recovering sparse signals with a certain family of nonconvex
  penalties and dc programming.
\newblock {\em IEEE Transactions on Signal Processing}, 57(12):4686--4698.

\bibitem[Golub and Van~Loan, 2013]{golub2013matrix}
Golub, G. and Van~Loan, C. (2013).
\newblock Matrix computations 4th edition the johns hopkins university press.
\newblock {\em Baltimore, MD}.

\bibitem[Han et~al., 2017]{han2017bilateral}
Han, J., Song, K., Nie, F., and Li, X. (2017).
\newblock Bilateral k-means algorithm for fast co-clustering.
\newblock In {\em Proceedings of the Thirty-First AAAI Conference on Artificial
  Intelligence}, pages 1969--1975.

\bibitem[Harris et~al., 2020]{harris2020array}
Harris, C.~R., Millman, K.~J., van~der Walt, S.~J., Gommers, R., Virtanen, P.,
  Cournapeau, D., Wieser, E., Taylor, J., Berg, S., Smith, N.~J., Kern, R.,
  Picus, M., Hoyer, S., van Kerkwijk, M.~H., Brett, M., Haldane, A., del
  R{\'{i}}o, J.~F., Wiebe, M., Peterson, P., G{\'{e}}rard-Marchant, P.,
  Sheppard, K., Reddy, T., Weckesser, W., Abbasi, H., Gohlke, C., and Oliphant,
  T.~E. (2020).
\newblock Array programming with {NumPy}.
\newblock {\em Nature}, 585(7825):357--362.

\bibitem[Hastie et~al., 2015]{hastie2015statistical}
Hastie, T., Tibshirani, R., and Wainwright, M. (2015).
\newblock {\em Statistical learning with sparsity: the lasso and
  generalizations}.
\newblock CRC press.

\bibitem[Hunter, 2007]{matplotlib2007hunter}
Hunter, J.~D. (2007).
\newblock Matplotlib: A 2d graphics environment.
\newblock {\em Computing in Science \& Engineering}, 9(3):90--95.

\bibitem[Hyodo et~al., 2015]{hyodo2015testing}
Hyodo, M., Shutoh, N., Nishiyama, T., and Pavlenko, T. (2015).
\newblock Testing block-diagonal covariance structure for high-dimensional
  data.
\newblock {\em Statistica Neerlandica}, 69(4):460--482.

\bibitem[Jiang et~al., 2013]{jiang2013graph}
Jiang, B., Ding, C., and Tang, J. (2013).
\newblock Graph-laplacian pca: Closed-form solution and robustness.
\newblock In {\em Proceedings of the IEEE Conference on Computer Vision and
  Pattern Recognition}, pages 3492--3498.

\bibitem[Krishnapuram et~al., 2005]{krishnapuram2005sparse}
Krishnapuram, B., Carin, L., Figueiredo, M.~A., and Hartemink, A.~J. (2005).
\newblock Sparse multinomial logistic regression: Fast algorithms and
  generalization bounds.
\newblock {\em IEEE transactions on pattern analysis and machine intelligence},
  27(6):957--968.

\bibitem[Kumar et~al., 2019]{kumar2019unified}
Kumar, S., Ying, J., Cardoso, J. V. d.~M., and Palomar, D. (2019).
\newblock A unified framework for structured graph learning via spectral
  constraints.
\newblock {\em arXiv preprint arXiv:1904.09792}.

\bibitem[Lange et~al., 2000]{lange2000optimization}
Lange, K., Hunter, D.~R., and Yang, I. (2000).
\newblock Optimization transfer using surrogate objective functions.
\newblock {\em Journal of computational and graphical statistics}, 9(1):1--20.

\bibitem[Levin et~al., 2019]{levin2019recovering}
Levin, K., Lodhia, A., and Levina, E. (2019).
\newblock Recovering low-rank structure from multiple networks with unknown
  edge distributions.
\newblock {\em arXiv preprint arXiv:1906.07265}.

\bibitem[Lewis, 1999]{lewis1999nonsmooth}
Lewis, A.~S. (1999).
\newblock Nonsmooth analysis of eigenvalues.
\newblock {\em Mathematical Programming}, 84(1):1--24.

\bibitem[Li et~al., 2020]{li2020graph}
Li, Y., Mark, B., Raskutti, G., Willett, R., Song, H., and Neiman, D. (2020).
\newblock Graph-based regularization for regression problems with alignment and
  highly correlated designs.
\newblock {\em SIAM journal on mathematics of data science}, 2(2):480--504.

\bibitem[Liu et~al., 2012]{liu2012high}
Liu, H., Han, F., Yuan, M., Lafferty, J., Wasserman, L., et~al. (2012).
\newblock High-dimensional semiparametric gaussian copula graphical models.
\newblock {\em The Annals of Statistics}, 40(4):2293--2326.

\bibitem[Liu et~al., 2017]{liu2017folded}
Liu, H., Yao, T., Li, R., and Ye, Y. (2017).
\newblock Folded concave penalized sparse linear regression: sparsity,
  statistical performance, and algorithmic theory for local solutions.
\newblock {\em Mathematical programming}, 166(1):207--240.

\bibitem[Loh and Wainwright, 2015]{loh2015regularized}
Loh, P.-L. and Wainwright, M.~J. (2015).
\newblock Regularized m-estimators with nonconvexity: Statistical and
  algorithmic theory for local optima.
\newblock {\em The Journal of Machine Learning Research}, 16(1):559--616.

\bibitem[Lu et~al., 2018]{lu2018subspace}
Lu, C., Feng, J., Lin, Z., Mei, T., and Yan, S. (2018).
\newblock Subspace clustering by block diagonal representation.
\newblock {\em IEEE transactions on pattern analysis and machine intelligence},
  41(2):487--501.

\bibitem[Lu et~al., 2007]{lu2007lower}
Lu, M., Zhang, L.-z., and Tian, F. (2007).
\newblock Lower bounds of the laplacian spectrum of graphs based on diameter.
\newblock {\em Linear algebra and its applications}, 420(2-3):400--406.

\bibitem[Mansinghka et~al., 2006]{mansinghka2006structured}
Mansinghka, V., Kemp, C., Tenenbaum, J., and Griffiths, T. (2006).
\newblock Structured priors for structure learning.
\newblock In {\em Proceedings of the Twenty-Second Conference on Uncertainty in
  Artificial Intelligence}, pages 324--331.

\bibitem[Marlin and Murphy, 2009]{marlin2009sparse}
Marlin, B.~M. and Murphy, K.~P. (2009).
\newblock Sparse gaussian graphical models with unknown block structure.
\newblock In {\em Proceedings of the 26th Annual International Conference on
  Machine Learning}, pages 705--712.

\bibitem[Massias et~al., 2020]{massias2020dual}
Massias, M., Vaiter, S., Gramfort, A., and Salmon, J. (2020).
\newblock Dual extrapolation for sparse generalized linear models.
\newblock {\em Journal of Machine Learning Research}, 21(234):1--33.

\bibitem[Negahban et~al., 2012]{negahban2012unified}
Negahban, S.~N., Ravikumar, P., Wainwright, M.~J., Yu, B., et~al. (2012).
\newblock A unified framework for high-dimensional analysis of $ m $-estimators
  with decomposable regularizers.
\newblock {\em Statistical science}, 27(4):538--557.

\bibitem[Nie et~al., 2017]{nie2017learning}
Nie, F., Wang, X., Deng, C., and Huang, H. (2017).
\newblock Learning a structured optimal bipartite graph for co-clustering.
\newblock In {\em Advances in Neural Information Processing Systems}, pages
  4129--4138.

\bibitem[Nie et~al., 2016]{nie2016constrained}
Nie, F., Wang, X., Jordan, M.~I., and Huang, H. (2016).
\newblock The constrained laplacian rank algorithm for graph-based clustering.
\newblock In {\em Thirtieth AAAI Conference on Artificial Intelligence}.

\bibitem[Ortega and Rheinboldt, 2000]{ortega2000iterative}
Ortega, J.~M. and Rheinboldt, W.~C. (2000).
\newblock {\em Iterative solution of nonlinear equations in several variables}.
\newblock SIAM.

\bibitem[Pavlenko et~al., 2012]{pavlenko2012covariance}
Pavlenko, T., Bj{\"o}rkstr{\"o}m, A., and Tillander, A. (2012).
\newblock Covariance structure approximation via glasso in high-dimensional
  supervised classification.
\newblock {\em Journal of Applied Statistics}, 39(8):1643--1666.

\bibitem[Pedregosa et~al., 2011]{scikit2011pedrogosa}
Pedregosa, F., Varoquaux, G., Gramfort, A., Michel, V., Thirion, B., Grisel,
  O., Blondel, M., Prettenhofer, P., Weiss, R., Dubourg, V., Vanderplas, J.,
  Passos, A., Cournapeau, D., Brucher, M., Perrot, M., and Duchesnay, E.
  (2011).
\newblock Scikit-learn: Machine learning in {P}ython.
\newblock {\em Journal of Machine Learning Research}, 12:2825--2830.

\bibitem[Petersen and Pedersen, 2012]{petersen2012matrix}
Petersen, K.~B. and Pedersen, M.~S. (2012).
\newblock The matrix cookbook (version: November 15, 2012).

\bibitem[Porter et~al., 2009]{porter2009communities}
Porter, M.~A., Onnela, J.-P., and Mucha, P.~J. (2009).
\newblock Communities in networks.
\newblock {\em Notices of the AMS}, 56(9):1082--1097.

\bibitem[Rad et~al., 2011]{rad2011lower}
Rad, A.~A., Jalili, M., and Hasler, M. (2011).
\newblock A lower bound for algebraic connectivity based on the
  connection-graph-stability method.
\newblock {\em Linear algebra and its applications}, 435(1):186--192.

\bibitem[Ravikumar et~al., 2011]{ravikumar2011high}
Ravikumar, P., Wainwright, M.~J., Raskutti, G., Yu, B., et~al. (2011).
\newblock High-dimensional covariance estimation by minimizing l1-penalized
  log-determinant divergence.
\newblock {\em Electronic Journal of Statistics}, 5:935--980.

\bibitem[Rothman et~al., 2009]{rothman2009generalized}
Rothman, A.~J., Levina, E., and Zhu, J. (2009).
\newblock Generalized thresholding of large covariance matrices.
\newblock {\em Journal of the American Statistical Association},
  104(485):177--186.

\bibitem[Sharpnack et~al., 2012]{sharpnack2012sparsistency}
Sharpnack, J., Singh, A., and Rinaldo, A. (2012).
\newblock Sparsistency of the edge lasso over graphs.
\newblock In {\em Artificial Intelligence and Statistics}, pages 1028--1036.
  PMLR.

\bibitem[Smola and Kondor, 2003]{smola2003kernels}
Smola, A.~J. and Kondor, R. (2003).
\newblock Kernels and regularization on graphs.
\newblock In {\em Learning theory and kernel machines}, pages 144--158.
  Springer.

\bibitem[Spielman, 2004]{spielman2004spectral}
Spielman, D. (2004).
\newblock Lecture notes for spectral graph theory and its applications.
\newblock Available at
  \url{https://ocw.mit.edu/courses/mathematics/18-409-topics-in-theoretical-computer-science-an-algorithmists-toolkit-fall-2009/readings/MIT18_409F09_spiel_lec2.pdf}.

\bibitem[Sun et~al., 2015]{sun2015inferring}
Sun, S., Wang, H., and Xu, J. (2015).
\newblock Inferring block structure of graphical models in exponential
  families.
\newblock In {\em Artificial Intelligence and Statistics}, pages 939--947.
  PMLR.

\bibitem[Sun et~al., 2016]{sun2016majorization}
Sun, Y., Babu, P., and Palomar, D.~P. (2016).
\newblock Majorization-minimization algorithms in signal processing,
  communications, and machine learning.
\newblock {\em IEEE Transactions on Signal Processing}, 65(3):794--816.

\bibitem[Tam and Dunson, 2020]{tam2020fiedler}
Tam, E. and Dunson, D. (2020).
\newblock Fiedler regularization: Learning neural networks with graph sparsity.
\newblock In {\em International Conference on Machine Learning}, pages
  9346--9355. PMLR.

\bibitem[Tan et~al., 2014]{tan2014learning}
Tan, K.~M., London, P., Mohan, K., Lee, S.-I., Fazel, M., and Witten, D.
  (2014).
\newblock Learning graphical models with hubs.
\newblock {\em arXiv preprint arXiv:1402.7349}.

\bibitem[Tan et~al., 2015]{tan2015cluster}
Tan, K.~M., Witten, D., and Shojaie, A. (2015).
\newblock The cluster graphical lasso for improved estimation of gaussian
  graphical models.
\newblock {\em Computational statistics \& data analysis}, 85:23--36.

\bibitem[Tibshirani et~al., 2005]{tibshirani2005sparsity}
Tibshirani, R., Saunders, M., Rosset, S., Zhu, J., and Knight, K. (2005).
\newblock Sparsity and smoothness via the fused lasso.
\newblock {\em Journal of the Royal Statistical Society: Series B (Statistical
  Methodology)}, 67(1):91--108.

\bibitem[Tibshirani et~al., 2011]{tibshirani2011solution}
Tibshirani, R.~J., Taylor, J., et~al. (2011).
\newblock The solution path of the generalized lasso.
\newblock {\em The annals of statistics}, 39(3):1335--1371.

\bibitem[Vershynin, 2018]{vershynin2018high}
Vershynin, R. (2018).
\newblock {\em High-dimensional probability: An introduction with applications
  in data science}, volume~47.
\newblock Cambridge university press.

\bibitem[Virtanen et~al., 2020]{scipy2020virtanen}
Virtanen, P., Gommers, R., Oliphant, T.~E., Haberland, M., Reddy, T.,
  Cournapeau, D., Burovski, E., Peterson, P., Weckesser, W., Bright, J., {van
  der Walt}, S.~J., Brett, M., Wilson, J., Millman, K.~J., Mayorov, N., Nelson,
  A. R.~J., Jones, E., Kern, R., Larson, E., Carey, C.~J., Polat, {\.I}., Feng,
  Y., Moore, E.~W., {VanderPlas}, J., Laxalde, D., Perktold, J., Cimrman, R.,
  Henriksen, I., Quintero, E.~A., Harris, C.~R., Archibald, A.~M., Ribeiro,
  A.~H., Pedregosa, F., {van Mulbregt}, P., and {SciPy 1.0 Contributors}
  (2020).
\newblock {{SciPy} 1.0: Fundamental Algorithms for Scientific Computing in
  Python}.
\newblock {\em Nature Methods}, 17:261--272.

\bibitem[Von~Luxburg, 2007]{von2007tutorial}
Von~Luxburg, U. (2007).
\newblock A tutorial on spectral clustering.
\newblock {\em Statistics and computing}, 17(4):395--416.

\bibitem[Wainwright, 2019]{wainwright2019high}
Wainwright, M.~J. (2019).
\newblock {\em High-dimensional statistics: A non-asymptotic viewpoint},
  volume~48.
\newblock Cambridge University Press.

\bibitem[Waskom, 2021]{seaborn2021waskom}
Waskom, M.~L. (2021).
\newblock seaborn: statistical data visualization.
\newblock {\em Journal of Open Source Software}, 6(60):3021.

\bibitem[Yu et~al., 2015]{yu2015useful}
Yu, Y., Wang, T., and Samworth, R.~J. (2015).
\newblock A useful variant of the davis--kahan theorem for statisticians.
\newblock {\em Biometrika}, 102(2):315--323.

\bibitem[Yuan and Lin, 2006]{yuan2006model}
Yuan, M. and Lin, Y. (2006).
\newblock Model selection and estimation in regression with grouped variables.
\newblock {\em Journal of the Royal Statistical Society: Series B (Statistical
  Methodology)}, 68(1):49--67.

\bibitem[Zhang et~al., 2012]{zhang2012general}
Zhang, C.-H., Zhang, T., et~al. (2012).
\newblock A general theory of concave regularization for high-dimensional
  sparse estimation problems.
\newblock {\em Statistical Science}, 27(4):576--593.

\bibitem[Zhou et~al., 2006]{zhou2006learning}
Zhou, D., Huang, J., and Sch{\"o}lkopf, B. (2006).
\newblock Learning with hypergraphs: Clustering, classification, and embedding.
\newblock {\em Advances in neural information processing systems},
  19:1601--1608.

\bibitem[Zou and Li, 2008]{zou2008one}
Zou, H. and Li, R. (2008).
\newblock One-step sparse estimates in nonconcave penalized likelihood models.
\newblock {\em Annals of statistics}, 36(4):1509.

\end{thebibliography}

\end{document}